\documentclass[11pt,reqno]{amsart}

\usepackage{amssymb, amscd, stmaryrd}
\usepackage{enumerate}
\usepackage{cases}
\usepackage{amsfonts}
\usepackage[dvips,dvipdf]{graphicx}
\usepackage[usenames,dvipsnames]{color}

\usepackage{color}
\usepackage{color, colortbl}
\usepackage{multirow}
\usepackage[table]{xcolor}
\usepackage{comment}
\usepackage{subfigure}
\usepackage{epstopdf}
\usepackage{tikz}
\usepackage{chngcntr}
\usepackage{url}
\usepackage{hyperref}
\usepackage{cleveref}

\usepackage{geometry}
\hypersetup{
  bookmarksnumbered = true,
  bookmarksopen=false,
  pdfborder=0 0 0,         
  pdffitwindow=true,      
  pdfnewwindow=true, 
  colorlinks=true,           
  linkcolor=blue,            
  citecolor=magenta,    
  filecolor=magenta,     
  urlcolor=cyan              
}

\makeatletter
\newcommand\figcaption{\def\@captype{figure}\caption}
\newcommand\tabcaption{\def\@captype{table}\caption}
\makeatother

\counterwithin*{equation}{section}

\newcommand{\pa}{\partial}

\newcommand{\maT}{\mathcal T}

\newcommand{\be}{\begin{eqnarray}}
\newcommand{\ee}{\end{eqnarray}}
\newcommand{\ben}{\begin{eqnarray*}}
\newcommand{\een}{\end{eqnarray*}}

\def\be{\begin{equation}}
\def\ee{\end{equation}}
\def\bes{\begin{equation*}}
\def\ees{\end{equation*}}
\def\beq{\begin{eqnarray}}
\def\eeq{\end{eqnarray}}
\def\beqs{\begin{eqnarray*}}
\def\eeqs{\end{eqnarray*}}
\def\bal{\begin{aligned}}
\def\eal{\end{aligned}}
\def\bsqs{\begin{subequations}}
\def\esqs{\end{subequations}}

\newtheorem{theorem}{Theorem}[section]

\newtheorem{corollary}[theorem]{Corollary}

\newtheorem{lemma}[theorem]{Lemma}
\newtheorem{lem}[theorem]{Lemma}

\theoremstyle{definition}
\newtheorem{definition}[theorem]{Definition}
\theoremstyle{definition}
\newtheorem{algorithm}[theorem]{Algorithm}
\theoremstyle{remark}
\newtheorem{remark}[theorem]{Remark}
\newtheorem{example}[theorem]{Example}

\makeatletter
\newcommand{\definetitlefootnote}[1]{%
  \newcommand\addtitlefootnote{%
    \makebox[0pt][l]{$^{\bigstar}$}%
    \footnote{\protect\@titlefootnotetext}
  }%
  \newcommand\@titlefootnotetext{\spaceskip=\z@skip $^{\bigstar}$#1}%
}
\makeatother

\author[H. Li, P. Yin]{Hengguang Li$^\dagger$, Peimeng Yin$^{\ddagger,*}$}
\address{$^\dagger$ Department of Mathematics, Wayne State University, Detroit, MI 48202, USA}
\email{li@wayne.edu}
\address{$^\ddagger$ Department of Mathematical Sciences, The University of Texas at El Paso, El Paso, TX 79968, USA.}
\email{pyin@utep.edu}
\thanks{* Corresponding author}

\begin{document}

\title[FEM for sixth order equation]{A $C^0$ finite element algorithm for the sixth order problem with simply supported boundary conditions}

\date{\today}

\begin{abstract}
In this paper, we investigate a sixth order elliptic equation with the simply supported boundary conditions in a polygonal domain. We propose a new method that decouples the sixth order problem into a system of second order equations. Unlike the direct decomposition, which yields three Poisson problems but is restricted to polygonal domains with the largest interior angle no more than ${\pi}/{2}$, we rigorously analyze and construct extra Poisson problems to confine the solution into the same function space as that of the original sixth order problem. Consequently, the proposed method can be applied to general polygonal domains. In turn, we also present a $C^0$ finite element algorithm to discretize the new resulting system and establish optimal error estimates for the numerical solution on quasi-uniform meshes. Finally, numerical experiments are performed to validate the theoretical findings.
\end{abstract}

\maketitle


\section{Introduction}
Consider the sixth order elliptic problem, also known as the triharmonic problem
\begin{eqnarray}\label{eqnbi}
-\Delta^3u=f \quad {\rm{in}} \ \Omega,\qquad \quad u=\Delta u=\Delta^2 u=0 \quad {\rm{on}} \ \pa\Omega,
\end{eqnarray}
where $\Omega\subset \mathbb R^2$ is a polygonal domain and $f \in H^{-1}(\Omega)$. The boundary conditions in (\ref{eqnbi}) are commonly referred to as simply supported boundary conditions \cite{droniou2019}. 
The sixth order partial differential equations (PDEs) arise from various mathematical models, including applications in differential geometry \cite{ugail2011}, thin film equations \cite{barrett2004}, and the phase field crystal model \cite{backofen2007, cheng2008, hu2009, wise2009}.  
The conforming finite element approximation for (\ref{eqnbi}) necessitates $H^3$ conforming finite elements, typically involving intricate constructions of the finite element space and the variational formulation \cite{zenisek1970, bramble1970, Ciarlet78, BS02}.
Recently, a nonconforming $H^3$ finite element was proposed in \cite{wu2019}, where the element is composed of $H^1$ conforming finite elements and additional bubble functions. $C^0$ interior penalty discontinuous Galerkin (IPDG) and $C^1$-IPDG methods were proposed in \cite{gudi2011} for the sixth-order elliptic equations with clamped boundary conditions. 
To balance the weak continuity and the complexity associated with choosing penalty parameters, a family of $\mathcal{P}_m$ interior nonconforming finite element methods was proposed in \cite{wu2017}. Additionally, a mixed finite element method was introduced in \cite{droniou2019}, based on low order $H^1$ conforming finite elements, with an optimal error estimate under an appropriate regularity assumption.

The direct mixed finite element method, employing $C^0$ finite elements, offers an appealing approach for addressing high-order elliptic problems, such as the biharmonic problem \cite{zulehner2015ciarlet, li2022ima, zhang2022solving, li2025analysis} and the sixth-order problem \eqref{eqnbi}. This is primarily due to the boundary conditions, which facilitate the derivation of three entirely decoupled Poisson equations. This suggests that a plausible numerical solution could be attained by simply employing a finite element Poisson solver within the mixed formulation. However, while the implementation of the mixed finite element method is straightforward, its solution may not always be reliable, as the solution obtained from the Poisson problem might reside in a different Sobolev space compared to that of the original sixth order problem \eqref{eqnbi}. This discrepancy is evident in the fact that the solution to the Poisson problem typically belongs to $H^1(\Omega)$, whereas that of the sixth order problem \eqref{eqnbi} usually belongs to $H^3(\Omega)$. 
This phenomenon was identified in the context of the biharmonic equation with Navier boundary conditions, known as the Sapongyan paradox \cite{MR2270884, MR2479119}.
To confine the solution of the Poisson problem to $H^2(\Omega)$, an additional Poisson problem needs to be solved \cite{li2022ima}, particularly when the polygonal domain features a reentrant corner. 
For the sixth order problem \eqref{eqnbi}, achieving confinement of the solution to $H^3(\Omega)$ is not a trivial task.


The direct mixed formulation, which decomposes the problem into three Poisson equations, actually defines a weak solution in a larger function space compared to that of equation (\ref{eqnbi}). This mismatch in function spaces does not impact the solution in a polygonal domain where the largest interior angle is no more than $\pi/2$. However, when the largest interior angle exceeds $\pi/2$, the direct mixed method allows for additional singular functions, leading to a solution different from that of equation (\ref{eqnbi}). 
To confine the solution to the correct function spaces, we propose a modified mixed formulation aiming at eliminating the singular functions.
More specifically, we first rigorously establish that the space of the singular functions, or equivalently, their image space under the Laplace operator, is finite-dimensional. 
In particular, the dimension of the singular function space associated with a corner depends on the corresponding interior angle: it is $0$ if the angle lies in $(0,\pi/2]$, $1$ if in $(\pi/2,\pi)$, $2$ if in $(\pi, 3\pi/2]$, and $3$ if in $(3\pi/2, 2\pi)$.
Subsequently, we identify a basis for the singular function space, or equivalently, its image space. 
Finally, we formulate the modified mixed formulation by removing the solution component that resides in the singular function space.
The resulting formulation is shown to be well-posed, and the solution is equivalent to the original problem.

In turn, we introduce a numerical algorithm to solve the proposed mixed formulation, utilizing piecewise linear $C^0$ finite elements on quasi-uniform meshes.
Meanwhile, we conduct an error analysis on the finite element approximations for both the auxiliary functions and the solution $u$.
For the auxiliary functions, the errors in the $H^1$ norm are standard and have a convergence rate $h^{\min\{\frac{\pi}{\omega},1\}}$, where $\omega$ is the largest interior angle of the polygonal domain; the $L^2$ error estimates can be obtained using the duality argument.
For the approximation to the solution $u$, the error in the $H^1$ norm is bounded by: (i) the $H^1$ interpolation error of the solution $u$; (ii) the $H^{-1}$ error for the auxiliary functions; and (iii) the $H^1$ errors and the weighted $L^2$ error for the approximations to the additional intermediate Poisson problems that confine the solution to the correct function space.
Depending on the largest interior angle, the convergence rate for the $H^1$ error of the numerical solution is dominated by either the degree of the polynomials or the singularity of the intermediate functions.

In summary, we propose a $C^0$ finite element algorithm that reduces the sixth-order problem with simply supported boundary conditions to a system of second-order equations. The key contributions of this work are outlined as follows:
\begin{itemize}
    \item Compared to existing penalty methods and nonconforming approaches, the proposed method is simple and intuitive in its formulation, and a plausible numerical solution can be obtained using only a standard $C^0$ finite element Poisson solver.
    \item The direct mixed formulation, which decomposes the original problem into three Poisson problems, fails to maintain equivalence with the original problem when the largest interior angle exceeds $\pi/2$. In contrast, by carefully confining the intermediate functions to the appropriate function space, the proposed method remains valid for general polygonal domains, regardless of whether any interior angle exceeds $\pi/2$ or not.
    \item We rigorously derive optimal error estimates for the proposed method on quasi-uniform meshes using $C^0$ linear finite element polynomials.
    \item Based on the largest interior angle of the domain, we conduct numerical tests to compare the solutions obtained from the direct mixed finite element method and the proposed method. In addition, we evaluate the convergence rate of the proposed method.
\end{itemize}

The rest of the paper is organized as follows. In Section \ref{sec-2}, according to the general regularity theory for second order elliptic equations \cite{Kondratiev67, Grisvard1992, KMR197, KMR01, nazarov1994}, we introduce the weak solution of the sixth order problem (\ref{eqnbi}). Additionally, we discuss the orthogonal space of the image of the operator $-\Delta$ in $H_0^1(\Omega)$ and identify basis functions of this space. We then propose a modified mixed formulation and demonstrate the equivalence of the solution to that of the original sixth order problem.
In Section \ref{sec-3}, we present the finite element algorithm and derive error estimates on quasi-uniform meshes for both the solution $u$ and the auxiliary functions.
Finally, in Section \ref{sec-5}, we present numerical test results to validate the theory.



Throughout the paper, the generic constant $C>0$ in our estimates may vary across different occurrences. Its value depends on the computational domain but remains independent of the functions involved or the mesh level in the finite element algorithms.

\section{The sixth order problem}\label{sec-2}

\subsection{Well-posedness of the solution} 
Denote by $H^m(\Omega)$, $m\geq 0$, the Sobolev space consisting of functions whose $i$th derivatives are square integrable for $0\leq i\leq m$. Let $L^2(\Omega):=H^0(\Omega)$. If $m$ is not an integer, then it defines the fractional Sobolev space.
Denote by $\mathcal{D}(\Omega)$ the space of infinitely differentiable functions in $\Omega$ with compact support. We define $H^s_0(\Omega)$ to be the closure of $\mathcal{D}(\Omega)$ in $H^s(\Omega)$.
Recall that $H^s_0(\Omega)\subset H^s(\Omega)$ for $0<s\leq 1$ is the subspace consisting of functions with zero traces on the boundary $\pa\Omega$ \cite{lions1972}. We shall denote the norm $\|\cdot\|_{L^2(\Omega)}$ by $\|\cdot\|$ when there is no ambiguity about the underlying domain.
Recall that for   $D\subseteq \mathbb{R}^d$,
the fractional order Sobolev space $H^s(D)$ consists of distributions $v$ in $D$ satisfying
$$
\|v\|^2_{H^s(D)}:=\|v\|^2_{H^m(D)} + \sum_{|\alpha|= m}\int_{D}\int_{D} \frac{|\partial^\alpha v(x) - \partial^\alpha v(y)|^2 }{|x-y|^{d+2t}} dxdy <\infty,
$$
where $\alpha=(\alpha_1, \ldots, \alpha_d)\in\mathbb Z^d_{\geq0}$ is a multi-index such that $\pa^\alpha=\pa_{x_1}^{\alpha_1}\ldots\pa^{\alpha_d}_{x_d}$ and $|\alpha|=\sum_{i=1}^d\alpha_i$. 

We define the space
\begin{equation}\label{space}
V=\{ \phi \ | \ \phi\in H^3(\Omega), \phi|_{\partial \Omega}=0, \ \Delta \phi|_{\partial \Omega}=0\},
\end{equation}
then the variational formulation for equation (\ref{eqnbi}) is to find $u \in V$ such that,
\begin{eqnarray}\label{eqn.firstbi}
a(u, \phi):=\int_\Omega \nabla\Delta u \cdot \nabla\Delta \phi dx=\int_\Omega f \phi dx=(f, \phi),\quad \forall \phi\in V.
\end{eqnarray}

For (\ref{eqn.firstbi}), we have the following result.
\begin{lem}
Given $f \in H^{-1}(\Omega)$ for the variational formulation \eqref{eqn.firstbi}, there exists at most one solution in $V$.
\end{lem}
\begin{proof}
We postpone the proof of the existence of the solution to Theorem \ref{solequthm}. Assume that (\ref{eqn.firstbi}) has two solutions $u_1$ and $u_2$ in $V$. Let $\delta{u}= u_1-u_2$. Then we have
\be\label{Galerkin}
a(\delta{u}, \phi) = 0, \quad \phi\in V.
\ee
Note that $\delta {u} \in V $ implies $\Delta \delta{u} \in H_0^1(\Omega)$. In addition $\delta u \in H^2(\Omega) \cap H_0^1(\Omega)$. Then, by the Poincar\'e-type inequality, 
$$
\|\nabla \Delta \delta{u}\| \geq C_0 \|\Delta \delta{u}\|_{H^1(\Omega)} \geq C_0 \|\Delta \delta{u}\| \geq C\|\delta{u}\|_{H^2(\Omega)},
$$
where \cite[Theorem 2.2.3]{Grisvard1992} has been used in the last inequality.
By setting $\phi=\delta{u}$ in \eqref{Galerkin}, it follows
$$
0=a(\delta{u}, \delta{u}) = \| \nabla\Delta \delta{u}\|^2 \geq C\|\delta{u}\|_{H^2(\Omega)} = 0.
$$
Thus $\delta{u}=0$, which implies $u_1=u_2$ in $H^2(\Omega)$, and therefore $u_1=u_2$ in $V$. 
\end{proof}

%

 \subsection{The direct mixed formulation}

Intuitively, we can decouple (\ref{eqnbi}) into a system of three Poisson problems by introducing auxiliary functions $w$ and $v$, satisfying:
\begin{eqnarray}\label{eqn7}
\left\{\begin{array}{ll}
-\Delta w=f \quad {\rm{in}} \ \Omega,\\
\hspace{0.6cm}w=0 \quad {\rm{on}} \ \pa\Omega;
\end{array}\right.
\qquad
\left\{\begin{array}{ll}
-\Delta v=w \quad {\rm{in}} \ \Omega,\\
\hspace{0.6cm}v=0 \quad {\rm{on}} \ \pa\Omega;
\end{array}\right.
\qquad {\rm{and}}\qquad
\left\{\begin{array}{ll}
-\Delta \bar{u}=v \quad {\rm{in}} \ \Omega,\\
\hspace{0.6cm}\bar{u}=0 \quad {\rm{on}} \ \pa\Omega.
\end{array}\right.
\end{eqnarray}
We refer to (\ref{eqn7}) as the direct mixed formulation. Note that numerical solvers for the Poisson problems (\ref{eqn7}) are readily available, while numerical approximation of the sixth order problem (\ref{eqnbi}) is generally a daunting task.
The weak formulation of (\ref{eqn7}) is to find $\bar{u}, v, w\in H_0^1(\Omega)$ such that
\begin{subequations}\label{weaksys}
\begin{align}
A( w, \phi ) = & ( f, \phi), \quad \forall \phi \in H_0^1(\Omega), \\
A( v, \psi ) = & ( w, \psi), \quad \forall \psi \in H_0^1(\Omega), \\
A( \bar{u}, \tau ) = & ( v, \tau), \quad \forall \tau \in H_0^1(\Omega),
\end{align}
\end{subequations}
where
$$
A(\phi,\psi)=\int_\Omega \nabla \phi \cdot \nabla \psi dx.
$$
Assuming that the source term $f$ in (\ref{weaksys}) and (\ref{eqn.firstbi}) satisfies $f\in H^{-1}(\Omega)\subset V^*$, the solutions $\bar u, v, w$ of the Poisson problems in (\ref{weaksys}) are well-defined \cite{Evans}. The important question is whether the solution $\bar u$ in (\ref{weaksys}) is the same as the solution $u$ in (\ref{eqn.firstbi}).

%
To address this question, it is imperative to delve into the solution structure of the Poisson problem within a polygonal domain. This exploration will be undertaken in the subsequent subsection.

\subsection{Image of the Laplace operator in $H_0^1(\Omega)$ and its orthogonal space} 
\begin{figure}
\begin{center}
\begin{tikzpicture}[scale=0.2]
\draw[thick]
(-14,-11) -- (2,-11) -- (0,-2) -- (14,-2) -- (-3,12) -- (-14,-11);
\draw (7,-3) node {$\theta = 0$};
\draw (4.5,-8) node {$\theta = \omega$};
\draw (-2.5,10) node {$\Omega$};
\draw[thick][densely dotted] (8,-2) arc (0:283:8);
\draw[densely dotted] (4,-2) arc (0:283:4);
\draw (-6,-2) node {$K_\omega^R$};
\draw (-0.2,0.2) node {$K_\omega^{\tau R}$};
\draw[thick] (0,-2) node {$\bullet$} node[anchor = north west] {$Q$};
\end{tikzpicture}
\end{center}
\vspace*{-15pt}
    \caption{Domain $\Omega$ containing a reentrant corner.}
    \label{fig:Omega}
\end{figure}
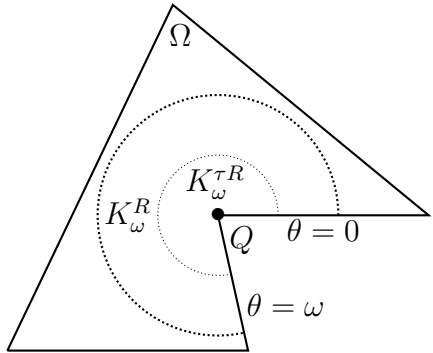

Assume that the polygonal domain $\Omega$ has at most one interior angle greater than $ \frac{\pi}{2}$. Let $\omega$ be the largest interior angle with the vertex $Q$. Without loss of generality, we set $Q$ as the origin and represent polar coordinates centered at the vertex $Q$ as $(r, \theta)$, where the interior angle $\omega$ is spanned by two half lines $\theta=0$ and $\theta=\omega$.
We construct a sector $K_\omega^R \subset \Omega$ at $Q$ with radius $R>0$ as
$$
K_\omega^R = \{(r\cos \theta, r\sin\theta) \in \Omega \ |\ 0 < r< R, 0< \theta < \omega\}.
$$
A sketch drawing of the domain $\Omega$ is depicted in Figure \ref{fig:Omega}.

To begin with, we introduce a general Poisson problem
\be\label{poissoneq}
-\Delta z = g  \text{ in } \Omega, \qquad z = 0 \text{ on } \partial \Omega.
\ee
Recall the space $V$ in \eqref{space}. For any function $\phi \in V$, it can be verified that $ -\Delta \phi  \in H_0^1(\Omega)$.
Then we have the following result.
\begin{lemma}\label{mapimage}
The mapping $-\Delta  : V \rightarrow H_0^1(\Omega)$ is injective and has a closed range,
where the subspace $V$ is given in \eqref{space}.
\end{lemma}
\begin{proof}
Let $z_1, z_2$ be functions in $V \subset H_0^1(\Omega)$ satisfying $\Delta z_1=\Delta z_2$. Then the function $g=-\Delta z_1=-\Delta z_2 \in H_0^1(\Omega)$.
By the Lax-Milgram Theorem for the Poisson problem \eqref{poissoneq}, it follows $z_1=z_2$ in $H_0^1(\Omega)$, and hence $z_1=z_2$ in $V$, demonstrating the injective nature of the mapping. 

Denote the image of the mapping by $\mathcal{M} \subset H_0^1(\Omega)$. 
Consider a sequence $\{g_i\}_{i=1}^\infty$ in $\mathcal{M}$ 
satisfying $g_i:=-\Delta z_i \rightarrow g$ for $z_i \in H^3(\Omega)$, which implies that $g_i\in \mathcal{M}$ is Cauchy and $g\in H_0^1(\Omega)$. 
We now show $\mathcal{M}$ is closed, namely, $g\in \mathcal{M}$.
By the regularity result for the elliptic equation, it holds
\be
\|z_m-z_n\|_{H^3(\Omega)} \leq C \|g_m-g_n\|_{H^1(\Omega)},
\ee
which implies $\{z_i\}_{i=1}^\infty$ is also Cauchy in $V$. Since the subspace $V$ is complete, it follows $z_i \rightarrow z \in V$, thus $-\Delta z_i \rightarrow -\Delta z \in \mathcal{M}$. Namely, $g=-\Delta z \in \mathcal{M}$.
Therefore, the space $\mathcal{M}$ is closed.
\end{proof}

Recall the image $\mathcal{M}$ of the mapping $-\Delta$ in $H_0^1(\Omega)$. Let $\mathcal{M}^\perp$ be its orthogonal complement in $H_0^1(\Omega)$. Namely, for any function $v \in H_0^1(\Omega)$, there exist unique $v_\mathcal{M} \in \mathcal M$ and $v_\perp \in \mathcal M^\perp$ such that
\be\label{gorth0}
v= v_\mathcal{M} + v_\perp, 
\ee
and
\be\label{gorth}
(\nabla v_\mathcal{M}, \nabla v_\perp) = 0.
\ee
In other words, $\mathcal{M} \oplus \mathcal{M}^\perp = H_0^1(\Omega)$.
By the definition of $\mathcal{M}$, the condition \eqref{gorth} is equivalent to 
\bes
(\nabla \Delta z, \nabla v_\perp) =0, \quad \forall z \in V.
\ees
In the following, we will show that the space $\mathcal{M}^\perp$ is finite-dimensional, allowing for the determination of its basis. 

Denote the $\ell$th side of $\partial \Omega$ by $\bar\Gamma_\ell$, where $\Gamma_\ell$ is open.
For $\forall \phi, \psi \in H^4(\Omega)$, Green's formula gives
\be\label{Green}
\int_{\Omega} \phi \Delta^2 \psi dx - \int_{\Omega} \Delta^2 \phi \psi dx = \sum_{\ell} \int_{\Gamma_\ell} \phi \partial_{\mathbf{n}} (\Delta \psi) - \partial_{\mathbf{n}} \phi \Delta \psi + \Delta \phi \partial_{\mathbf{n}} \psi -\partial_{\mathbf{n}} (\Delta \phi) \psi ds,
\ee
where $\mathbf{n}$ is the outward normal derivative.

We denote by $\mathcal{D}(\Omega)$ the space of infinitely differentiable functions with compact support in $\Omega$. Then we can show the following result.
\begin{lemma}\label{vmperp}
A function $v$ belongs to $\mathcal{M}^\perp$ if and only if 
$v\in H_0^1(\Omega)$ is the solution of the following (adjoint) boundary value problem 
\be\label{adjoint}
\Delta^2 v =0 \text{ in } \Omega, \qquad v =0, \quad \Delta v =0 \text{ on } \partial \Omega.
\ee
\end{lemma}
\begin{proof}
($\Rightarrow$) 
By \eqref{gorth}, it holds for $\forall v \in \mathcal{M}^\perp$ and $z \in V$, 
\be\label{gorth++}
(-\nabla \Delta z, \nabla v)=0.
\ee
In particular, for $\forall z \in \mathcal{D}(\Omega) \subset V$,
\be\label{gdzgv-}
\bal
(-\nabla \Delta z, \nabla v) = 0 = (z,\Delta^2  v),
\eal
\ee
which implies $\Delta^2 v =0$ in $\Omega$.

Define $D(\Delta^2, H^{-1}(\Omega))$ to be the maximal extension of the biharmonic operator in $H_0^1(\Omega)$:
$$
D(\Delta^2, H^{-1}(\Omega)): = \{ v \in H_0^1(\Omega):\ \Delta^2 v \in H^{-1}(\Omega) \}.
$$
Note that $v \in \mathcal{M}^\perp \subset D(\Delta^2, H^{-1}(\Omega))$. 
Now, suppose $v \in H^4(\Omega) \cap \mathcal{M}^\perp$. By Green's formula, it holds for $z \in V$,
\be\label{gdzgv}
\bal
(-\nabla \Delta z, \nabla v) = 0 =  (z,\Delta^2  v)  - \sum_{\ell} \int_{\Gamma_\ell} z \partial_{\mathbf{n}} (\Delta v) - \partial_{\mathbf{n}} z \Delta v + \Delta z \partial_{\mathbf{n}} v ds.
\eal
\ee
Then \eqref{gdzgv-}, \eqref{gdzgv} together with the boundary condition $z = \Delta z =0$ on every $\Gamma_\ell$ yields the boundary value condition $\Delta v = 0$ on $\Gamma_\ell$.
Given that $H^4(\Omega)$ is dense in $D(\Delta^2, H^{-1}(\Omega))$ \cite{lions1972}, the density argument asserts that the same boundary condition also holds for any $v \in \mathcal{M}^\perp \subset H_0^1(\Omega)$. Consequently, \eqref{adjoint} holds.\\
($\Leftarrow$) For $v\in H_0^1(\Omega)$ satisfying \eqref{adjoint}, it follows $v \in D(\Delta^2, H^{-1}(\Omega))$.
Suppose $v \in H^4(\Omega) \cap D(\Delta^2, H^{-1}(\Omega))$. By \eqref{adjoint} and Green's formula, \eqref{gdzgv} also holds for $\forall z \in V$. Since $H^4(\Omega)$ is dense in $D(\Delta^2, H^{-1}(\Omega))$, the equality $(-\nabla \Delta z, \nabla v) = 0$ also holds for $v \in D(\Delta^2, H^{-1}(\Omega))$ and implies $v \in \mathcal{M}^\perp$. 
\end{proof}

One of the main goals of this section is to show that $\mathcal{M}^\perp$ is finite dimensional and to identify the basis of $\mathcal{M}^\perp$.
Next, we introduce some pertinent functions in the domain $\Omega$.
\begin{definition}\label{defs}
Given $R>0$ such that $K_\omega^R \subset \Omega$. Let $N \geq 0$ be the largest integer satisfying $ N< \frac{2\omega}{\pi}$ with values specified in \Cref{TabomegaN}. Additionally, let $\tau \in (0,1)$ be a given parameter. \\
(i) For $1\leq i\leq N$, we define the $H^{-1}(\Omega)$ functions,
\begin{eqnarray}\label{ssol}
\xi_i(r, \theta; \tau, R):=\chi_i(r, \theta; \tau, R)+\zeta_i(r, \theta; \tau, R),
\end{eqnarray}
where
\begin{eqnarray}\label{sl2}
\bal
\chi_i(r, \theta; \tau, R)=&\eta(r; \tau, R)r^{-\frac{i\pi}{\omega}}\sin\left(\frac{i\pi}{\omega}\theta\right),
 \eal
\end{eqnarray}
with the cut-off function $\eta(r; \tau, R) \in C^\infty(\Omega)$ satisfying 
$\eta(r; \tau, R)=1$ for $0\leq r\leq \tau R$ and $\eta(r; \tau, R)=0$ for $r>R$, and
 $\zeta_i \in H_0^1(\Omega)$ is obtained by solving 
\be\label{l2part}
-\Delta \zeta_i = \Delta \chi_i  \text{ in } \Omega, \qquad \zeta_i = 0 \text{ on } \partial \Omega.
\ee
(ii) For $1\leq i\leq N$, we define $\sigma_i \in H_0^1(\Omega)$ satisfying
\be\label{H1part}
-\Delta \sigma_i = \xi_i  \text{ in } \Omega, \qquad \sigma_i = 0 \text{ on } \partial \Omega.
\ee
\end{definition}


\begin{table}[!htbp]\tabcolsep0.03in
\caption{The range of $\frac{\pi}{\omega}$ and the value of $N$ for different $\omega$ in Definition \ref{defs}.}
\begin{tabular}[c]{|c|c|c|c|c|}
\hline
$\omega$ & $(0,\frac{\pi}{2}]$ & $(\frac{\pi}{2},\pi)$ & $(\pi,\frac{3\pi}{2}] $ & $(\frac{3\pi}{2}, 2\pi)$ \\
\hline
$\frac{\pi}{\omega}$ & $(2,\infty)$ & $(1, 2)$ & $(\frac{2}{3},1)$ & $(\frac{1}{2},\frac{2}{3})$ \\
\hline
$N$ & 0 & 1 & 2 & 3 \\
\hline
\end{tabular}\label{TabomegaN}
\end{table}



\begin{remark}
If $N=0$, both the function sets $\{\xi_i\}_{i=1}^N$ and $\{\sigma_i\}_{i=1}^N$ are empty. The functions $\xi_i$, $i=1,\ldots, N$ defined in \eqref{ssol} are not in $H^1(\Omega)$.
\end{remark}

Define $D(\Delta, H^{-1}(\Omega))$ to be the maximal extension of the Laplace operator in $H^{-1}(\Omega)$ \cite{lions1972},
$$
D(\Delta, H^{-1}(\Omega)): = \{ v \in H^{-1}(\Omega):\ \Delta v \in H^{-1}(\Omega) \}.
$$
For the functions $\xi_i$ in \Cref{defs}, the following properties hold. 
\begin{lem}\label{slem}
Given $\eta\in C^\infty(\Omega)$ in Definition \ref{defs}, the functions $\xi_i \in D(\Delta, H^{-1}(\Omega))$, $i=1, \ldots, N$, are uniquely defined and satisfy
\be\label{seqn}
-\Delta \xi_i =0 \text{ in } \Omega, \qquad \xi_i =0 \text{ on } \partial \Omega.
\ee
Moreover, $\xi_i$ depends on the domain $\Omega$, but not on $\tau$ and $R$. Namely, for any positive numbers $\tau_1, \tau_2$ and $R_1, R_2$, it holds
\be\label{paraindep}
\xi_i(r, \theta) :=\xi_i(r, \theta; \tau_1, R_1) = \xi_i(r, \theta; \tau_2, R_2).
\ee
\end{lem}
\begin{proof}
For $\chi_i$ given in (\ref{sl2}) with $1\leq i\leq N$, it can be verified that $\chi_i \in C^\infty(\Omega \setminus K_\omega^\delta)$ for any $\delta>0$ and $\chi_i = 0$ for $(r\cos \theta, r\sin\theta) \in \Omega \setminus K_\omega^R$. Moreover, $\Delta \chi_i = 0$ if $r<\tau R$ and $r>R$.
These imply that $\Delta \chi_i\in C^{\infty}(\Omega) \subset L^2(\Omega)$.
Given $\eta \in C^\infty(\Omega)$, the explicit function $\chi_i $ belonging to $ H^{-1}(\Omega)$  in (\ref{sl2}) is uniquely defined, so is $\Delta \chi_i$.
In addition, $\zeta_i\in  H^1_0(\Omega)$ is uniquely defined via  (\ref{l2part}).
Therefore, $\xi_i$ in (\ref{ssol}) is uniquely defined due to the uniqueness of $\chi_i$ and $\zeta_i$.

Taking $-\Delta$ on both side of (\ref{ssol}) yields
\bes
\bal
-\Delta \xi_i = -\left( \Delta \chi_i + \Delta \zeta_i \right) = 0,
\eal
\ees
where (\ref{l2part}) have been applied.
In addition, $\xi_i=0$ on $\partial \Omega$ is obtained by $\chi_i=0$ and $\zeta_i=0$ on $\partial \Omega$.

Next, we prove (\ref{paraindep}). By taking $\delta  \in (0, \min\{ \tau_1 R_1, \tau_2 R_2 \})$, it follows $K_\omega^\delta \subset K_\omega^{\tau_1 R_1} \cap K_\omega^{\tau_2 R_2} \subset \Omega$. By (\ref{sl2}), we have
$$
\chi_i(r, \theta; \tau_1, R_1) - \chi_i(r, \theta; \tau_2, R_2) = 0, \quad (r\cos \theta, r\sin\theta) \in K_\omega^\delta.
$$
Recall that $\chi_i(r, \theta; \tau_j, R_j) \in C^\infty(\Omega \setminus K_\omega^\delta)$, $j=1,2$. Then it follows
$$
\chi_i(r, \theta; \tau_1, R_1) - \chi_i(r, \theta; \tau_2, R_2) \in C^\infty(\Omega).
$$
Since $\zeta_i(r, \theta; \tau_j, R_j)\in H_0^1(\Omega)$, $j=1,2$,  we have
\begin{align*}
\tilde{\xi}_i := & \xi_i(r, \theta; \tau_1, R_1) - \xi_i(r, \theta; \tau_2, R_2) \\
= &\zeta_i(r, \theta; \tau_1, R_1) - \zeta_i(r, \theta; \tau_2, R_2) + \left(\chi_i(r, \theta; \tau_1, R_1) - \chi_i(r, \theta; \tau_2, R_2) \right) \in H_0^1(\Omega).
\end{align*}
Meanwhile, from (\ref{seqn}), we have
\be\label{lapxi}
\Delta \tilde{\xi}_i = \Delta \xi_i(r, \theta; \tau_1, R_1) - \Delta \xi_i(r, \theta; \tau_2, R_2)=0 \text{ in } \Omega, \qquad \tilde{\xi}_i =0 \text{ on } \partial \Omega.
\ee
By applying the Lax-Milgram Theorem to \eqref{lapxi}, it is established that $\tilde{\xi}_i = 0$, indicating the validity of (\ref{paraindep}).
\end{proof}

\begin{remark}
\Cref{slem} implies that $\xi_i(r, \theta; \tau, R)$ in \Cref{defs} can be replaced by $\xi_i(r,\theta)$. Moreover, the $H^{-1}(\Omega)$  functions $\xi_i(r, \theta)\not\equiv0$, because otherwise we have $ \chi_i =- \zeta_i \in H^1_0(\Omega)$, which contradicts  the fact that $\chi_i \not \in H^{1}(\Omega)$.
\end{remark}


Subsequently, for the functions $\sigma_i$ in $H_0^1(\Omega)$ defined in \Cref{defs}, the following property is satisfied.
\begin{lem}\label{slemH1}
The functions $\sigma_i \in D(\Delta^2, H^{-1}(\Omega))$, $i=1, \ldots, N$, in Definition \ref{defs} are uniquely defined and satisfy
\be\label{seqnH1}
\Delta^2 \sigma_i =0 \text{ in } \Omega, \qquad \sigma_i =0, \quad \Delta \sigma_i =0 \text{ on } \partial \Omega.
\ee
\end{lem}

\begin{proof}
Note that $\sigma_i$ is obtained through the Poisson problem (\ref{H1part}) with $\xi_i$ as the source term. From Lemma \ref{slem}, $\xi_i$ is uniquely defined, which yields the uniqueness $\sigma_i$.
Applying $-\Delta$ to \eqref{H1part} in conjunction with (\ref{seqn}) yields (\ref{seqnH1}).
\end{proof}

For both functions $\xi_i$ and $\sigma_i$, we have the following results.
\begin{lemma}\label{linin}
(a) The functions $\xi_i(=\Delta \sigma_i)$, $i=1,2,\ldots, N$, are linearly independent.
(b) The functions $\sigma_i, \nabla \sigma_i$, $i=1,2,\ldots, N$, are also linearly independent, respectively.
\end{lemma}
\begin{proof}
(a) $\xi_i(r, \theta)\not\equiv0$, because otherwise we have $ \chi_i =- \zeta_i \in H^1_0(\Omega)$, which contradicts  the fact that $\chi_i\notin H^1_0(\Omega)$.
We assume that
\be\label{sumxi}
\sum_{i=1}^N C_i \xi_i =0,
\ee
where $C_i, \ i=1,2,\ldots, N$ are some constants.
Plugging \eqref{ssol} into \eqref{sumxi} gives
\bes
\sum_{i=1}^N C_i \chi_i = - \sum_{i=1}^N C_i \zeta_i \in H_0^1(\Omega).
\ees
Note that $\chi_i \not \in H_0^1(\Omega), \ i=1,2,\ldots, N$. Therefore, it holds
\be\label{sindep}
\sum_{i=1}^N C_i \chi_i=0.
\ee
Multiplying (\ref{sindep}) by $r^{-\frac{\pi}{\omega}}$, we have $C_N r^{-\frac{\pi}{\omega}} \chi_N=-\sum_{i=1}^{N-1} C_i r^{-\frac{\pi}{\omega}} \chi_i \in H^{-1}(\Omega)$, which contradicts the fact that
$C_N r^{-\frac{\pi}{\omega}} \chi_N \not \in H^{-1}(\Omega)$. Thus, it follows
$
C_N =0.
$
For $i=2, \ldots, N$, multiplying (\ref{sindep}) by $r^{-\frac{i\pi}{\omega}}$, the same argument yields
$
C_{N+1-i}=0.
$
Thus, $\xi_i$, $i=1,2,\ldots, N$, are linearly independent.\\
(b) We assume $\sum_{i=1}^N C'_i \sigma_i =0$ for some constants $C'_i$ and apply $-\Delta$ to both sides of the equation, it follows
$$
\sum_{i=1}^N C'_i \xi_i =0.
$$
By (a), we have $C'_i=0$, $i=1,\ldots, N$, which implies $\sigma_i$, $i=1,2,\ldots, N$, are linearly independent. The linear independence of $\nabla \sigma_i$ can be proved similarly.
\end{proof}


\begin{corollary}\label{spantoM}
The space $\text{span}\{\sigma_i, i=1,\ldots, N\} \subset \mathcal{M}^\perp$, and the dimension of $\mathcal{M}$ satisfies $\dim(\mathcal{M}^\perp) \geq N$.
\end{corollary}
The proof follows from \Cref{vmperp}, \Cref{slemH1}, and \Cref{linin}.

\begin{lemma}\label{xitosig}
For any function $v \in H_0^1(\Omega)$, it holds
\be\label{grotheq}
\langle v, \xi_i \rangle = (\nabla v, \nabla\sigma_i),\quad \forall i =1, \ldots, N,
\ee
where $\xi_i$ and $\sigma_i$ are given in \Cref{defs}. 
In particular, if $v\in \mathcal{M}$, it holds
\be\label{grotheq+}
\langle v, \xi_i \rangle = (\nabla v, \nabla\sigma_i)=0,\quad \forall i =1, \ldots, N.
\ee
\end{lemma}
\begin{proof}
Multiplying (\ref{H1part}) by $v \in H_0^1(\Omega)$ and applying Green's formula yield \eqref{grotheq}.
Since $\sigma_i \in \mathcal{M}^\perp$, \eqref{grotheq+} follows from \eqref{gorth}.
\end{proof}

To determine the dimension of $\mathcal{M}^\perp$, we let $\lambda_{i}^2$ be the eigenvalues to the following one dimensional problem
\be\label{egneqn}
-\partial_{\theta \theta}\phi_{i} = \lambda_{i}^2 \phi_{i} \quad \text{ in } (0,\omega), \qquad \phi(0)=\phi(\omega) =0.
\ee
For $i\geq 1$, it is clear that when $\lambda_i>0$,
\be\label{eigen}
\lambda_{i} = \frac{i\pi}{\omega}, \quad \phi_{i}= \sqrt{\frac{2}{\omega}} \sin \left(\frac{i\pi}{\omega} \theta\right).
\ee
In addition, we also recall the following result for the Poisson problem (\ref{poissoneq}).
\begin{lemma}\label{Poissonreg}
Assume that $g\in H^1(\Omega)$, and $\lambda_i=\frac{i\pi}{\omega}$, $1\leq i \leq N$ for $N$ given in \Cref{TabomegaN}, is not an integer, namely $\omega \not = \frac{\pi}{2}$ and $\omega \not = \frac{3\pi}{2}$.
Then the solution $z$ of the Poisson problem (\ref{poissoneq}) from the space $H^{1+\alpha}(\Omega) \cap H_0^1(\Omega)$ for $\alpha<\frac{\pi}{\omega}$ possesses the asymptotic representation in the neighborhood of $Q$,
\be\label{zasym}
z(x) = \tilde{z}(x) + \eta(r)\sum_{i=1}^{N} d_i (i\pi)^{-\frac{1}{2}}r^{\frac{i\pi}{\omega}}\sin\left(\frac{i\pi \theta}{\omega}\right),
\ee
where $\tilde{z}(x) \in H^3(\Omega)\cap H_0^1(\Omega)$ and the coefficients $d_i$ are defined by
\be\label{dgxi}
d_i = \langle g,\xi_i \rangle, \quad i=1, \ldots, m.
\ee
Moreover, it follows that
\be\label{preg0}
\|\tilde{z}\|_{H^{3}(\Omega)} + \sum_{i=1}^N|d_i| \leq C\|g\|_{H^1(\Omega)}.
\ee
\end{lemma}
\begin{proof}
The proof can be found in Theorem 3.4 in \cite{nazarov1994} and Section 2.7 in \cite{Grisvard1992}.
\end{proof}

Based on \Cref{Poissonreg}, we can identify the dimension of $\mathcal{M}^\perp$ as follows.

\begin{lem}\label{lem92}
Under the condition in \Cref{Poissonreg}.
The dimension of $\mathcal{M}^\perp$ is equal to the cardinality of the set $\{\lambda_i: 0<\lambda_{i}<2\}$, namely
\ben
\dim(\mathcal{M}^\perp) = \text{card }\{\lambda_{i} : 0<\lambda_{i}<2\} = N,
\een
where the condition $0<\lambda_{i}<2$ corresponds to $1\leq i\leq N$.
\end{lem}
\begin{proof}
For $\forall v \in \mathcal{M}^\perp$, by \eqref{gorth} it holds
\be\label{gvorth}
(\nabla g, \nabla v) = 0, \quad \forall g\in \mathcal{M}.
\ee
For the Poisson problem \eqref{poissoneq} with $g\in \mathcal{M} \subset H_0^1(\Omega)$, \Cref{mapimage} implies that its solution $z\in V\subset H^3(\Omega)$.
By \Cref{Poissonreg}, $z \in H^3(\Omega)$ is equivalent to the fact that 
for $\lambda_i\in (0,2)$, the coefficients
\be\label{specdi}
d_i=\langle g,\xi_i \rangle=(\nabla g, \nabla\sigma_i)=0,
\ee
where we have used \eqref{grotheq} in the second equality. If $\lambda_i$ is not an integer, $\lambda_i \in (0,2)$ corresponds to the integer $i \in [1,N]$. \eqref{gvorth} and \eqref{specdi} imply that 
$\mathcal{M}^\perp \subset \text{span}\{\sigma_i,\ i=1,\ldots, N\}$, which together with \Cref{spantoM} gives the conclusion.
\end{proof}

The cases of $\omega=\frac{\pi}{2}$ and $\frac{3\pi}{2}$ are not covered by \Cref{Poissonreg}. To address this limitation, we introduce the following additional result.
\begin{remark}\label{regpi}
The asymptotic representation of the solution $z$ to problem (\ref{poissoneq}) typically involves two types of singular functions depending on $\lambda_i=i\pi/\omega$:
\begin{subequations}\label{Sjcoef}
\begin{align}
S_i = & (i\pi)^{-\frac{1}{2}} r^{\frac{i\pi}{\omega}}\sin\left(\frac{i\pi \theta}{\omega}\right) \quad \text{when } \lambda_i  \text{ is not an integer},\\
S_i = &  r^{\frac{i\pi}{\omega}}\left(\ln r \sin\left(\frac{i\pi \theta}{\omega}\right) + \theta \cos\left(\frac{i\pi \theta}{\omega}\right) \right)\quad \text{otherwise}.
\end{align}
\end{subequations}
Specially, the coefficient of the term in (\ref{Sjcoef}b) depends locally on the restriction of the data $g$ to any neighborhood of the corner \cite{Grisvard1985}.
If $g\in H_0^1(\Omega)$, the solution of problem (\ref{poissoneq}) has the expansion \cite{Grisvard1985}
\be
z - \sum_{0<\lambda_i<2} d_iS_i \in  H^3(\Omega),
\ee
where $d_i$ is given by (\ref{dgxi}).
In other words, when the source term $g\in H_0^1(\Omega)$, the singular function $S_i$ in (\ref{Sjcoef}b) with $\lambda_i = i\pi/\omega = 2$ vanishes in the asymptotic representation of $z$.
\end{remark}

\begin{corollary}\label{vdecom}
The dimension of $\mathcal{M}^\perp$  satisfies $\dim(\mathcal{M}^\perp)=N$. Moreover, 
$$
\text{span}\{\sigma_i, \ i=1, \ldots, N\} = \mathcal{M}^\perp.
$$
\end{corollary}
\begin{proof}
The proof follows from \Cref{linin}, \Cref{spantoM}, \Cref{lem92}, and \Cref{regpi}.
\end{proof}

For $\forall v \in H_0^1(\Omega)$, \Cref{vdecom} and \eqref{gorth} imply that $(\nabla v_\mathcal{M}, \nabla \sigma_i)=0$, $1\leq i \leq N$ and that there exists a unique decomposition,
\be\label{mapv}
v=v_\mathcal{M} + \sum_{i=1}^N c_i \sigma_i,
\ee
where $v_\mathcal{M} \in \mathcal{M}$ and the coefficients $c_i$ are uniquely determined by the linear system,
\be\label{scoef}
\sum_{i=1}^N c_i (\nabla\sigma_i, \nabla\sigma_j)=(\nabla v, \nabla \sigma_j), \quad j=1, \ldots, N.
\ee
By \Cref{xitosig}, it holds that for $\forall \phi \in H_0^1(\Omega)$,
\be\label{sigweak}
(\nabla\sigma_j, \nabla \phi) = \langle \xi_j, \phi\rangle, \quad j=1,\ldots, N.
\ee
Therefore, the linear system (\ref{scoef0}) is equivalent to the following linear system
\be\label{scoef0}
\bal
\sum_{i=1}^N c_i \langle\sigma_i, \xi_j\rangle=\langle v, \xi_j \rangle, \quad j=1, \ldots, N.
\eal
\ee

\begin{lem}\label{UniC}
The linear system (\ref{scoef0}) or (\ref{scoef}) admits a unique solution $c_i$, $i=1,\ldots,N$.
\end{lem}
\begin{proof}
Since (\ref{scoef0}) and (\ref{scoef}) are equivalent, we only need to consider (\ref{scoef}), which is a finite-dimensional linear system. The existence of the solution is equivalent to the uniqueness.
Let $\bar{c}_i$ be the difference between two possible solutions; it follows
$$
\left(\sum_{i=1}^N \bar{c}_i\nabla\sigma_i, \nabla\sigma_j \right)=\sum_{i=1}^N \bar{c}_i (\nabla\sigma_i, \nabla\sigma_j)=0, \quad j=1, \ldots, N.
$$
A linear combination in terms of $\nabla\sigma_j$ gives
$$
\left(\sum_{i=1}^N \bar{c}_i\nabla\sigma_i, \sum_{j=1}^N \bar{c}_j\nabla\sigma_j \right)=0,
$$
which means $\left\|\sum_{i=1}^N \bar{c}_i\nabla\sigma_i\right\|=0$, thus we have
$$
\sum_{i=1}^N \bar{c}_i\nabla\sigma_i =0.
$$
\Cref{linin} indicates $\bar{c}_i=0$, $i=1,\ldots,N$. Thus, the conclusion holds.
\end{proof}

\subsection{The modified mixed formulation}
Based on the discussion above, we propose a modified mixed formulation for (\ref{eqnbi}),
\begin{eqnarray}\label{eqnnew}
\left\{\begin{array}{ll}
-\Delta w=f \quad {\rm{in}} \ \Omega,\\
\hspace{0.6cm}w=0 \quad {\rm{on}} \ \pa\Omega;
\end{array}\right.
\qquad
\left\{\begin{array}{ll}
-\Delta v=w \quad {\rm{in}} \ \Omega,\\
\hspace{0.6cm}v=0 \quad {\rm{on}} \ \pa\Omega;
\end{array}\right.
\qquad\left\{\begin{array}{ll}
-\Delta \tilde{u}=v - \sum_{i=1}^N c_i \sigma_i \quad {\rm{in}} \ \Omega,\\
\hspace{0.6cm}\tilde{u}=0 \quad {\rm{on}} \ \pa\Omega,
\end{array}\right.
\end{eqnarray}
where $\sigma_i$ are given in (\ref{H1part}) and $c_i$ are given by (\ref{scoef}).

The modified mixed weak formulation for (\ref{eqnnew}) is to find $w,v,\tilde{u}\in H_0^1(\Omega)$ such that
\begin{subequations}\label{weaksysnew}
\begin{align}
A(w,\phi) = & (f, \phi), \\
A(v,\phi) = & (w, \psi), \\
A(\tilde{u}, \tau) = & \left(v - \sum_{i=1}^N c_i \sigma_i, \tau \right),
\end{align}
\end{subequations}
for any $\phi, \psi, \tau \in H_0^1(\Omega)$.

Next, we show that $\tilde u$ is the weak solution to the variational formulation \eqref{eqn.firstbi}.
\begin{theorem}\label{solequthm}
Given $f \in H^{-1}(\Omega)$, let $\tilde{u}$ be the solution of the modified mixed weak formulation \eqref{weaksysnew}. Then $\tilde{u}$ is equivalent to the solution of the weak formulation \eqref{eqn.firstbi}, namely, $u = \tilde{u}$ in $V$, and vice versa.
\end{theorem}
\begin{proof}
Note that $v, \sigma_i\in H_0^1(\Omega)$. Thus $v - \sum_{i=1}^N c_i \sigma_i \in H_0^1(\Omega)$. By \eqref{scoef0}, it holds $d_j=\langle v - \sum_{i=1}^N c_i \sigma_i, \xi_j\rangle=0$, $j=1,\ldots, N$. Therefore, by applying Lemma \ref{Poissonreg} and \Cref{regpi} to the last Poisson equation in (\ref{eqnnew}), it follows $\tilde{u} \in H^3(\Omega) \cap H_0^1(\Omega)$. Since $\Delta \tilde{u}|_{\partial \Omega}=-(v - \sum_{i=1}^N c_i \sigma_i)|_{\partial \Omega}=0$, it follows $\tilde{u}\in V$.

On the other hand,
$$
-\Delta^3 \tilde{u} = \Delta^2 v - \sum_{i=1}^N c_i \Delta^2 \sigma_i = -\Delta (-\Delta v) = -\Delta w = f,
$$
where we have used the result (\ref{seqnH1}). Thus, we have $\tilde{u}\in V$ satisfying (\ref{eqn.firstbi}).
Finally, by the uniqueness of the solution of (\ref{eqn.firstbi}) in $V$, the conclusion holds.
\end{proof}

Therefore, by Theorem \ref{solequthm}, the solution $u$ of the sixth order problem (\ref{eqnbi}) satisfies
\begin{eqnarray}\label{eqnnew+}
\left\{\begin{array}{ll}
-\Delta w=f \quad {\rm{in}} \ \Omega,\\
\hspace{0.6cm}w=0 \quad {\rm{on}} \ \pa\Omega;
\end{array}\right.
\qquad
\left\{\begin{array}{ll}
-\Delta v=w \quad {\rm{in}} \ \Omega,\\
\hspace{0.6cm}v=0 \quad {\rm{on}} \ \pa\Omega;
\end{array}\right.
\qquad\left\{\begin{array}{ll}
-\Delta u=v - \sum_{i=1}^N c_i \sigma_i \quad {\rm{in}} \ \Omega,\\
\hspace{0.6cm}u=0 \quad {\rm{on}} \ \pa\Omega,
\end{array}\right.
\end{eqnarray}
The corresponding weak formulation is to find $w,v,u\in H_0^1(\Omega)$ such that for any $\phi, \psi, \tau\in H_0^1(\Omega)$,
\begin{subequations}\label{weaknew}
	\begin{align}
A(w,\phi) = & (f, \phi), \\
A(v,\psi) = & (w, \psi), \\
A(u, \tau) = & \left(v - \sum_{i=1}^N c_i \sigma_i, \tau \right),
	\end{align}
\end{subequations}
where $c_i$, $i=1,\ldots, N$, are given in (\ref{scoef}).

\begin{remark}
For the following cases, the modified mixed formulation (\ref{eqnnew+}) is identical to the direct mixed formulation (\ref{eqn7}):
(i) $N=0$, which happens if $\omega\leq \frac{\pi}{2}$ as shown in Table \ref{TabomegaN}; (ii) the boundary of domain $\Omega$ is sufficiently smooth; (iii) $c_i=0$, $i=1,\ldots, N$ in (\ref{scoef0}) or (\ref{scoef}), which is possible for some source term $f$ such that the solution $v\in \mathcal{M}$ in  (\ref{eqnnew+}).
\end{remark}

\begin{lem}\label{cmap}
The mapping $v\rightarrow v_\mathcal{M}$ in (\ref{mapv}) defines a norm non-increasing mapping $H_0^1(\Omega)\rightarrow \mathcal{M}$ in the sense
$$
\|\nabla v_\mathcal{M}\| \leq \|\nabla v\|.
$$
\end{lem}
\begin{proof}
Multiplying (\ref{mapv}) by $-\Delta v_\mathcal{M}$, integrating over the domain $\Omega$, and applying Green's Theorem give
\be\label{invvm}
(\nabla v, \nabla v_\mathcal{M}) = (\nabla v_\mathcal{M}, \nabla v_\mathcal{M}) + \left(\sum_{i=1}^N c_i \nabla \sigma_i, \nabla v_\mathcal{M} \right).
\ee
Note that
$$
(\nabla v_\mathcal{M}, \nabla \sigma_j)=\left(\nabla \left(v-\sum_{i=1}^N c_i \sigma_i\right), \nabla \sigma_j \right)=(\nabla  v, \nabla \sigma_j)- \sum_{i=1}^N c_i\left( \nabla\sigma_i, \nabla \sigma_j \right) = 0, \quad j=1,\ldots, N,
$$
where we have used (\ref{scoef}) in the last equality.
For the last term in (\ref{invvm}), it follows
$$
\left(\sum_{i=1}^N c_i \nabla \sigma_i, \nabla v_\mathcal{M} \right)=0.
$$
Then, applying H\"older's inequality to (\ref{invvm}), it follows
$$
\|\nabla v_\mathcal{M}\|^2 =(\nabla v, \nabla v_\mathcal{M}) \leq \|\nabla v_\mathcal{M}\|\|\nabla v\|,
$$
which gives the conclusion.
\end{proof}
In addition,  we have the following regularity result.
\begin{theorem}\label{wellposed}
\label{lemsreg}
Given $f \in H^{-1}(\Omega)$, for
$w,u,v$ in (\ref{eqnnew+}), it follows
\begin{subequations}\label{weakreg0}
\begin{align}
\|w\|_{H^1(\Omega)} \leq & C \|f\|_{H^{-1}(\Omega)},\\
\|v\|_{H^1(\Omega)} \leq & C \|f\|_{H^{-1}(\Omega)},\\
\|u\|_{H^3(\Omega)} \leq & C \|f\|_{H^{-1}(\Omega)}.
\end{align}
\end{subequations}
\end{theorem}
\begin{proof}
The estimate (\ref{weakreg0}a) is a direct consequence of the fact that the Laplace operator is an isomorphism between $H^1_0(\Omega)$ and $H^{-1}(\Omega)$.
In a similar fashion,
$$
\|v\|_{H^1(\Omega)} \leq \|w\|_{H^{-1}(\Omega)} \leq C\|w\|_{H^1(\Omega)}  \leq C \|f\|_{H^{-1}(\Omega)},
$$
which gives the estimate (\ref{weakreg0}b).
By Theorem \ref{solequthm}, it follows $u \in V$. Moreover, (\ref{preg0}) gives
\be\label{preg}
\|u\|_{H^{3}}\leq C\left\|v - \sum_{i=1}^N c_i \sigma_i \right\|_{H^1(\Omega)} = C\| v_\mathcal{M}\|_{H^1(\Omega)} \leq C\| v\|_{H^1(\Omega)}\leq C \|f\|_{H^{-1}(\Omega)},
\ee
where we have used Lemma \ref{cmap} and Poincar\'e inequailty.
\end{proof}


\section{The finite element method}\label{sec-3}

In this section, we introduce a linear $C^0$ finite element method for solving the sixth order problem (\ref{eqnbi}). Subsequently, we conduct a finite element error analysis.

\subsection{The finite element algorithm}\label{fem}
Let $\mathcal{T}_n$ denote a triangulation of $\Omega$ consisting of shape-regular triangles, and let $S_n \subset H^1_0(\Omega)$ be the $C^0$ Lagrange linear finite element space associated with $\mathcal{T}_n$.
Then we proceed to propose the  finite element algorithm.




\begin{algorithm}\label{femalg}
We define the finite element solution of the sixth order problem (\ref{eqnbi}) by employing the decoupling presented in (\ref{weaknew}) as follows.
\begin{itemize}
\item{Step 1.}  Find the finite element solution $w_n\in S_n$ of the Poisson equation
\be\label{femw}
A(w_n, \phi)=(f, \phi), \qquad \forall \phi\in S_n.
\ee
\item{Step 2.}  Find the finite element solution $v_n\in S_n$ of the Poisson equation
\be\label{femv}
A(v_n, \psi)=(w_n, \psi), \qquad \forall \psi\in S_n.
\ee
\item{Step 3.} With $\chi_i$, $i=1,\ldots, N$ defined in (\ref{sl2}), we compute the finite element solution $\zeta_{i,n}\in S_n$ of the Poisson equation
\be\label{femb}
A(\zeta_{i,n}, \phi)=(\Delta \chi_i, \phi), \qquad \forall \phi\in S_n,
\ee
and set $\xi_{i,n}=\zeta_{i,n}+ \chi_i$.
\item{Step 4.}  Find the finite element solution $\sigma_{i,n}\in S_n$, $i=1,\ldots, N$ of the Poisson equation
\be\label{fems}
A(\sigma_{i,n}, \phi)=(\xi_{i,n}, \phi), \qquad \forall \phi\in S_n.
\ee
\item{Step 5.} Find the coefficient $c_{i,n}\in \mathbb R$ by solving the linear system
\be\label{scoefh}
\sum_{i=1}^N c_{i,n} \langle \sigma_{i,n}, \xi_{j,n} \rangle= \langle v_n, \xi_{j,n} \rangle, \quad j=1, \ldots, N.
\ee
\item{Step 6.} Find the finite element solution $u_n\in S_n$ of the Poisson equation
\be\label{femu}
A(u_n, \tau)=\left(v_n - \sum_{i=1}^N c_{i,n} \sigma_{i,n}, \tau \right), \qquad \forall \tau\in S_n.
\ee
\end{itemize}
\end{algorithm}

\begin{remark}
According to (\ref{femb}), $\zeta_{i,n}\in S_n$,  while $\xi_{i,n} \in H^{-1}(\Omega)$ but $\xi_{i,n} \not\in S_n$. 
In addition, the finite element approximations in Algorithm \ref{femalg} are well defined based on the Lax-Milgram Theorem.
\end{remark}

For the functions in \eqref{femalg}, the following results hold.
\begin{lemma}\label{linin2}
(a) The $H^{-1}(\Omega)$ functions $\xi_{i,n}$, $i=1,2,\ldots, N$, are linearly independent.\\
(b) The functions $\sigma_{i,n}, \nabla \sigma_{i,n}$, $i=1,2,\ldots, N$, are also linearly independent, respectively.
\end{lemma}
\begin{proof}
(a) The proof is similar to the proof of Theorem \ref{linin}(a).\\
(b) We assume that
$\sum_{i=1}^N C'_i \sigma_{i,n} =0$ for some constants $C'_i$. The combination of (\ref{fems}) gives
$$
\left(\sum_{i=1}^N C'_i\xi_{i,n}, \phi\right)=  A\left(\sum_{i=1}^N C'_i\sigma_{i,n}, \phi\right)=0.
$$
By (a), we have $C'_i=0$, $i=1,\ldots, N$, which implies $\sigma_{i,n}$, $i=1,2,\ldots, N$, are linearly independent. The linear independence of $\nabla \sigma_{i,n}$ can be proved similarly.
\end{proof}

\subsection{Optimal error estimates on quasi-uniform meshes}

Suppose that the mesh $\maT_n$ consists of quasi-uniform triangles with size $h$. Recall the interpolation error estimates \cite{Ciarlet78} on $\maT_n$ for any $z \in H^{1+s}(\Omega)$, $s>0$, 
\be\label{interr}
\| z - z_I \|_{H^m(\Omega)} \leq Ch^{{\min\{s+1, 2 \}}-m}\|z\|_{H^{\min\{s+1, 2 \}}(\Omega)},
\ee
where $m= 0, 1$ and $z_I\in S_n$ represents the nodal interpolation of $z$. Let $z_n \in S_n$ be the finite element solution of the Poisson equation (\ref{poissoneq}) in the polygonal domain, if $z \in H^{1+s}(\Omega)$, $s>0$, the standard error estimate  \cite{Ciarlet78, li2022graded} yields
\be\label{zerr}
\|z-z_n\|_{H^1(\Omega)} \leq Ch^{\min\{s, 1\}}\|z\|_{H^{1+\min\{s, 1\}}(\Omega)}, \quad \|z-z_n\| \leq Ch^{2\min\{s, 1\}}\|z\|_{H^{1+\min\{s, 1\}}(\Omega)}.
\ee


Given $g\in L^2(\Omega)$ in (\ref{poissoneq}), it is well known that the solution $z \in H^{1+\alpha}(\Omega)$ with $\alpha<\frac{\pi}{\omega}$ (see e.g., \cite{Grisvard1985, Grisvard1992,li2022graded}).
Note that $f, \Delta \chi_i\in L^2(\Omega)$ in Poisson equations (\ref{l2part}) and (\ref{eqnnew+}), so it follows $w, \zeta_i \in H^{1+\alpha}(\Omega)$. Note that $\xi_i \in H^{-1}(\Omega)$, but Step 3 in Algorithm \ref{femalg} indicates $\xi_i-\xi_{i,n}=\zeta_i-\zeta_{i,n}$.
Therefore, we have the following error estimates.
\begin{lemma}\label{lemwerrl2}
Given $w_n$ and $\xi_{i,n}$ in Algorithm \ref{femalg}, it follows
\begin{subequations}\label{werrs+}
\begin{align}
&\|w-w_n\|_{H^1(\Omega)} \leq Ch^{\min\{\alpha, 1\}}\|w\|_{H^{1+\min\{\alpha, 1\}}(\Omega)},\\
&\|w-w_n\| \leq Ch^{2\min\{\alpha, 1\}}\|w\|_{H^{1+\min\{\alpha, 1\}}(\Omega)}, \\
& \|\xi_i-\xi_{i,n}\|_{H^1(\Omega)} \leq Ch^{\min\{\alpha, 1\}}\|\zeta_i\|_{H^{1+\min\{\alpha, 1\}}(\Omega)},\\
& \|\xi_i-\xi_{i,n}\|_{H^{-1}(\Omega)} \leq C\|\xi_i-\xi_{i,n}\| \leq Ch^{2\min\{\alpha, 1\}}\|\zeta_i\|_{H^{1+\min\{\alpha, 1\}}(\Omega)}.
\end{align}
\end{subequations}
\end{lemma}

Note that the basis $\{\sigma_i\}_{i=1}^N$ given in \Cref{defs} is not orthogonal if $\omega>\pi$. For analysis convenience, we can apply Schmidt orthogonalization to obtain an orthogonal basis $\{\tilde{\sigma}_i\}_{i=1}^N$, 
\be\label{xibasis}
\bal
\tilde{\sigma}_1 = & \sigma_1,\\
\tilde{\sigma}_2 = & \sigma_2 - \frac{(\nabla\sigma_2, \nabla\tilde{\sigma}_1)}{\|\nabla\tilde{\sigma}_1\|^2}\tilde{\sigma}_1,\\
\tilde{\sigma}_3 = & \sigma_3 - \frac{(\nabla\sigma_3, \nabla\tilde{\sigma}_1)}{\|\nabla\tilde{\sigma}_1\|^2}\tilde{\sigma}_1- \frac{(\nabla\sigma_3, \nabla\tilde{\sigma}_2)}{\|\nabla\tilde{\sigma}_2\|^2}\tilde{\sigma}_2.\\
\eal
\ee
Namely, $(\nabla \tilde \sigma_i, \nabla \tilde \sigma_j) =\delta_{ij}$, where $\delta_{ij}$ is the the Kronecker delta function.
Furthermore, we denote $\{\tilde{\xi}_i\}_{i=1}^N$ by
\be\label{xibasis2}
\bal
\tilde{\xi}_1 = & \xi_1,\\
\tilde{\xi}_2 = & \xi_2 - \frac{(\nabla\sigma_2, \nabla\tilde{\sigma}_1)}{\|\nabla\tilde{\sigma}_1\|^2}\tilde{\xi}_1,\\
\tilde{\xi}_3 = & \xi_3 - \frac{(\nabla\sigma_3, \nabla\tilde{\sigma}_1)}{\|\nabla\tilde{\sigma}_1\|^2}\tilde{\xi}_1- \frac{(\nabla\sigma_3, \nabla\tilde{\sigma}_2)}{\|\nabla\tilde{\sigma}_2\|^2}\tilde{\xi}_2.\\
\eal
\ee
It can be verified that
\be\label{tsige}
-\Delta \tilde{\sigma}_i = \tilde{\xi}_i  \text{ in } \Omega, \qquad \tilde{\sigma}_i = 0 \text{ on } \partial \Omega,
\ee
and its weak formulation is to find $\tilde{\sigma}_i\in H_0^1(\Omega)$ such that $\forall \phi\in H_0^1(\Omega)$,
\be\label{tsigew}
A(\tilde{\sigma}_i, \phi)=\langle\tilde{\xi}_i, \phi\rangle.
\ee
With the new basis $\{\tilde\sigma_i\}_{i=1}^N$, the third Poisson problem in (\ref{eqnnew+}) can be equivalently written as
\begin{eqnarray}\label{eqnnew++}
\left\{\begin{array}{ll}
-\Delta u=v - \sum_{i=1}^N \tilde{c}_i \tilde{\sigma}_i \quad {\rm{in}} \ \Omega,\\
\hspace{0.6cm}u=0 \quad {\rm{on}} \ \pa\Omega,
\end{array}\right.
\end{eqnarray}
where the coefficients
\be\label{scoeft}
\tilde{c}_i=\frac{\langle v,\tilde{\xi}_i\rangle}{\langle\tilde{\sigma}_i,\tilde{\xi}_i\rangle}, \quad i=1,\ldots, N,
\ee
or equivalently,
\be\label{scoeft1}
\tilde{c}_i=\frac{(\nabla v,\nabla\tilde{\sigma}_i)}{\|\nabla\tilde{\sigma}_i\|^2}, \quad i=1,\ldots, N.
\ee
Correspondingly, the weak formulation (\ref{weaknew}c) becomes
\be\label{weaknew+}
A(u, \tau) = \left(v - \sum_{i=1}^N \tilde{c}_i \tilde{\sigma}_i, \tau \right).
\ee

Similarly, we apply the Schmidt orthogonalization to obtain an orthogonal basis $\{\tilde{\sigma}_{i,n}\}_{i=1}^N$,
\be\label{xibasish}
\bal
\tilde{\sigma}_{1,n} = & \sigma_{1,n},\\
\tilde{\sigma}_{2,n} = & \sigma_{2,n} - \frac{(\nabla\sigma_{2,n}, \nabla\tilde{\sigma}_{1,n})}{\|\nabla\tilde{\sigma}_{1,n}\|^2}\tilde{\sigma}_{1,n},\\
\tilde{\sigma}_{3,n} = & \sigma_{3,n} - \frac{(\nabla\sigma_{3,n}, \nabla\tilde{\sigma}_{1,n})}{\|\nabla\tilde{\sigma}_{1,n}\|^2}\tilde{\sigma}_{1,n}- \frac{(\nabla\sigma_{3,n}, \nabla\tilde{\sigma}_{2,n})}{\|\nabla\tilde{\sigma}_{2,n}\|^2}\tilde{\sigma}_{2,n}.\\
\eal
\ee
Namely, $(\nabla \tilde \sigma_{i,n}, \nabla \tilde \sigma_{j,n}) =\delta_{ij}$, where $\delta_{ij}$ is the the Kronecker delta function.
Similarly, we take $\{\tilde{\xi}_i\}_{i=1}^N$,
\be\label{xibasis2h}
\bal
\tilde{\xi}_{1,n} = & \xi_{1,n},\\
\tilde{\xi}_{2,n} = & \xi_{2,n} - \frac{(\nabla\sigma_{2,n}, \nabla\tilde{\sigma}_{1,n})}{\|\nabla\tilde{\sigma}_{1,n}\|^2}\tilde{\xi}_{1,n},\\
\tilde{\xi}_{3,n} = & \xi_{3,n} - \frac{(\nabla\sigma_{3,n}, \nabla\tilde{\sigma}_{1,n})}{\|\nabla\tilde{\sigma}_{1,n}\|^2}\tilde{\xi}_{1,n}- \frac{(\nabla\sigma_{3,n}, \nabla\tilde{\sigma}_{2,n})}{\|\nabla\tilde{\sigma}_{2,n}\|^2}\tilde{\xi}_{2,n}.\\
\eal
\ee
For the orthogonal basis $\{\tilde \sigma_{i,n}\}_{i=1}^N$,
\be\label{sigweak2}
A(\tilde{\sigma}_{i,n}, \phi) = \langle\tilde{\xi}_{i,n}, \phi\rangle,
\ee
and the last two steps of Algorithm \ref{femalg} can be modified as
\begin{itemize}
\item{Step 5'.} Find the coefficient $\tilde{c}_{i,n}\in \mathbb R$,
\be\label{scoefh2}
\tilde{c}_{i,n} = \frac{\langle v_n, \tilde \xi_{i,n} \rangle}{\langle \tilde\sigma_{i,n}, \tilde\xi_{i,n} \rangle}, \quad i=1, \ldots, N.
\ee
\item{Step 6'.} Find the finite element solution $u_n\in S_n$ of the Poisson equation
\be\label{femu2}
A(u_n, \tau)=\left(v_n - \sum_{i=1}^N \tilde{c}_{i,n} \tilde{\sigma}_{i,n}, \tau \right), \qquad \forall \tau\in S_n.
\ee
\end{itemize}


To show the error estimates, we prepare the following results.

\begin{lem}\label{rfrac}
(i) Assume that $0 \leq s \leq 1$. Then for $\phi\in H_0^s(\Omega)\subset H_0^1(\Omega)$ it follows that $r^{-s}\phi \in L^2(\Omega)$ and
\be\label{ralpha}
\|r^{-s}\phi\|\leq C \|\phi\|_{H^{s}(\Omega)}\leq C \|\phi\|_{H^{1}(\Omega)}.\\
\ee
(ii) If $\gamma \in [0,1)$, $s' \leq 1+\gamma$, and $\phi \in H_0^{s'-\gamma}(\Omega)\subset H_0^1(\Omega)$, then we have $r^{-s'}\phi \in H^{-\gamma}(\Omega)$ and
\be\label{ralphab}
\|r^{-s'}\phi\|_{H^{-\gamma}(\Omega)}\leq C \| r^{-s'+\gamma}\phi \|. 
\ee
\end{lem}
\begin{proof}
(i) The estimate (\ref{ralpha}) follows from \cite[Theorem 1.2.15]{Grisvard1992}.

(ii) Since $s'-\gamma \leq 1$, then we have $r^{-s'+\gamma}\phi \in L^2(\Omega)$ by (i) and it holds $\| r^{-s'+\gamma}\phi \| \leq C \|\phi\|_{H^{s'-\gamma}(\Omega)}$.
For (\ref{ralphab}), we have
\begin{align*}
\|r^{-s'}\phi\|_{H^{-\gamma}(\Omega)} := & \sup_{\psi\in H_0^{\gamma}(\Omega)} \frac{\langle r^{-s'}\phi,\psi \rangle}{\|\psi\|_{H^{\gamma}(\Omega)}}
=  \sup_{\psi\in H_0^{\gamma}(\Omega)} \frac{\langle r^{-s'+\gamma}\phi,r^{-\gamma}\psi \rangle}{\|\psi\|_{H^{\gamma}(\Omega)}}\\
\leq &  \sup_{\psi\in H_0^{\gamma}(\Omega)} \frac{\| r^{-s'+\gamma}\phi \| \|r^{-\gamma}\psi \|}{\|\psi\|_{H^{\gamma}(\Omega)}}
\leq C \| r^{-s'+\gamma}\phi \|,
\end{align*}
where (\ref{ralpha}) is used for $\psi$ in the last inequality.
\end{proof}

%

Next, we introduce some regularity results for a general Poisson problem \eqref{poissoneq}.
\begin{lem}\label{zreglem}
For $g\in H^{\min\{\alpha-1,s\}}(\Omega)$ for any $s \in (-1,0]$ and $\alpha \in (\frac{1}{2}, \frac{\pi}{\omega})$, then (\ref{poissoneq}) admits a unique solution $z\in H^{\min\{\alpha+1, s+2\}}(\Omega)$ and it holds
\be\label{eregpoi}
\|z\|_{H^{\min\{\alpha+1,s+2\}}(\Omega)} \leq C\|g\|_{H^{\min\{\alpha-1,s\}}(\Omega)}.
\ee
\end{lem}
\begin{proof}
The proof follows from \cite[Theorem 3.1]{nazarov1994}.
\end{proof}

\begin{lem}\label{zregdel}
For $\beta_i \in (-1,1-\frac{i\pi}{\omega})$, $i=1,\ldots, N$, and $\alpha \in (\frac{1}{2}\lfloor \frac{2\pi}{\omega}\rfloor, \frac{\pi}{\omega})$ with $\omega \in (\frac{\pi}{2}, 2\pi)$, if $\phi\in H_0^1(\Omega)$ and $g=r^{2\min\{\alpha-1,\beta_i\}}\phi$, then (\ref{poissoneq}) admits a unique solution $z\in H^{\min\{2+\beta_i,1+\alpha\}}(\Omega)$ and holds the estimates
\be\label{zesti}
\|z\|_{H^{\min\{2+\beta_i,1+\alpha\}}(\Omega)} \leq C \|r^{\min\{\alpha-1,\beta_i\}}\phi\|.
\ee
Here, $\lfloor \cdot \rfloor$ represents the floor function.
\end{lem}
\begin{proof}
(1) If $\omega \in \left(\frac{\pi}{2}, \pi\right)$, namely $\frac{\pi}{\omega} \in (1,2)$, then $-1<\beta_N=\beta_1<1-\frac{\pi}{\omega} <0$ and $\alpha-1 > \frac{1}{2}\lfloor \frac{2\pi}{\omega}\rfloor -1 \geq 0$. Consequently, it holds $\min\{\alpha-1,\beta_i\}=\beta_i \in (-1,0)$.\\
(2) If $\omega \in (\pi, 2\pi)$, it follows $\alpha \in \left(\frac{1}{2}, \frac{\pi}{\omega}\right) \subset \left(\frac{1}{2},1\right)$, implying $\alpha-1 \in \left(-\frac{1}{2},0\right)$. This, together with the assumption on $\beta_i$, implies $-1<\min\{\alpha-1,\beta_i\}\leq \alpha-1 <0$.\\
Combining (1) and (2), we conclude that for $\omega \in \left(\frac{\pi}{2}, 2\pi\right) \setminus {\pi}$,
\be\label{alphabeta}
\min\{\alpha-1,\beta_i\} \in (-1,0).
\ee

For $\forall \phi \in H_0^1(\Omega)$, taking $s=-\min\{\alpha-1, \beta_i\} \in (0,1)$ in \Cref{rfrac}(i) yields $r^{\min\{\alpha-1, \beta_i\}}\phi = r^{-s}\phi \in L^2(\Omega)$ and
\be\label{rphiL2}
\|r^{\min\{\alpha-1, \beta_i\}}\phi\| = \|r^{-s}\phi\|  \leq C\|\phi\|_{H^s(\Omega)} \leq C\|\phi\|_{H^1(\Omega)}.
\ee
By taking $s'=-2\min\{\alpha-1,\beta_i\}$ and $\gamma=-\min\{\alpha-1,\beta_i\}$ in Lemma \ref{rfrac}(ii), it follows $g = r^{-s'}\phi \in H^{-\gamma}(\Omega) = H^{\min\{\alpha-1, \beta_i\}}(\Omega)$ and
\be\label{resti0}
\|g\|_{H^{\min\{\alpha-1, \beta_i\}}(\Omega)}=\|r^{-s'}\phi\|_{H^{-\gamma}(\Omega)} \leq C\|r^{-s'+\gamma}\phi\| = C \|r^{\min\{\alpha-1, \beta_i\}}\phi\|.
\ee
\eqref{rphiL2} and \eqref{resti0} imply that for $\forall \phi \in H_0^1(\Omega)$,
\be
\|g\|_{H^{\min\{\alpha-1, \beta_i\}}(\Omega)} \leq C \|\phi\|_{H^1(\Omega)}.
\ee

By Lemma \ref{zreglem}, the Poisson problem (\ref{poissoneq}) admits a unique solution $z\in H^{\min\{2+\beta_i,1+\alpha\}}(\Omega)$ and
\be\label{eregpoi2}
\|z\|_{H^{\min\{2+\beta_i,1+\alpha\}}(\Omega)} \leq C\|g\|_{H^{\min\{\alpha-1, \beta_i\}}(\Omega)},
\ee
which, combined with (\ref{resti0}), yields the estimate (\ref{zesti}).
\end{proof}

By (\ref{ssol}) in Definition \ref{defs}, we have $\xi_i \in H^{\beta_i}\subset H^{-1}(\Omega)$, where $-1<\beta_i<1-\frac{i\pi}{\omega}$, $i=1,\ldots,N$ satisfying $\beta_1>\ldots>\beta_N$. Applying \Cref{zreglem} to the Poisson problem (\ref{H1part}), it follows $\sigma_i \in H^{\min\{2+\beta_i, 1+\alpha\}}$, which is further specified in Table \ref{sigreg}.
\begin{table}[!htbp]\tabcolsep0.03in
\caption{The regularity of $\sigma_i$ in different cases. ($--$ means no such term.)}
\begin{tabular}[c]{|c|c|c|c|c|}
\hline
$\omega$ & $(0,\frac{\pi}{2}]$ & $(\frac{\pi}{2},\pi)$ & $(\pi,\frac{3\pi}{2}] $ & $(\frac{3\pi}{2}, 2\pi)$ \\
\hline
$\sigma_1$ & $--$ & $H^{2+\beta_1}(\Omega)$ & $H^{1+\alpha}(\Omega)$ & $H^{1+\alpha}(\Omega)$ \\
\hline
$\sigma_2$ & $--$ & $--$ & $H^{2+\beta_2}(\Omega)$ & $H^{1+\alpha}(\Omega)$ \\
\hline
$\sigma_3$ & $--$ & $--$ & $--$ & $H^{2+\beta_3}(\Omega)$ \\
\hline
\end{tabular}\label{sigreg}
\end{table}

Then for the finite element solution $\sigma_{i,n}$ in (\ref{fems}), we have the following result.
\begin{lem}\label{lemsigerr}
For $\sigma_{i,n}$ in Algorithm \ref{femalg}, we have for $1\leq i\leq N$,
\begin{subequations}\label{smerrs}
\begin{align}
& \|\sigma_i-\sigma_{i,n}\|_{H^1(\Omega)} \leq Ch^{\min\{1+\beta_i, \alpha\}},\\
& \|\sigma_i-\sigma_{i,n}\| \leq Ch^{\min\{1+\beta_i+\min\{\alpha,1\}, 2\alpha\}},\\
&\|r^{\min\{\alpha-1, \beta_i\}}(\sigma_j-\sigma_{j,n})\| \leq Ch^{\min\{1+\beta_i, \alpha\}+\min\{1+\beta_j, \alpha\}},
\end{align}
\end{subequations}
where $1\leq j \leq N$.
\end{lem}
\begin{proof}
The difference of weak formulation of (\ref{H1part}) and (\ref{fems}) gives
\be\label{siggal}
A(\sigma_i-\sigma_{i,n}, \phi) = (\xi_i-\xi_{i,n}, \phi).
\ee
Let $\sigma_{i,I}\in S_n$ be the nodal interpolation of $\sigma_i$. Set $\epsilon_i = \sigma_{i,I} - \sigma_i, \ e_i = \sigma_{i,I} - \sigma_{i,n}$ and take $\phi=e_i$ in the equation above, we have
\begin{align*}
A(e_i,e_i) = A(\epsilon_i,e_i) + (\xi_i-\xi_{i,n},e_i),
\end{align*}
which implies
\begin{align*}
\|e_i\|_{H^1(\Omega)} \leq \|\epsilon_i\|_{H^1(\Omega)} + \|\xi_i-\xi_{i,n}\|_{H^{-1}(\Omega)}.
\end{align*}
Using the triangle inequality, it follows
\begin{align*}
\|\sigma_i-\sigma_{i,n}\|_{H^1(\Omega)} \leq & \|e_i\|_{H^1(\Omega)}+\|\epsilon_i\|_{H^1(\Omega)} \leq C\left(\|\epsilon_i\|_{H^1(\Omega)} + \|\xi_i-\xi_{i,n}\|_{H^{-1}(\Omega)}\right) \leq Ch^{\min\{1+\beta_i, \alpha\}},
\end{align*}
where we have used the projection error (\ref{interr}) and (\ref{werrs+}d).

To obtain the error in $L^2$ norm, we consider the Poisson problem (\ref{poissoneq}).
By the Aubin-Nitsche Lemma in \cite[Theorem 3.2.4]{Ciarlet78}, we have
\be\label{wel2}
\bal
\|\sigma_i-\sigma_{i,n}\| \leq  C \|\sigma_i-\sigma_{i,n}\|_{H^1(\Omega)} \sup_{g \in L^2(\Omega)} \left( \frac{\inf_{\psi \in S_n}\|z-\psi\|_{H^1(\Omega)}}{\|g\|}  \right).
\eal
\ee
By the regularity (\ref{eregpoi}a), we have
\be\label{projl2}
\bal
\inf_{\psi \in S_n}\|z-\psi\|_{H^1(\Omega)} \leq \|z-z_I\|_{H^1(\Omega)} \leq Ch^{\min\{\alpha,1\}}\|z\|_{H^{1+\min\{\alpha,1\}}(\Omega)}\leq Ch^{\min\{\alpha,1\}}\|g\|.
\eal
\ee
Plugging (\ref{projl2}) and (\ref{smerrs}a) into (\ref{wel2}) gives the estimate (\ref{smerrs}b).

We take $g=r^{2\min\{\alpha-1, \beta_i\}}(\sigma_j-\sigma_{j,n})$, $i=1,\ldots, N$ in (\ref{poissoneq}), since $\sigma_j-\sigma_{j,n} \in H_0^1(\Omega)$, so Lemma \ref{zregdel} indicates that $z \in H^{\min\{2+\beta_i, 1+\alpha\}}(\Omega)$. By (\ref{interr}), we have the interpolation error
\be\label{interrz2}
\| z - z_I \|_{H^1(\Omega)} \leq Ch^{\min\{1+\beta_i,\alpha\}}\|z\|_{H^{\min\{2+\beta_i, 1+\alpha\}}(\Omega)}.
\ee
The weak formulation of (\ref{poissoneq}) with given $g$ is find to $z\in H_0^1(\Omega)$ such that
$$
\langle r^{2\min\{\alpha-1, \beta_i\}}(\sigma_j-\sigma_{j,n}), \psi \rangle = A(z,\psi), \quad \forall \psi \in H_0^1(\Omega).
$$
Set $\psi=\sigma_j-\sigma_{j,n}$ and subtract (\ref{siggal}) with $\phi=z_I$ from the equation above, it follows
\begin{align*}
\|r^{\min\{\alpha-1, \beta_i\}}(\sigma_j-\sigma_{j,n})\|^2 = & A(\sigma_j-\sigma_{j,n}, z-z_I)+(\xi_j-\xi_{j,n}, z_I)\\
= & A(\sigma_j-\sigma_{j,n}, z-z_I)+(\xi_j-\xi_{j,n}, z_I-z)+(\xi_j-\xi_{j,n}, z)\\
\leq & \|\sigma_j-\sigma_{j,n}\|_{H^1(\Omega)}\|z-z_I\|_{H^1(\Omega)} + \|\xi_j-\xi_{j,n}\|_{H^{-1}(\Omega)}(\|z_I-z\|_{H^1(\Omega)}+\|z\|_{H^1(\Omega)}).
\end{align*}
By the estimates in (\ref{werrs+}d), (\ref{smerrs}a), (\ref{interrz2}), and the regularity result in \Cref{zregdel}, it holds
\begin{align*}
\|r^{\min\{\alpha-1, \beta_i\}}(\sigma_j-\sigma_{j,n})\|^2
\leq & Ch^{\min\{1+\beta_i, \alpha\}+\min\{1+\beta_j, \alpha\}} \|z\|_{H^{\min\{2+\beta_i, 1+\alpha\}}(\Omega)} \\
\leq & Ch^{\min\{1+\beta_i, \alpha\}+\min\{1+\beta_j, \alpha\}} \|r^{\min\{\alpha-1, \beta_i\}}(\sigma_j-\sigma_{j,n})\|,
\end{align*}
which gives the error estimate (\ref{smerrs}c).
\end{proof}


Lemma \ref{lemwerrl2} and Lemma \ref{lemsigerr} imply that $\|\xi_{i,n}\|_{H^{-1}(\Omega)}, \|\sigma_{i,n}\|$ and $\|\nabla \sigma_{i,n}\|$, $i=1,\ldots, N$ are uniformly bounded when $h\leq h_0$ for some threshold $h_0$.

\begin{lem}\label{sigest}
For the basis $\tilde \sigma_i$ and the corresponding finite element solution $\tilde \sigma_{i,n}$, we have
\be\label{smerrs2}
\|\tilde\sigma_i-\tilde\sigma_{i,n}\|_{H^1(\Omega)} \leq Ch^{\min\{1+\beta_i, \alpha\}}, \quad i=1,\ldots,N.
\ee
\end{lem}
\begin{proof}
By Lemma \ref{lemsigerr}, it is obvious that
$$
\|\tilde\sigma_1-\tilde\sigma_{1,n}\|_{H^1(\Omega)}=\|\sigma_1-\sigma_{1,n}\|_{H^1(\Omega)} \leq Ch^{\min\{1+\beta_1, \alpha\}}.
$$
We assume that the conclusion holds for $i\leq j-1$,
\be\label{sigassu}
\|\tilde\sigma_i-\tilde\sigma_{i,n}\|_{H^1(\Omega)} \leq Ch^{\min\{1+\beta_i, \alpha\}}.
\ee
A quick calculation gives that
\begin{align*}
\nabla \tilde\sigma_j-\nabla \tilde\sigma_{j,n} = & \nabla\sigma_j - \nabla\sigma_{j,n} - \sum_{i=1}^{j-1} \left(\frac{(\nabla\sigma_j, \nabla\tilde{\sigma}_i)}{\|\nabla\tilde{\sigma}_i\|^2}\nabla\tilde{\sigma}_i-\frac{(\nabla\sigma_{j,n}, \nabla\tilde{\sigma}_{i,n})}{\|\nabla\tilde{\sigma}_{i,n}\|^2}\nabla\tilde{\sigma}_{i,n} \right) \\
= & \left(\nabla\sigma_j - \nabla\sigma_{j,n}\right) - \sum_{i=1}^{j-1}\frac{ \left( (\nabla\sigma_j, \nabla\tilde{\sigma}_i)\|\nabla\tilde{\sigma}_{i,n}\|^2\nabla\tilde{\sigma}_i - (\nabla\sigma_{j,n}, \nabla\tilde{\sigma}_{i,n})\|\nabla\tilde{\sigma}_i\|^2\nabla\tilde{\sigma}_{i,n}   \right) }{\|\nabla\tilde{\sigma}_i\|^2\|\nabla\tilde{\sigma}_{i,n}\|^2}.
\end{align*}
We then have
\begin{align*}
\|\nabla \tilde\sigma_j-\nabla \tilde\sigma_{j,n}\| \leq \|\nabla\sigma_j - \nabla\sigma_{j,n}\| + \sum_{i=1}^{j-1}\frac{ \left\| (\nabla\sigma_j, \nabla\tilde{\sigma}_i)\|\nabla\tilde{\sigma}_{i,n}\|^2\nabla\tilde{\sigma}_i - (\nabla\sigma_{j,n}, \nabla\tilde{\sigma}_{i,n})\|\nabla\tilde{\sigma}_i\|^2\nabla\tilde{\sigma}_{i,n}   \right\| }{\|\nabla\tilde{\sigma}_i\|^2\|\nabla\tilde{\sigma}_{i,n}\|^2}.
\end{align*}
We know that $\nabla\tilde{\sigma}_i$ obtained through (\ref{xibasis}) depend only on $\Omega$. Therefore, we have
\be\label{sigbdd}
0<\gamma_1\leq \|\nabla\tilde{\sigma}_i\| \leq \gamma_2, \quad i=1,\ldots, N,
\ee
where $\gamma_1=\min_{1\leq i\leq N} \{\|\nabla\tilde{\sigma}_i\|\}$, and $\gamma_2=\max_{1\leq i\leq N} \{\|\nabla\tilde{\sigma}_i\|\}$. Let $h\leq h_0 \leq \min\left\{1, \left(\frac{\gamma_1}{2C}\right)^{\frac{1}{\min\{1+\beta_i, \alpha\}}}\right\}$, $i=1,\ldots, j-1$ in (\ref{smerrs2}), it follows that
\be\label{signbdd}
\frac{1}{2}\gamma_1\leq \|\nabla\tilde{\sigma}_{i,n}\| \leq \gamma_2-\frac{1}{2}\gamma_1, \quad i=1,\ldots, j-1.
\ee
(\ref{sigbdd}) and (\ref{signbdd}) implies
\be\label{xiubb}
\frac{1}{\|\nabla\tilde{\sigma}_i\|^2\|\nabla\tilde{\sigma}_{i,n}\|^2} \leq C,
\ee
where $C$ is a constant.
By Lemma \ref{lemsigerr}, it holds
\be\label{sigerr001}
\|\nabla\sigma_j - \nabla\sigma_{j,n}\|\leq\|\sigma_j-\sigma_{j,n}\|_{H^1(\Omega)} \leq Ch^{\min\{1+\beta_j, \alpha\}}.
\ee

To this end, we will get an error estimate for $\left\| (\nabla\sigma_j, \nabla\tilde{\sigma}_i)\|\nabla\tilde{\sigma}_{i,n}\|^2\nabla\tilde{\sigma}_i - (\nabla\sigma_{j,n}, \nabla\tilde{\sigma}_{i,n})\|\nabla\tilde{\sigma}_i\|^2\nabla\tilde{\sigma}_{i,n} \right\|$. Note that
\begin{align*}
& (\nabla\sigma_j, \nabla\tilde{\sigma}_i)\|\nabla\tilde{\sigma}_{i,n}\|^2\nabla\tilde{\sigma}_i - (\nabla\sigma_{j,n}, \nabla\tilde{\sigma}_{i,n})\|\nabla\tilde{\sigma}_i\|^2\nabla\tilde{\sigma}_{i,n} \\
& = (\nabla\sigma_j-\nabla\sigma_{j,n}, \nabla\tilde{\sigma}_i)\|\nabla\tilde{\sigma}_{i,n}\|^2\nabla\tilde{\sigma}_i + (\nabla\sigma_{j,n}, \nabla\tilde{\sigma}_i-\nabla\tilde{\sigma}_{i,n})\|\nabla\tilde{\sigma}_{i,n}\|^2\nabla\tilde{\sigma}_i\\
& + (\nabla\sigma_{j,n}, \nabla\tilde{\sigma}_{i,n})(\|\nabla\tilde{\sigma}_{i,n}\|^2-\|\nabla\tilde{\sigma}_{i}\|^2)\nabla\tilde{\sigma}_i  + (\nabla\sigma_{j,n}, \nabla\tilde{\sigma}_{i,n})\|\nabla\tilde{\sigma}_{i}\|^2(\nabla\tilde{\sigma}_i-\nabla\tilde{\sigma}_{i,n}) \\
& := T_1+T_2+T_3+T_4.
\end{align*}
By Lemma \ref{lemsigerr} again, we have
$$
\|T_1\| \leq \|\nabla\tilde{\sigma}_{i,n}\|^2\|\nabla\tilde{\sigma}_i\|^2\|\nabla\sigma_j-\nabla\sigma_{j,n}\| \leq Ch^{\min\{1+\beta_j, \alpha \}}.
$$
By assumption (\ref{sigassu}), we have
$$
\|T_2\| \leq \|\nabla\sigma_{j,n}\|\|\nabla\tilde{\sigma}_{i,n}\|^2\|\nabla\tilde{\sigma}_i\| \|\nabla\tilde{\sigma}_i-\nabla\tilde{\sigma}_{i,n}\| \leq Ch^{\min\{1+\beta_i, \alpha \}}.
$$
$$
\|T_3\|\leq \|\nabla\sigma_{j,n}\| \|\nabla\tilde{\sigma}_{i,n}\|\|\nabla\tilde{\sigma}_i \| \left( \|\nabla\tilde{\sigma}_{i}\|+\|\nabla\tilde{\sigma}_{i,n}\| \right)\left| \|\nabla\tilde{\sigma}_{i}\|-\|\nabla\tilde{\sigma}_{i,n}\| \right|\leq Ch^{\min\{1+\beta_i, \alpha \}},
$$
where we used the inequality $\left| \|\nabla\tilde{\sigma}_{i}\|-\|\nabla\tilde{\sigma}_{i,n}\| \right| \leq \|\nabla\tilde{\sigma}_{i}-\nabla\tilde{\sigma}_{i,n}\| $.

Last, we have
$$
\|T_4\| \leq \|\nabla\sigma_{j,n}\|\| \nabla\tilde{\sigma}_{i,n}\|\|\nabla\tilde{\sigma}_{i}\|^2\|\nabla\tilde{\sigma}_i - \nabla\tilde{\sigma}_{i,n}\|\leq Ch^{\min\{1+\beta_i, \alpha \}}.
$$
Thus, we have
\be\label{sigerr002}
\left\| (\nabla\sigma_j, \nabla\tilde{\sigma}_i)\|\nabla\tilde{\sigma}_{i,n}\|^2\nabla\tilde{\sigma}_i - (\nabla\sigma_{j,n}, \nabla\tilde{\sigma}_{i,n})\|\nabla\tilde{\sigma}_i\|^2\nabla\tilde{\sigma}_{i,n} \right\|\leq\sum_{l=1}^4\|T_l\| \leq Ch^{\min\{1+\beta_i, \alpha\}}.
\ee
Note that $\beta_i>\beta_j$, thus the combination of (\ref{sigerr001}) and (\ref{sigerr002}) gives
\be\label{sigassuj}
\|\tilde\sigma_j-\tilde\sigma_{j,n}\|_{H^1(\Omega)} \leq Ch^{\min\{1+\beta_j, \alpha\}}.
\ee
The method of induction leads to the conclusion.
\end{proof}


\begin{lemma}\label{rxibddlem}
For $i=1,\ldots, N$, it holds $r^{-\min\{\alpha-1, \beta_i\}}\tilde{\xi}_i \in L^2(\Omega)$ and
\be\label{rxibdd}
\|r^{-\min\{\alpha-1, \beta_i\}}\tilde{\xi}_i\| \leq C,
\ee
where $C$ depends on $\beta_i$ and $\Omega$.
\end{lemma}
\begin{proof}
By (\ref{ssol}), for $k \leq i$,
$$
r^{-\min\{\alpha-1, \beta_i\}} \xi_k= \eta(r; \tau, R)r^{-\min\{\alpha-1, \beta_i\}-\frac{k\pi}{\omega}}\sin\left(\frac{k\pi}{\omega}\theta\right)+ r^{-\min\{\alpha-1, \beta_i\}}\zeta_k:=T_{11}+T_{12}.
$$
Since $\beta_i \leq \beta_k <1-\frac{k\pi}{\omega}$, it follows
$$
-\min\{\alpha-1, \beta_i\}-\frac{k\pi}{\omega} > \beta_i -\min\{\alpha-1, \beta_i\}- 1 \geq -1.
$$
Therefore, $T_{11} \in L^2(\Omega)$, namely,
\be\label{t11bdd}
\|T_{11}\|\leq C.
\ee
For $T_{12}$, it follows
\be\label{t12bdd}
\|T_{12}\| \leq \|r^{-\min\{\alpha-1, \beta_i\}}\|_{L^\infty(\Omega)}\|\zeta_k\| \leq C\|\zeta_k\|_{H^1(\Omega)}.
\ee
(\ref{t11bdd}) and (\ref{t12bdd}) imply that
\be\label{rxibdd-}
\|r^{-\min\{\alpha-1, \beta_i\}}{\xi}_k\| \leq \|T_{11}\|+\|T_{12}\| \leq C.
\ee
By the construction of $\tilde \xi_i$ in (\ref{xibasis2}), we can obtain the estimate (\ref{rxibdd}).
\end{proof}

\begin{lem}
For $i=1,\ldots, N$, the orthogonal functions $\tilde{\sigma}_{i}$ in (\ref{xibasis}), $\tilde{\xi}_{i}$ in (\ref{xibasis2}), and their finite element approximations $\tilde{\sigma}_{i,n}$ in (\ref{xibasish}), $\tilde{\xi}_{i,n}$ in (\ref{xibasis2h}) satisfy
\begin{subequations}\label{tsigl22}
\begin{align}
& \|\tilde\xi_i-\tilde\xi_{i,n}\|_{H^{-1}(\Omega)} \leq Ch^{\min\{1+\beta_i+\min\{\alpha, 1\}, 2\alpha \}},\\
& \|\tilde\sigma_i-\tilde\sigma_{i,n}\| \leq Ch^{\min\{1+\beta_i+\min\{\alpha, 1\}, 2\alpha \}},\\
& \|r^{\min\{\alpha-1, \beta_k\}}(\tilde\sigma_i-\tilde\sigma_{i,n})\| \leq Ch^{\min\{1+\beta_k, \alpha\}+\min\{1+\beta_i, \alpha\}},
\end{align}
\end{subequations}
where $1 \leq k \leq N$.
\end{lem}
\begin{proof}
It is easy to verify that the estimates in (\ref{tsigl22}) hold when $i=1$, and we assume that they also hold for $i\leq j-1$ if $N\geq 2$.
Next, we prove the estimates in (\ref{tsigl22}) hold at $j$.
The proof for (\ref{tsigl22}b) is similar to that for (\ref{tsigl22}a), we will skip its proof.

Using the similar argument as in Lemma \ref{sigest}, we have that $\|\tilde{\xi}_i\|_{H^{-1}(\Omega)}, \|\tilde{\sigma}_i\|$, $i=1,\ldots, N$ are uniformly bounded. When $h\leq h_0$ for some $h_0$, it follows that $\|\tilde{\xi}_{i,n}\|_{H^{-1}(\Omega)}, \|\tilde{\sigma}_{i,n}\|$, $1\leq i\leq j-1$ are also uniformly bounded.

The difference of the $\tilde \xi_j- \tilde \xi_{j,n}$, $j=2,\ldots,N$, gives
\be\label{zdiff}
\bal
\tilde \xi_j- \tilde \xi_{j,n} = & \xi_j - \xi_{j,n} - \sum_{i=1}^{j-1} \left(\frac{(\nabla\sigma_j, \nabla\tilde{\sigma}_i)}{\|\nabla\tilde{\sigma}_i\|^2}\tilde{\xi}_i-\frac{(\nabla\sigma_{j,n}, \nabla\tilde{\sigma}_{i,n})}{\|\nabla\tilde{\sigma}_{i,n}\|^2}\tilde{\xi}_{i,n} \right) \\
= & \left(\xi_j - \xi_{j,n}\right) - \sum_{i=1}^{j-1}\frac{ \left( (\nabla\sigma_j, \nabla\tilde{\sigma}_i)\|\nabla\tilde{\sigma}_{i,n}\|^2\tilde{\xi}_i - (\nabla\sigma_{j,n}, \nabla\tilde{\sigma}_{i,n})\|\nabla\tilde{\sigma}_i\|^2\tilde{\xi}_{i,n}   \right) }{\|\nabla\tilde{\sigma}_i\|^2\|\nabla\tilde{\sigma}_{i,n}\|^2}.
\eal
\ee
By (\ref{xiubb}), $\frac{1}{\|\nabla\tilde{\sigma}_i\|^2\|\nabla\tilde{\sigma}_{i,n}\|^2}$ are uniformly bounded.
We denote by
\begin{align*}
& (\nabla\sigma_j, \nabla\tilde{\sigma}_i)\|\nabla\tilde{\sigma}_{i,n}\|^2\tilde{\xi}_i - (\nabla\sigma_{j,n}, \nabla\tilde{\sigma}_{i,n})\|\nabla\tilde{\sigma}_i\|^2\tilde{\xi}_{i,n} \\
& = (\nabla\sigma_j-\nabla\sigma_{j,n}, \nabla\tilde{\sigma}_i)\|\nabla\tilde{\sigma}_{i,n}\|^2\tilde{\xi}_i
  + (\nabla\sigma_{j,n}, \nabla\tilde{\sigma}_i-\nabla\tilde{\sigma}_{i,n})\|\nabla\tilde{\sigma}_{i,n}\|^2\tilde{\xi}_i\\
& + (\nabla\sigma_{j,n}, \nabla\tilde{\sigma}_{i,n})(\|\nabla\tilde{\sigma}_{i,n}\|^2-\|\nabla\tilde{\sigma}_{i}\|^2)\tilde{\xi}_i
  + (\nabla\sigma_{j,n}, \nabla\tilde{\sigma}_{i,n})\|\nabla\tilde{\sigma}_{i}\|^2(\tilde{\xi}_i-\tilde{\xi}_{i,n}) \\
& := T_1+T_2+T_3+T_4.
\end{align*}

By (\ref{zdiff}), it follows
\begin{align*}
\|\tilde \xi_j- \tilde \xi_{j,n}\|_{H^{-1}(\Omega)} \leq \|\xi_j - \xi_{j,n}\|_{H^{-1}(\Omega)} + \sum_{i=1}^{j-1}\frac{ \left\|(\nabla\sigma_j, \nabla\tilde{\sigma}_i)\|\nabla\tilde{\sigma}_{i,n}\|^2\tilde{\xi}_i - (\nabla\sigma_{j,n}, \nabla\tilde{\sigma}_{i,n})\|\nabla\tilde{\sigma}_i\|^2\tilde{\xi}_{i,n}   \right\|_{H^{-1}(\Omega)} }{\|\nabla\tilde{\sigma}_i\|^2\|\nabla\tilde{\sigma}_{i,n}\|^2}.
\end{align*}
From (\ref{werrs+}d), we have
\be\label{xierrtil}
\|\xi_j - \xi_{j,n}\|_{H^{-1}(\Omega)} \leq Ch^{2\min\{\alpha, 1\}}.
\ee
By taking $\phi=\sigma_j-\sigma_{j,n}\in H_0^1(\Omega)$ in (\ref{tsigew}), we have
\be\label{sigtoxi}
(\nabla\sigma_j-\nabla\sigma_{j,n}, \nabla\tilde{\sigma}_i)=\langle \sigma_j-\sigma_{j,n}, \tilde{\xi}_i \rangle,
\ee
which implies that
$$
T_1=\langle \sigma_j-\sigma_{j,n}, \tilde{\xi}_i \rangle \|\nabla\tilde{\sigma}_{i,n}\|^2\tilde{\xi}_i=\langle r^{\min\{\alpha-1, \beta_i\}}(\sigma_j-\sigma_{j,n}), r^{-\min\{\alpha-1, \beta_i\}}\tilde{\xi}_i \rangle \|\nabla\tilde{\sigma}_{i,n}\|^2\tilde{\xi}_i.
$$
By Lemma \ref{rxibddlem}, we have $r^{-\min\{\alpha-1, \beta_i\}}\tilde{\xi}_i \in L^2(\Omega)$.
Therefore, we have the estimate
\begin{align*}
\|T_1\|_{H^{-1}(\Omega)} \leq &  \|r^{\min\{\alpha-1, \beta_i\}}(\sigma_j-\sigma_{j,n})\| \|r^{-\min\{\alpha-1, \beta_i\}} \tilde{\xi}_i\|\|\nabla\tilde{\sigma}_{i,n}\|^2\|\tilde{\xi}_i\|_{H^{-1}(\Omega)} \\
\leq & Ch^{\min\{1+\beta_i, \alpha\}+\min\{1+\beta_j, \alpha\}} = Ch^{\min\{1+\beta_j+\min\{\alpha,1\}, 2\alpha\}},
\end{align*}
where we have used the estimate (\ref{smerrs}c).

Subtracting equation (\ref{sigweak2}) from equation (\ref{tsigew}) and setting $\phi=\sigma_{j,n}$ yields 
$$
(\nabla\sigma_{j,n}, \nabla\tilde{\sigma}_i-\nabla\tilde{\sigma}_{i,n})=\langle \sigma_{j,n}, \tilde\xi_i - \tilde\xi_{i,n} \rangle.
$$
Thus, we have by the assumption,
$$
\|T_2\|_{H^{-1}(\Omega)} \leq \|\sigma_{j,n}\|\| \tilde\xi_i - \tilde\xi_{i,n}\|_{H^{-1}(\Omega)}\|\nabla\tilde{\sigma}_{i,n}\|^2\|\tilde{\xi}_i\|_{H^{-1}(\Omega)} \leq Ch^{\min\{1+\beta_i+\min\{\alpha, 1\}, 2\alpha \}}.
$$
We have by (\ref{tsigew}) and (\ref{sigweak2}),
$$
\|\nabla\tilde{\sigma}_{i,n}\|^2-\|\nabla\tilde{\sigma}_{i}\|^2=\langle \tilde{\sigma}_{i,n}, \tilde{\xi}_{i,n} \rangle-\langle\tilde{\sigma}_{i}, \tilde{\xi}_{i}\rangle=\langle \tilde{\sigma}_{i,n}-\tilde{\sigma}_{i}, \tilde{\xi}_{i,n} \rangle+\langle\tilde{\sigma}_{i}, \tilde{\xi}_{i,n}-\tilde{\xi}_{i}\rangle:=T_{31}+T_{32}.
$$
By the assumption for (\ref{tsigl22}c), we have
$$
|T_{31}|=|\langle \tilde{\sigma}_{i,n}-\tilde{\sigma}_{i}, \tilde{\xi}_{i,n} \rangle| \leq \|r^{\min\{\alpha-1, \beta_i\}} (\tilde{\sigma}_{i,n}-\tilde{\sigma}_{i})\|\|r^{-\min\{\alpha-1, \beta_i\}} \tilde{\xi}_{i,n}\|\leq Ch^{\min\{1+\beta_i+\min\{\alpha, 1\}, 2\alpha \}}.
$$
For the second term, we have by the assumption for (\ref{tsigl22}a),
$$
|T_{32}|=|\langle\tilde{\sigma}_{i}, \tilde{\xi}_{i,n}-\tilde{\xi}_{i}\rangle|\leq \|\tilde{\sigma}_{i}\|_{H^{1}(\Omega)}\|\tilde{\xi}_{i,n}-\tilde{\xi}_{i}\|_{H^{-1}(\Omega)} \leq Ch^{\min\{1+\beta_i+\min\{\alpha, 1\}, 2\alpha \}}.
$$
The estimates of $|T_{31}|$ and $|T_{32}|$ imply that
$$
\|T_3\|_{H^{-1}(\Omega)} \leq Ch^{\min\{1+\beta_i+\min\{\alpha, 1\}, 2\alpha \}}.
$$
Again by the assumption for (\ref{tsigl22}a), we have
$$
\|T_4\|_{H^{-1}(\Omega)} \leq \|\nabla\sigma_{j,n}\|\|\nabla\tilde{\sigma}_{i,n}\|\|\nabla\tilde{\sigma}_{i}\|^2\|\tilde{\xi}_i-\tilde{\xi}_{i,n}\|_{H^{-1}(\Omega)} \leq Ch^{\min\{1+\beta_i+\min\{\alpha, 1\}, 2\alpha \}}.
$$
Note that $\beta_i>\beta_j$, we have
\be\label{subtl}
\sum_{l=1}^4 \|T_l\|_{H^{-1}(\Omega)} \leq Ch^{\min\{1+\beta_j+\min\{\alpha, 1\}, 2\alpha \}}.
\ee
The combination of (\ref{xierrtil}) and (\ref{subtl}) indicate that (\ref{tsigl22}a) holds at $j$, so that the method of induction state that (\ref{tsigl22}a) holds for $i=1,\ldots, N$.

Next, we prove the estimate (\ref{tsigl22}c) holds at $j$. For $j=2,\ldots,N$, we have
\be\label{zdiff+}
\bal
& r^{\min\{\alpha-1, \beta_k\}}(\tilde \sigma_j- \tilde \sigma_{j,n}) =  r^{\min\{\alpha-1, \beta_k\}}(\sigma_j - \sigma_{j,n}) \\
& - r^{\min\{\alpha-1, \beta_k\}}\sum_{i=1}^{j-1} \left(\frac{(\nabla\sigma_j, \nabla\tilde{\sigma}_i)}{\|\nabla\tilde{\sigma}_i\|^2}\tilde{\sigma}_i-\frac{(\nabla\sigma_{j,n}, \nabla\tilde{\sigma}_{i,n})}{\|\nabla\tilde{\sigma}_{i,n}\|^2}\tilde{\sigma}_{i,n} \right) \\
& = r^{\min\{\alpha-1, \beta_k\}}\left(\sigma_j - \sigma_{j,n}\right) \\
& - \sum_{i=1}^{j-1}\frac{ \left( (\nabla\sigma_j, \nabla\tilde{\sigma}_i)\|\nabla\tilde{\sigma}_{i,n}\|^2 r^{\min\{\alpha-1, \beta_k\}} \tilde{\sigma}_i - (\nabla\sigma_{j,n}, \nabla\tilde{\sigma}_{i,n})\|\nabla\tilde{\sigma}_i\|^2 r^{\min\{\alpha-1, \beta_k\}}\tilde{\sigma}_{i,n}   \right) }{\|\nabla\tilde{\sigma}_i\|^2\|\nabla\tilde{\sigma}_{i,n}\|^2},
\eal
\ee
where $1\leq k \leq N$.
By (\ref{xiubb}), $\frac{1}{\|\nabla\tilde{\sigma}_i\|^2\|\nabla\tilde{\sigma}_{i,n}\|^2}$ are uniformly bounded.
We denote by
\begin{align*}
& (\nabla\sigma_j, \nabla\tilde{\sigma}_i)\|\nabla\tilde{\sigma}_{i,n}\|^2 r^{\min\{\alpha-1, \beta_k\}}\tilde{\sigma}_i - (\nabla\sigma_{j,n}, \nabla\tilde{\sigma}_{i,n})\|\nabla\tilde{\sigma}_i\|^2 r^{\min\{\alpha-1, \beta_k\}}\tilde{\sigma}_{i,n} \\
& = (\nabla\sigma_j-\nabla\sigma_{j,n}, \nabla\tilde{\sigma}_i)\|\nabla\tilde{\sigma}_{i,n}\|^2r^{\min\{\alpha-1, \beta_k\}}\tilde{\sigma}_i
  + (\nabla\sigma_{j,n}, \nabla\tilde{\sigma}_i-\nabla\tilde{\sigma}_{i,n})\|\nabla\tilde{\sigma}_{i,n}\|^2r^{\min\{\alpha-1, \beta_k\}}\tilde{\sigma}_i\\
& + (\nabla\sigma_{j,n}, \nabla\tilde{\sigma}_{i,n})(\|\nabla\tilde{\sigma}_{i,n}\|^2-\|\nabla\tilde{\sigma}_{i}\|^2)r^{\min\{\alpha-1, \beta_k\}}\tilde{\sigma}_i
  + (\nabla\sigma_{j,n}, \nabla\tilde{\sigma}_{i,n})\|\nabla\tilde{\sigma}_{i}\|^2r^{\min\{\alpha-1, \beta_k\}}(\tilde{\sigma}_i-\tilde{\sigma}_{i,n}) \\
& := K_1+K_2+K_3+K_4.
\end{align*}
By (\ref{zdiff+}), it follows
\bes
\bal
& \|r^{\min\{\alpha-1, \beta_k\}}(\tilde \sigma_j- \tilde \sigma_{j,n})\| \leq \| r^{\min\{\alpha-1, \beta_k\}}\left(\sigma_j - \sigma_{j,n}\right) \| \\
& + \sum_{i=1}^{j-1}\frac{ \left\| (\nabla\sigma_j, \nabla\tilde{\sigma}_i)\|\nabla\tilde{\sigma}_{i,n}\|^2 r^{\min\{\alpha-1, \beta_k\}} \tilde{\sigma}_i - (\nabla\sigma_{j,n}, \nabla\tilde{\sigma}_{i,n})\|\nabla\tilde{\sigma}_i\|^2 r^{\min\{\alpha-1, \beta_k\}}\tilde{\sigma}_{i,n}   \right\| }{\|\nabla\tilde{\sigma}_i\|^2\|\nabla\tilde{\sigma}_{i,n}\|^2},
\eal
\ees
From (\ref{smerrs}b), we have
\be\label{xierrtil+}
\| r^{\min\{\alpha-1, \beta_k\}}\left(\sigma_j - \sigma_{j,n}\right) \| \leq Ch^{\min\{1+\beta_k, \alpha\}+\min\{1+\beta_j, \alpha\}}.
\ee
Similar to the estimate of $T_1$, we have by (\ref{sigtoxi}),
\begin{align*}
K_1= \langle r^{\min\{\alpha-1, \beta_i\}}(\sigma_j-\sigma_{j,n}), r^{-\min\{\alpha-1, \beta_i\}}\tilde{\xi}_i \rangle \|\nabla\tilde{\sigma}_{i,n}\|^2r^{\min\{\alpha-1, \beta_k\}}\tilde{\sigma}_i.
\end{align*}
Since $\tilde{\sigma}_i \in H_0^1(\Omega) \cap H^{\min\{2+\beta_i, 1+\alpha \}}(\Omega)$, so we have $r^{\min\{\alpha-1, \beta_k\}}\tilde{\sigma}_i \in L^2(\Omega)$ by Lemma \ref{rfrac} or from (\ref{rphiL2}). %
Therefore, we have the estimate
\begin{align*}
\|K_1\| \leq &  \|r^{\min\{\alpha-1, \beta_i\}}(\sigma_j-\sigma_{j,n})\| \|r^{-\min\{\alpha-1, \beta_i\}} \tilde{\xi}_i\|\|\nabla\tilde{\sigma}_{i,n}\|^2\|r^{\min\{\alpha-1, \beta_k\}}\tilde{\sigma}_i\| \\
\leq & Ch^{\min\{1+\beta_i, \alpha\}+\min\{1+\beta_j, \alpha\}}=Ch^{\min\{1+\beta_j+\min\{\alpha, 1\}, 2\alpha \}},
\end{align*}
the last equality is due to the fact that $\min\{1+\beta_i, \alpha\}=\alpha$ when $1\leq i<N$.
Similar to the estimate of $T_2$, we have
$$
\|K_2\|_{H^{-1}(\Omega)} \leq \|\sigma_{j,n}\|\| \tilde\xi_i - \tilde\xi_{i,n}\|_{H^{-1}(\Omega)}\|\nabla\tilde{\sigma}_{i,n}\|^2\|r^{\min\{\alpha-1, \beta_k\}}\tilde{\sigma}_i\| \leq Ch^{\min\{1+\beta_i+\min\{\alpha, 1\}, 2\alpha \}}.
$$
The estimates of $|T_{31}|$ and $|T_{32}|$ above also indicate that
$$
\|K_3\| \leq Ch^{\min\{1+\beta_i+\min\{\alpha, 1\}, 2\alpha \}}.
$$
By the assumption for (\ref{tsigl22}c), we have
$$
\|K_4\| \leq \|\nabla\sigma_{j,n}\|\|\nabla\tilde{\sigma}_{i,n}\|\|\nabla\tilde{\sigma}_{i}\|^2\|r^{\min\{\alpha-1, \beta_k\}}(\tilde{\sigma}_i-\tilde{\sigma}_{i,n})\| \leq Ch^{\min\{1+\beta_k, \alpha\}+\min\{1+\beta_i, \alpha\}}.
$$
Note again that $\beta_i>\beta_j$, we have
\be\label{subtl+}
\sum_{l=1}^4 \|K_l\| \leq Ch^{\min\{1+\beta_k, \alpha\}+\min\{1+\beta_j, \alpha\}}.
\ee
The combination of (\ref{xierrtil+}) and (\ref{subtl+}) indicate that (\ref{tsigl22}c) holds at $j$, so that the method of induction state that (\ref{tsigl22}c) holds for $i=1,\ldots, N$.
\end{proof}

Note that $v\in H^{1+\alpha}(\Omega)$, then we have the following estimates for $v_n$ in (\ref{femv}).
\begin{lem}\label{lemv1err}
Let $v_n \in S_n$ be the finite element approximation to (\ref{femv}), and $v$ be the solution to the Poisson equation in the mixed formulation (\ref{eqnnew+}). Then it follows
\begin{subequations}\label{verrs}
\begin{align}
&\|v-v_n\|_{H^1(\Omega)}\leq Ch^{\min\{\alpha, 1\}}, \\
&\|v-v_n\|\leq Ch^{2\min\{\alpha, 1\}}.
\end{align}
\end{subequations}
\end{lem}
\begin{proof}
Subtracting (\ref{weaknew}b) from (\ref{femv}) gives the Galerkin orthogonality
\be\label{galerkin}
A(v-v_n, \psi) = (w-w_n,\psi).
\ee
Let $v_I\in S_n$ be the nodal interpolation of $v$. Set $\epsilon = v_I - u, \ e = v_I - v_n$ and take $\psi=e$ in the equation above, we have
\begin{align*}
A(e,e) = A(\epsilon,e) + (w-w_n,e),
\end{align*}
which implies
\begin{align*}
\|e\|_{H^1(\Omega)} \leq \|\epsilon\|_{H^1(\Omega)} + \|w-w_n\|_{H^{-1}(\Omega)},
\end{align*}
Using the triangle inequality, it follows
\begin{align*}
\|v-v_n\|_{H^1(\Omega)} \leq & \|e\|_{H^1(\Omega)}+\|\epsilon\|_{H^1(\Omega)} \leq C\left(\|\epsilon\|_{H^1(\Omega)} + \|w-w_n\|_{H^{-1}(\Omega)}\right),\\
\leq & C\left(\|\epsilon\|_{H^1(\Omega)} + \|w-w_n\|\right) \leq Ch^{\min\{\alpha, 1\}},
\end{align*}
where we have used the projection error (\ref{interr}) and (\ref{werrs+}b). To obtain the $L^2$ error, we consider the problem (\ref{poissoneq}) with $g=v-v_n$,
then we have
$$
\|v-v_n\|^2=A(v-v_n, z).
$$
Subtract (\ref{galerkin}) from the above equation and set $\psi=z_I$, we have
\be\label{vl2est}
\bal
\|v-v_n\|^2= & A(v-v_n, z-z_I)+(w-w_n,z_I)\\
= & A(v-v_n, z-z_I)+(w-w_n,z_I-z)+(w-w_n,z)\\
\leq & \|v-v_n\|_{H^1(\Omega)}\|z-z_I\|_{H^1(\Omega)} + \|w-w_n\|\|z-z_I\| + \|w-w_n\|\|z\|\\
\leq & Ch^{2\min\{\alpha, 1\}}\|z\|_{H^{1+\min\{\alpha, 1\}}(\Omega)} \leq Ch^{2\min\{\alpha, 1\}}\|v-v_n\|,
\eal
\ee
where in the last inequality we have use the estimates (\ref{interr}), (\ref{werrs+}b), (\ref{verrs}a). By the regularity (\ref{eregpoi}a), we have
\be\label{vvnreg}
\|z\|_{H^{1+\min\{\alpha, 1\}}(\Omega)} \leq C\|v-v_n\|_{H^{\min\{\alpha, 1\}-1}(\Omega)}\leq C\|v-v_n\|.
\ee
(\ref{vl2est}) and (\ref{vvnreg}) give the $L^2$ error estimate (\ref{verrs}b).
\end{proof}

Next, we carry out the error estimate for the finite element approximation $u_n$ in (\ref{femu}).
\begin{theorem}\label{thmuerr}
Let $u_n \in S_n$ be the finite element approximation to (\ref{femu}), and $u$ be the solution to the sixth order problem (\ref{eqn.firstbi}). Then it follows
\be\label{uerr}
\|u-u_n\|_{H^1(\Omega)} \leq C_0h+\sum_{i=1}^NC_ih^{\min\{2(1+\beta_i),1\}}\leq Ch^\gamma,
\ee
where $-1<\beta_i<1-\frac{i\pi}{\omega}$, the convergence rate $\gamma=1$ if $N=0$, and $\gamma = \min\{2(1+\beta_N),1\}$ if $1\leq N \leq 3$, the constants $C$, $C_i$ depend on the coefficients $\tilde c_i$ in (\ref{scoeft}).
\end{theorem}
\begin{proof}
Subtracting (\ref{femu2}) from (\ref{weaknew+}) gives
\be\label{soldiff}
\bal
A(u-u_n, \tau)=&(v-v_n, \tau) - \sum_{i=1}^N \left( \tilde{c}_{i} \tilde{\sigma}_{i}-\tilde{c}_{i,n} \tilde{\sigma}_{i,n}, \tau\right)\\
=&(v-v_n, \tau)+\sum_{i=1}^N \left[\tilde{c}_{i,n}(\tilde{\sigma}_{i,n}-\tilde{\sigma}_{i}, \tau)+(\tilde{c}_{i,n}-\tilde{c}_{i})(\tilde{\sigma}_{i}, \tau)\right].
\eal
\ee
Let $u_I\in S_n$ be the nodal interpolation of $u$. Set $\epsilon = u_I - u, \ e = u_I - u_n$ and take $\tau=e$ in (\ref{soldiff}), we have
\ben
\bal
A(e, e)=A(\epsilon,e) + (v-v_n, e)+\sum_{i=1}^N \left[\tilde{c}_{i,n}(\tilde{\sigma}_{i,n}-\tilde{\sigma}_{i}, e)+(\tilde{c}_{i,n}-\tilde{c}_{i})(\tilde{\sigma}_{i}, e)\right].
\eal
\een
Thus, we have
$$
\|e\|_{H^1(\Omega)} \leq C \left( \|\epsilon\|_{H^1(\Omega)} +  \|v-v_n\|_{H^{-1}(\Omega)} + \sum_{i=1}^{N} \left[ |\tilde{c}_{i,n}|\|\tilde{\sigma}_{i}-\tilde{\sigma}_{i,n}\|_{H^{-1}(\Omega)} + |\tilde{c}_{i}-\tilde{c}_{i,n}|\|\tilde \sigma_{i}\|_{H^{-1}(\Omega)} \right]  \right).
$$
Using the triangle inequality and the inequality above, we have
\be\label{errbdd}
\bal
& \|u-u_n\|_{H^1(\Omega)} \leq \|e\|_{H^1(\Omega)} +  \|\epsilon\|_{H^1(\Omega)}\\
&\leq  C \left( \|\epsilon\|_{H^1(\Omega)} +  \|v-v_n\|_{H^{-1}(\Omega)} + \sum_{i=1}^{N} \left[  |\tilde{c}_{i,n}|\|\tilde{\sigma}_{i}-\tilde{\sigma}_{i,n}\|_{H^{-1}(\Omega)} + |\tilde{c}_{i}-\tilde{c}_{i,n}|\|\tilde \sigma_{i}\|_{H^{-1}(\Omega)} \right]  \right).
\eal
\ee
We shall estimate every term in (\ref{errbdd}). Recall the solution $u\in H^3(\Omega)$. By the interpolation error estimate (\ref{interr}),
\be\label{uprojerr}
\|\epsilon\|_{H^1(\Omega)} = \|u-u_I\|_{H^1(\Omega)} \leq Ch\|u\|_{H^2(\Omega)}.
\ee
Recall that $\frac{\pi}{\omega}>\frac{1}{2}$. Thus, choosing $\alpha=1/2<\frac{\pi}{\omega}$ in (\ref{verrs}b), we have
\be\label{eqn.11}
 \|v-v_n\|_{H^{-1}(\Omega)} \leq \|v-v_n\|\leq Ch.
\ee
By (\ref{tsigl22}b), we have
$$
\|\tilde{\sigma}_{i}-\tilde{\sigma}_{i,n}\|_{H^{-1}(\Omega)} \leq \|\tilde{\sigma}_{i}-\tilde{\sigma}_{i,n}\| \leq Ch^{\min\{1+\beta_i+\min\{\alpha, 1\}, 2\alpha \}}.
$$
To obtain the error estimate for the third term in (\ref{errbdd}), we still need to show that $|\tilde{c}_{i,n}|$ is uniformly bounded. By (\ref{scoefh2}), we have
\be\label{scoefh2+}
|\tilde{c}_{i,n}| = \left|\frac{\langle v_n, \tilde \xi_{i,n} \rangle}{\langle \tilde\sigma_{i,n}, \tilde\xi_{i,n} \rangle}\right|
=\left| \frac{(\nabla\tilde\sigma_{i,n}, \nabla v_n)}{(\nabla\tilde\sigma_{i,n}, \nabla\tilde\sigma_{i,n})} \right| \leq \frac{\| v_n\|_{H^1(\Omega)}\|\tilde\sigma_{i,n}\|_{H^1(\Omega)}}{\|\tilde\sigma_{i,n}\|^2_{H^1(\Omega)}} \leq \frac{\| v_n\|_{H^1(\Omega)}}{\|\tilde\sigma_{i,n}\|_{H^1(\Omega)}},
\ee
where we have used H\"older's inequality.
By the regularity result (\ref{weakreg0}) and the estimate (\ref{verrs}), we have $\|v_n\|_{H^1(\Omega)} \leq C\|f\|$ when $h\leq h_0$ for some $h_0$, which together with (\ref{signbdd}) implies that (\ref{scoefh2+}) is uniformly bounded.


Subtracting (\ref{scoefh2}) from (\ref{scoeft}) or (\ref{scoeft1}) gives
$$
\tilde{c}_{i}-\tilde{c}_{i,n} = \frac{(\nabla v,\nabla(\tilde\sigma_{i}-\tilde\sigma_{i,n}))}{\|\nabla \tilde\sigma_{i}\|^2}+\frac{(\nabla \tilde\sigma_{i,n}, \nabla (v- v_n))}{\|\nabla \tilde\sigma_{i}\|^2} + \frac{\|\nabla \tilde\sigma_{i,n}\|^2-\|\nabla \tilde\sigma_{i}\|^2}{\|\nabla \tilde\sigma_{i}\|^2\|\nabla \tilde\sigma_{i,n}\|^2}(\nabla v_n, \nabla \tilde\sigma_{i,n}):=T_1+T_2+T_3.
$$
By setting $\psi=(\tilde\sigma_{i}-\tilde\sigma_{i,n}) \in H_0^1(\Omega)$ in (\ref{weaknew}), we obtain
$$
(\nabla v,\nabla(\tilde\sigma_{i}-\tilde\sigma_{i,n})) = (w, \tilde\sigma_{i}-\tilde\sigma_{i,n}).
$$
Thus, we have by (\ref{tsigl22})
$$
\|T_1\|\leq \frac{\|w\|}{\|\nabla \tilde\sigma_{i}\|^2} \|\tilde\sigma_{i}-\tilde\sigma_{i,n}\|\leq Ch^{\min\{1+\beta_i+\min\{\alpha, 1\}, 2\alpha \}}.
$$
Subtracting equation (\ref{femv}) from equation (\ref{weaknew}b) and setting $\psi=\tilde\sigma_{i,n}$, we obtain
$$
(\nabla \tilde\sigma_{i,n}, \nabla (v- v_n))=(w-w_n, \nabla \sigma_{i,n}).
$$
Then we have by (\ref{werrs+}b) and taking $\alpha = \frac{1}{2}$,
$$
\|T_2\| \leq \frac{1}{\|\nabla \tilde\sigma_{i}\|}\|w-w_n\| \leq Ch^{2\min\{\alpha, 1\}}=Ch.
$$
Note that
\begin{align*}
\|\nabla \tilde\sigma_{i,n}\|^2-\|\nabla \tilde\sigma_{i}\|^2=& (\nabla \tilde\sigma_{i,n}-\nabla \tilde\sigma_{i}, \nabla \tilde\sigma_{i,n})+(\nabla \tilde\sigma_{i,n}-\nabla \tilde\sigma_{i}, \nabla \tilde\sigma_{i})\\
=& \langle \tilde{\xi}_{i,n}-\tilde{\xi}_{i} , \tilde\sigma_{i,n}\rangle + \langle \tilde \sigma_{i,n}-\tilde\sigma_{i}, \tilde{\xi}_{i} \rangle.
\end{align*}
By (\ref{tsigl22}a), we have
\be\label{eqn.112}
\|\tilde\xi_i-\tilde\xi_{i,n}\|_{H^{-1}(\Omega)} \leq Ch^{\min\{1+\beta_i+\min\{\alpha, 1\}, 2\alpha \}}.
\ee
It is easy to check
$$
|\langle \tilde{\xi}_{i,n}-\tilde{\xi}_{i} , \tilde\sigma_{i,n}\rangle| \leq \|\tilde{\xi}_{i,n}-\tilde{\xi}_{i}\|_{H^{-1}(\Omega)}\|\tilde\sigma_{i,n}\|_{H^{1}(\Omega)} \leq Ch^{\min\{1+\beta_i+\min\{\alpha, 1\}, 2\alpha\}}.
$$
Note that
\begin{align*}
|\langle \tilde \sigma_{i,n}-\tilde\sigma_{i},  \tilde{\xi}_{i} \rangle| = & |\langle  r^{\min\{\alpha-1, \beta_i\}}(\tilde\sigma_i-\tilde\sigma_{i,n}), r^{-\min\{\alpha-1, \beta_i\}} \tilde{\xi}_{i} \rangle| \\
\leq &  \|r^{\min\{\alpha-1, \beta_i\}}(\tilde\sigma_i-\tilde\sigma_{i,n})\|\|r^{-\min\{\alpha-1, \beta_i\}} \tilde{\xi}_{i}\| \leq Ch^{2\min\{1+\beta_i, \alpha\}}.
\end{align*}
The last two inequalities imply that
$$
\|T_3\| \leq Ch^{2\min\{1+\beta_i, \alpha\}}.
$$
Thus, we have
\be\label{Tforu}
|\tilde{c}_{i}-\tilde{c}_{i,n} |\leq \sum_{l=1}^3 \|T_l\| \leq Ch^{2\min\{1+\beta_i, \alpha\}}.
\ee
Plugging (\ref{uprojerr}), (\ref{eqn.11}), (\ref{eqn.112}) and (\ref{Tforu}) with $\alpha=\frac{1}{2}$ into (\ref{errbdd}), the conclusion holds.
\end{proof}

\begin{remark}\label{rk35}
For the following cases, we have $\min\{2(1+\beta_i),1\}=1$, (i) $1\leq i<N$; (ii) $i=N$ and $\beta_i\geq-\frac{1}{2}$. To better view $\|u-u_n\|_{H^1(\Omega)}$ in (\ref{uerr}), we explicitly show the value of $\min\{2(1+\beta_i),1\}$ and the value of $\gamma$ in Table \ref{uh1errtab} and \Cref{Rate}.
\begin{table}[!htbp]\tabcolsep0.03in
\caption{The value of $\min\{2(1+\beta_i),1\}$ and $\gamma$ in \Cref{thmuerr} for different $\omega$.}
\begin{tabular}[c]{|c|c|c|c|c|c|c|}
\hline
$\omega$ & $(0,\frac{\pi}{2}]$ & $(\frac{\pi}{2},\frac{2\pi}{3}]$ & $(\frac{2\pi}{3},\pi)$ & $(\pi,\frac{4\pi}{3}] $ & $(\frac{4\pi}{3},\frac{3\pi}{2}] $ & $(\frac{3\pi}{2}, 2\pi)$ \\
\hline
$\min\{2(1+\beta_1),1\}$ & $--$ & $2(1+\beta_1)$ & $1$ & $1$ & $1$ & $1$ \\
\hline
$\min\{2(1+\beta_2),1\}$ & $--$ & $--$ & $--$  & $2(1+\beta_2)$ & $1$ & $1$ \\
\hline
$\min\{2(1+\beta_3),1\}$ & $--$ & $--$ & $--$  & $--$ & $--$ & $2(1+\beta_3)$ \\
\hline
$\gamma$ & $1$ & $2(1+\beta_1)$ & $1$  & $2(1+\beta_2)$ & $1$ & $2(1+\beta_3)$ \\
\hline
\end{tabular}\label{uh1errtab}
\end{table}

\begin{figure}[h]
\centering
\subfigure{\includegraphics[width=0.49\textwidth]{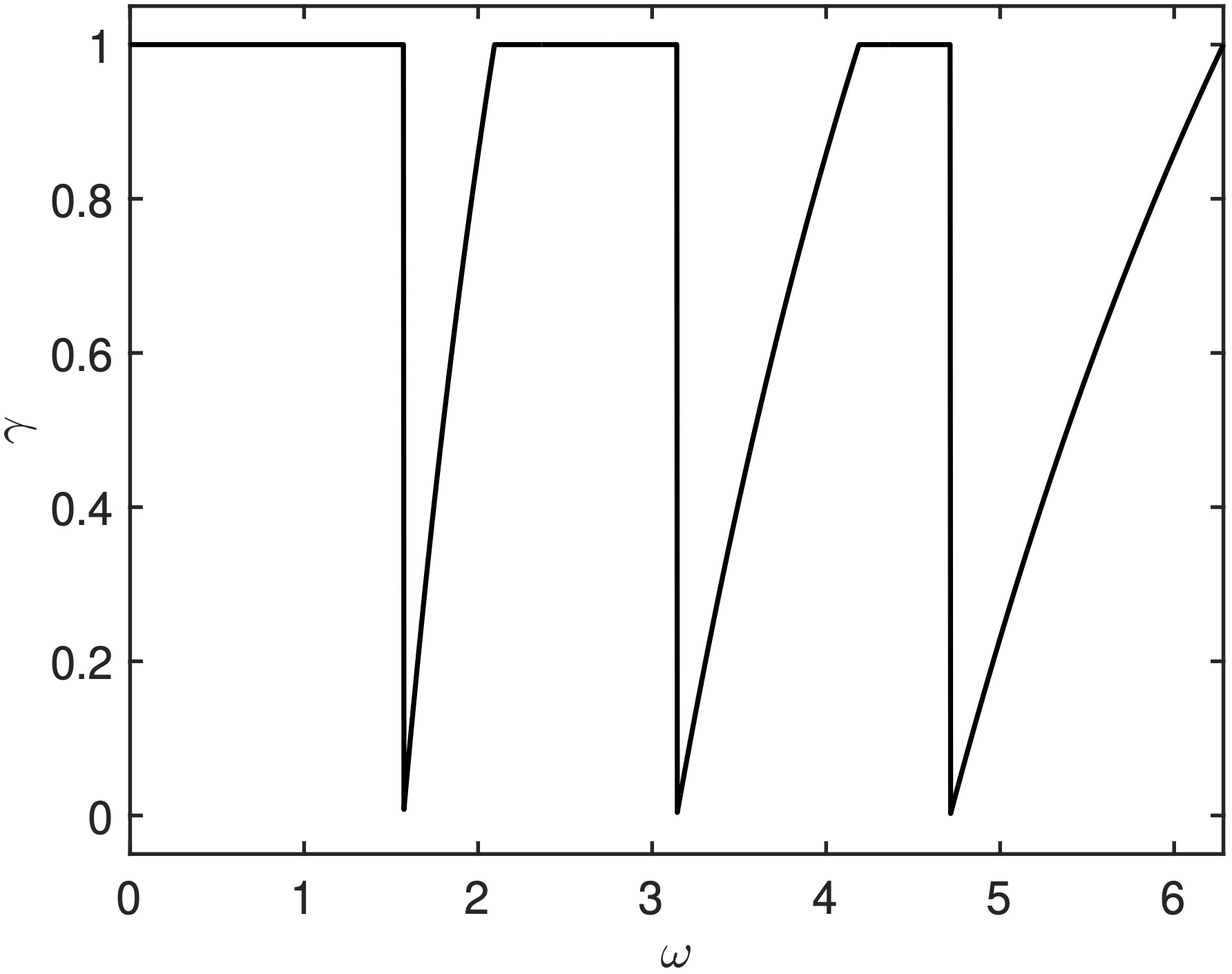}}
\caption{The $H^1$ convergence rate $\gamma$ in \Cref{thmuerr} for different $\omega$.}\label{Rate}
\end{figure}

\end{remark}

\section{Numerical illustrations}\label{sec-5}
In this section, we present numerical test results to validate our theoretical predictions for Algorithm \ref{femalg} solving the sixth order problem (\ref{eqnbi}). For comparison, we also implement the finite element method for the direct mixed formulation (\ref{eqn7}), referred to as the \textit{direct mixed finite element method}. We will utilize the following convergence rate as an indicator of the actual convergence rate
of the exact solutions $u, v, w$ in (\ref{eqnnew+}) are given, then calculate the convergence rate by
\begin{eqnarray}\label{rate0}
{\mathcal R}=\log_2\frac{|\phi-\phi_{j-1}|_{H^1(\Omega)}}{|\phi-\phi_j|_{H^1(\Omega)}},
\end{eqnarray}
otherwise,
\begin{eqnarray}\label{rate}
{\mathcal R}=\log_2\frac{|\phi_j-\phi_{j-1}|_{H^1(\Omega)}}{|\phi_{j+1}-\phi_j|_{H^1(\Omega)}}.
\end{eqnarray}
Here, $\phi_j$ represents the finite element solution on the mesh $\mathcal{T}_j$, obtained after $j$ refinements of the initial triangulation $\mathcal{T}_0$. It can be either $u_j$, $v_j$, or $w_j$, depending on the underlying Poisson problem.
In particular, suppose the actual convergence rate is $|\phi-\phi_j|_{H^1(\Omega)}=O(h^\beta)$ for $\beta>0$. Then, for the $P_1$ finite element method, the rate in (\ref{rate}) is also a good approximation of the exponent $\beta$ as the level of refinements $j$ increases \cite{LN18}.

We use the following cut-off function in Algorithm \ref{femalg}:
\begin{eqnarray*}
\eta(r; \tau, R)=
\left\{\begin{array}{ll}
0, & \text{if } r \geq R, \\
1, & \text{if } r \leq \tau R, \\
\frac{1}{2} - \frac{15}{16}\left( \frac{2r}{R(1-\tau)}-\frac{1+\tau}{1-\tau}\right) + \frac{5}{8}\left( \frac{2r}{R(1-\tau)}-\frac{1+\tau}{1-\tau}\right)^3 - \frac{3}{16}\left( \frac{2r}{R(1-\tau)}-\frac{1+\tau}{1-\tau}\right)^5 , & \text{otherwise.}
\end{array}\right.
\end{eqnarray*}
We set the default parameters $R=\frac{32}{5}, \tau = \frac{1}{8}$. If a different $R$ is used, it will be specified.

\begin{example}\label{P1S}
We solve the problem (\ref{eqnbi}) on different domains using both the direct mixed finite element method and Algorithm \ref{femalg} on quasi-uniform meshes obtained by midpoint refinements with the given initial mesh. We start with a ``wrong solution" $u\not \in H^3(\Omega)$,
\be\label{wrongsolu}
u(r,\theta)= \tilde{\eta}(r; \tau, R) r^{\frac{\pi}{\omega}}\sin\left(\frac{\pi}{\omega} \theta\right),
\ee
where $\tilde{\eta}(r; \tau, R)$ is also a cut-off function
\begin{eqnarray}\label{cutoffexact}
\tilde{\eta}(r; \tau, R)=
\left\{\begin{array}{ll}
0, & \text{if } r > R, \\
1, & \text{if } r < \tau R, \\
\frac{1}{2}+\sum_{i=0}^{6} C_i\left( \frac{2r}{R(1-\tau)} -\frac{1+\tau}{1-\tau} \right)^{2i+1}, & \text{otherwise,}
\end{array}\right.
\end{eqnarray}
with $R=\frac{32}{5}, \tau = \frac{1}{8}$, and the coefficients $C_i$ are determined by solving the linear system
$$
\tilde{\eta}^{(i)}(R; \tau, R) = 0, \quad i=0,\ldots,6.
$$
The source term $f$ is obtained by calculating
$$
f=-\Delta(\Delta(\Delta u)),
$$
and it can be verified that $f \in L^2(\Omega)$. Note that $u\not \in H^3(\Omega)$ and therefore $u$ is not the solution of the weak formulation (\ref{eqn.firstbi}) because the ``true solution" should be a function in $H^3(\Omega)$. The purpose of this example is to test the convergence of the finite element method for the direct mixed formulation and Algorithm \ref{femalg} to the ``spurious solution" $u$ in (\ref{wrongsolu}).

\noindent\textbf{Test case 1.} Take $\Omega$ as the triangle $\vartriangle QQ_1Q_2$ with $Q(0,0)$, $Q_1(16,0)$ and $Q_2(-8, 8\sqrt{3})$. The domain $\Omega$ with the initial mesh is shown in Figure \ref{Ex1Q2}(a), and the ``spurious solution" $u$ is shown in Figure \ref{Ex1Q2}(b). Here, $\omega=\angle Q_1QQ_2=\frac{2\pi}{3} \in (\frac{\pi}{2},\pi)$.

The direct mixed finite element solution $u_{10}^U$ and the difference $|u-u_{10}^U|$ are shown in Figure \ref{Ex1Q2}(c) and Figure \ref{Ex1Q2}(d), respectively. The error $\|u-u_j^U\|_{H^1(\Omega)}$ is shown in Table \ref{Ex1Q2Err}. These results indicate that the direct mixed finite element solution converges to the ``spurious solution" $u \not \in H^3(\Omega)$.
On the other hand, since $\omega=\frac{2\pi}{3} \in (\frac{\pi}{2},\pi)$, so it follows $N=1$ in Algorithm \ref{femalg} by checking Table \ref{TabomegaN}. The solution $u_{10}^A$ from Algorithm \ref{femalg} and the difference $|u-u_{10}^A|$ are shown in Figure \ref{Ex1Q2}(e) and Figure \ref{Ex1Q2}(f), respectively. The error $\|u-u_j^A\|_{H^1(\Omega)}$ is shown in Table \ref{Ex1Q2Err}. These results imply that the solution of Algorithm \ref{femalg} does not converge to the ``spurious solution", since the solution of Algorithm \ref{femalg} converges to the solution in $H^3(\Omega)$ as stated in Theorem \ref{solequthm}.

\begin{figure}[h]
\centering
\subfigure[]{\includegraphics[width=0.30\textwidth]{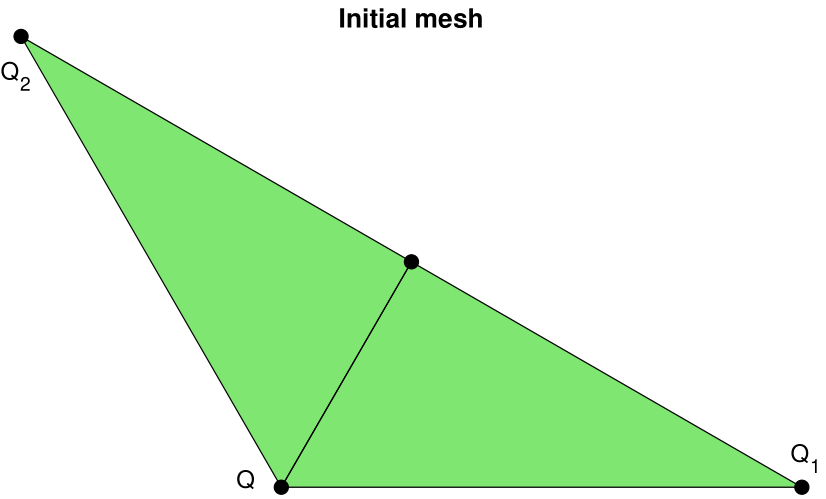}}
\subfigure[]{\includegraphics[width=0.30\textwidth]{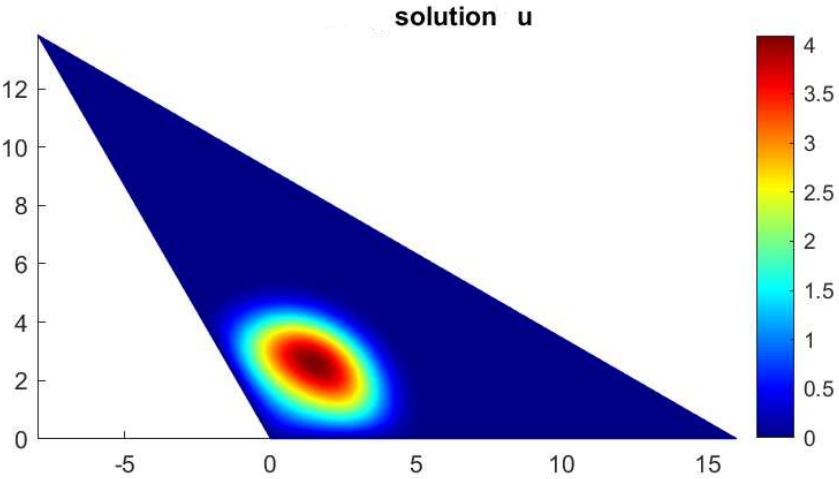}}
\subfigure[]{\includegraphics[width=0.30\textwidth]{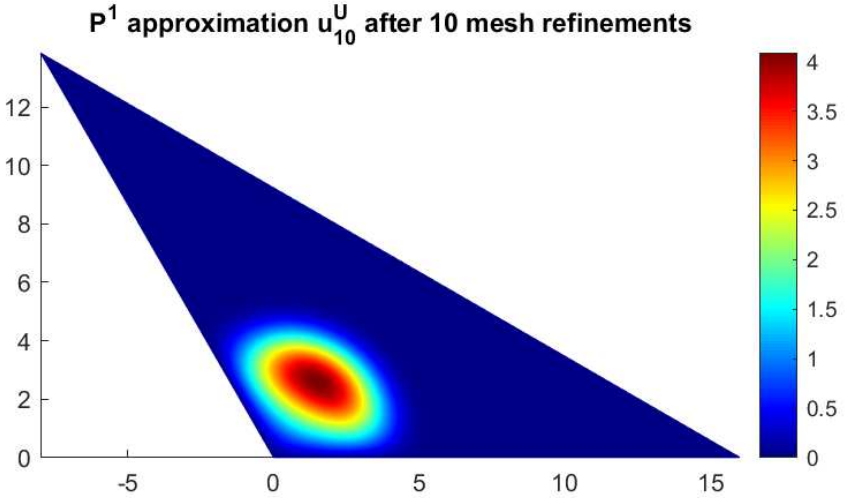}}
\subfigure[]{\includegraphics[width=0.30\textwidth]{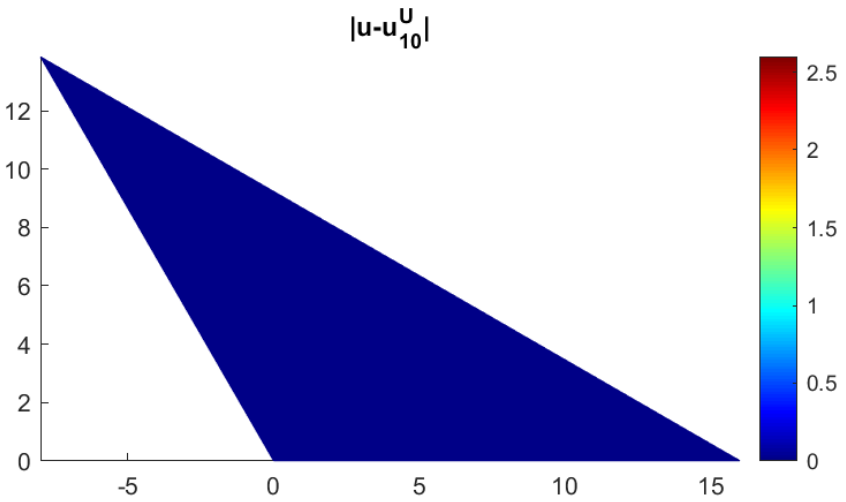}}
\subfigure[]{\includegraphics[width=0.30\textwidth]{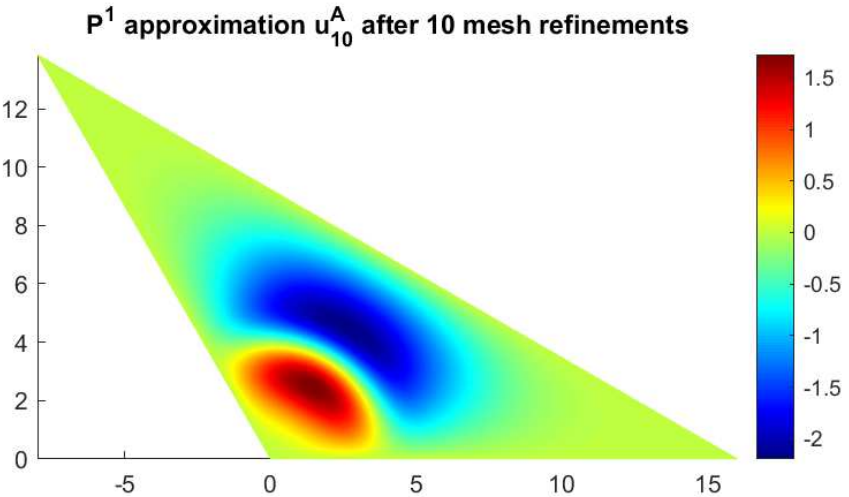}}
\subfigure[]{\includegraphics[width=0.30\textwidth]{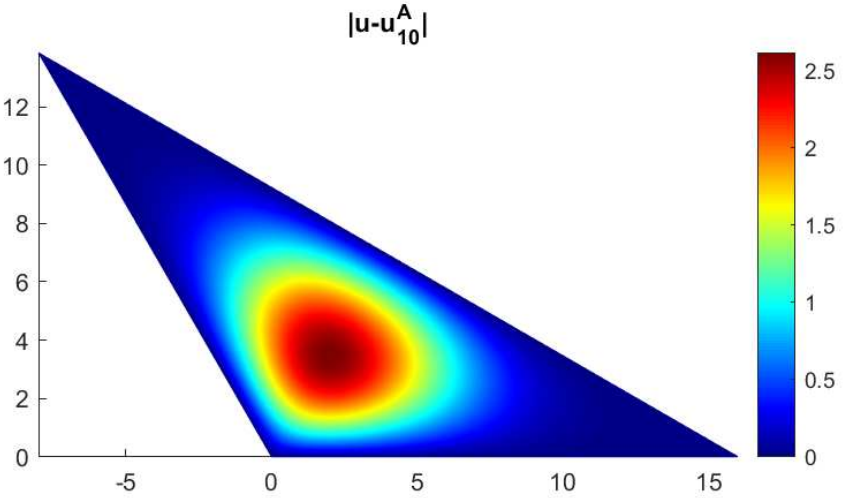}}
\caption{Example \ref{P1S} Test case 1:  (a) the domain and the initial mesh; (b) the ``spurious solution" $u$; (c) the direct mixed finite element solution $u^U_{10}$; (d) the difference $|u-u^U_{10}|$; (e) the solution $u^A_{10}$ from Algorithm \ref{femalg}; (f) the difference $|u-u^A_{10}|$.}\label{Ex1Q2}
\end{figure}

\begin{table}[!htbp]\tabcolsep0.03in
\caption{The $H^1$ error of the numerical solutions on quasi-uniform meshes.}
\begin{tabular}[c]{|c|c|c|c|c|}
\hline
\multirow{2}{*}{} & $j=7$ & $j=8$ & {$j=9$} & {$j=10$}  \\
\cline{2-5}
\hline
$\|u-u_j^U\|_{H^1(\Omega)}$  & 2.74964e-01 & 1.35594e-01 & 6.77391e-02 & 3.38605e-02 \\
\hline
$\|u-u_j^A\|_{H^1(\Omega)}$  &  6.07564 & 6.02331 & 6.00958 & 6.00306 \\
\cline{1-5}
\end{tabular}\label{Ex1Q2Err}
\end{table}

\noindent\textbf{Test case 2.} Here, we consider the domain $\Omega$ to be the polygon with vertices $Q(0,0)$, $Q_1(\frac{16\sqrt{3}}{3},0)$, $Q_2(\frac{16-8\sqrt{2}}{1+\sqrt{3}}, \frac{16-8\sqrt{2}}{1+\sqrt{3}}+8\sqrt{2})$, $Q_3(-8\frac{\sqrt{2}+\frac{2\sqrt{3}}{2}}{1+\frac{1}{\sqrt{3}}}, 8\sqrt{2}-8\frac{\sqrt{2}+\frac{2\sqrt{3}}{2}}{1+\frac{1}{\sqrt{3}}})$ and $Q_4(-\frac{8}{3}, -4\sqrt{3})$. Then we have $\omega=\angle Q_1QQ_4 \approx 1.383\pi \in (\pi, \frac{3\pi}{2})$. The domain $\Omega$ with the initial mesh is shown in Figure \ref{Ex1Q3}(a), and the ``spurious solution" $u$ is shown in Figure \ref{Ex1Q3}(b).

The direct mixed finite element solution $u_{10}^U$ and the difference $|u-u_{10}^U|$ are shown in Figure \ref{Ex1Q3}(c) and Figure \ref{Ex1Q3}(d), respectively. The error $\|u-u_j^U\|_{H^1(\Omega)}$ is shown in Table \ref{Ex1Q3Err}. These results imply that the direct mixed finite element solution converges to the ``spurious solution" $u \not \in H^3(\Omega)$.
On the other hand, since $\omega \in (\pi, \frac{3\pi}{2})$, we have $N=2$ in Algorithm \ref{femalg}. The solution $u_{10}^A$ of Algorithm \ref{femalg} and the difference $|u-u_{10}^A|$ are shown in Figure \ref{Ex1Q3}(e) and Figure \ref{Ex1Q3}(f), respectively. The error $\|u-u_j^A\|_{H^1(\Omega)}$ is shown in Table \ref{Ex1Q3Err}. These results imply that the solution of Algorithm \ref{femalg} does not converge to the ``spurious solution".

\begin{figure}[h]
\centering
\subfigure[]{\includegraphics[width=0.30\textwidth]{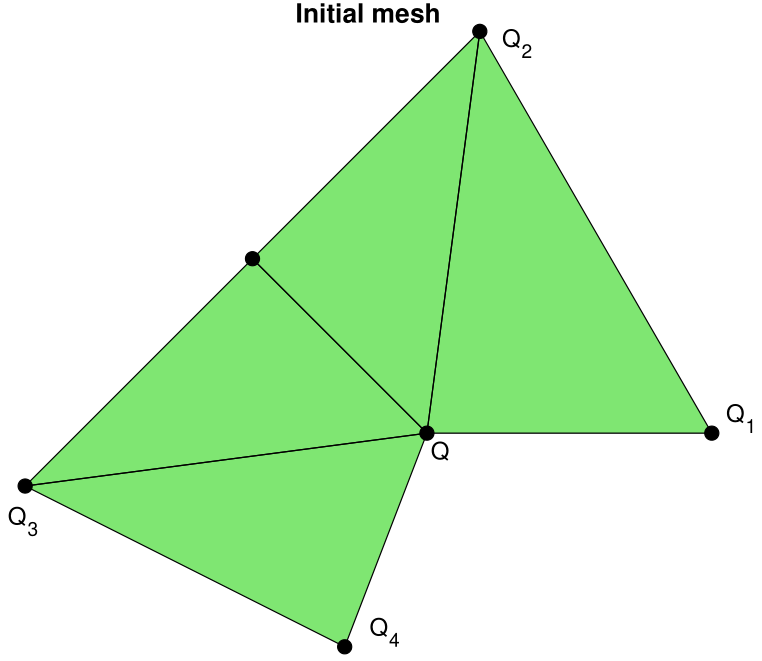}}
\subfigure[]{\includegraphics[width=0.30\textwidth]{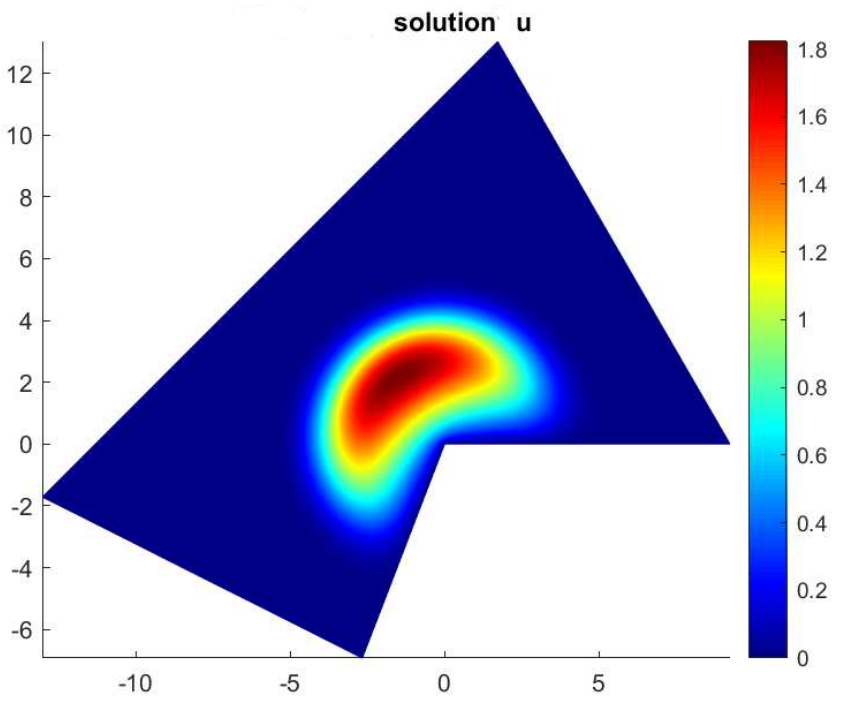}}
\subfigure[]{\includegraphics[width=0.30\textwidth]{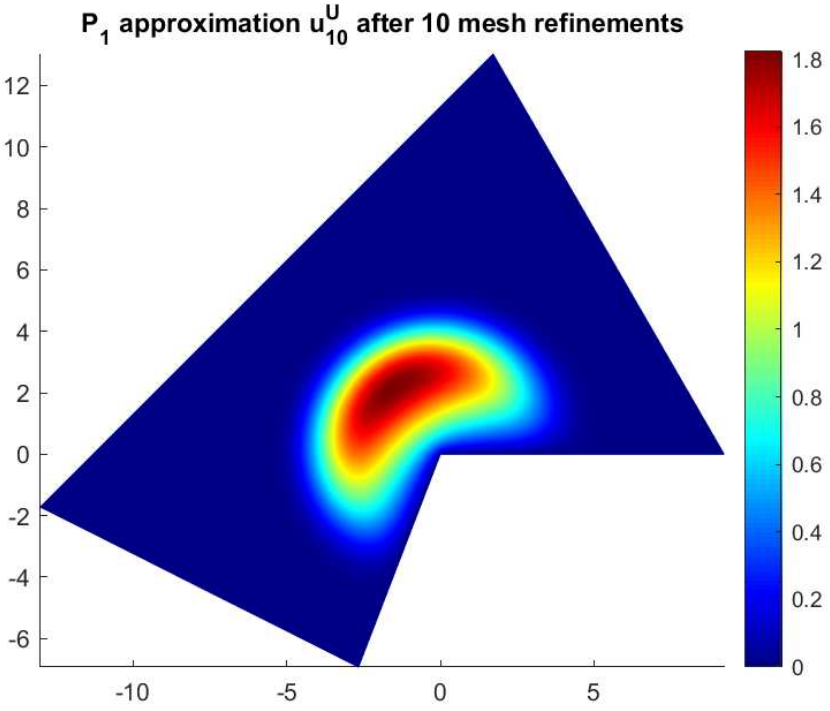}}
\subfigure[]{\includegraphics[width=0.30\textwidth]{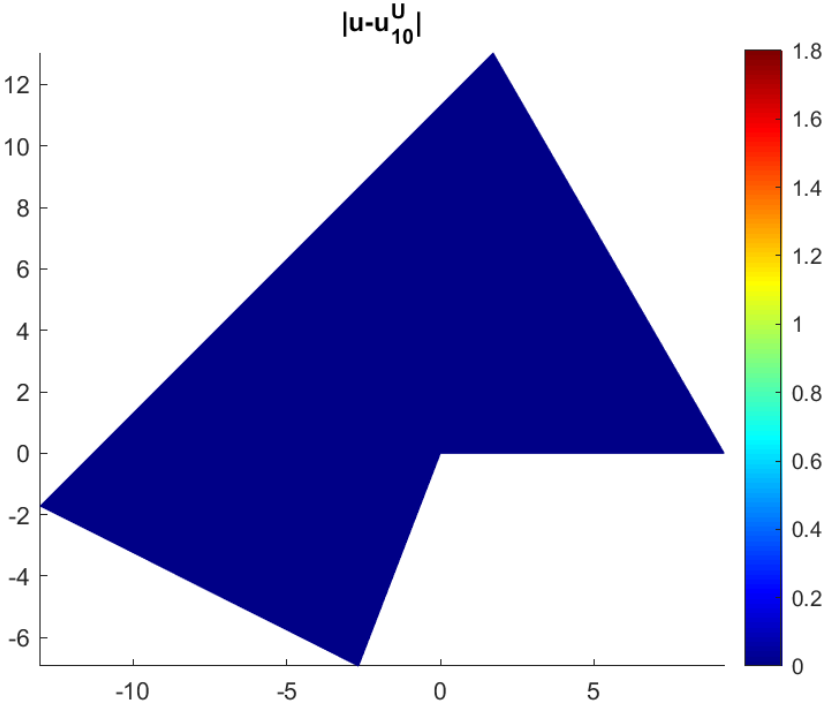}}
\subfigure[]{\includegraphics[width=0.30\textwidth]{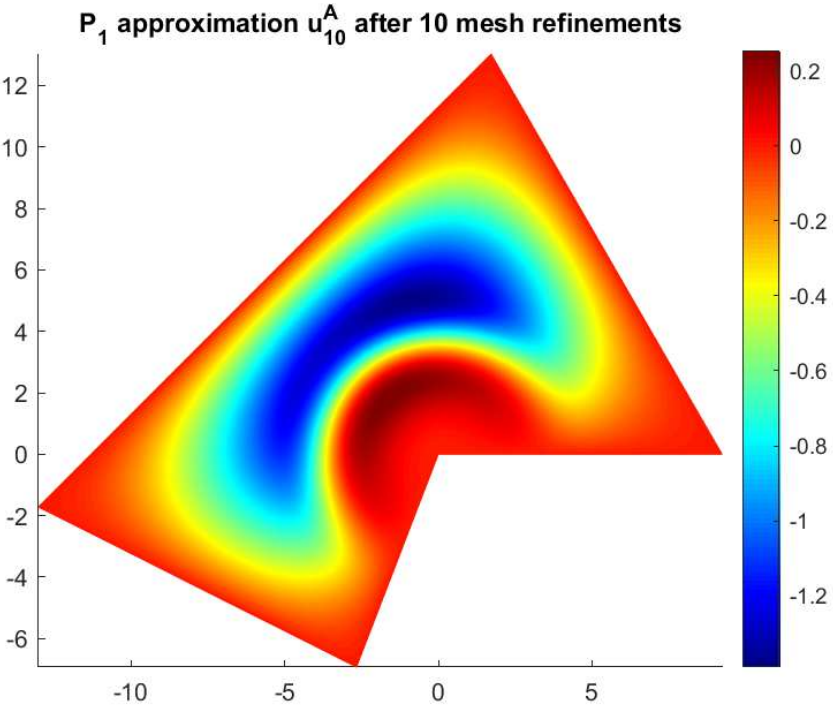}}
\subfigure[]{\includegraphics[width=0.30\textwidth]{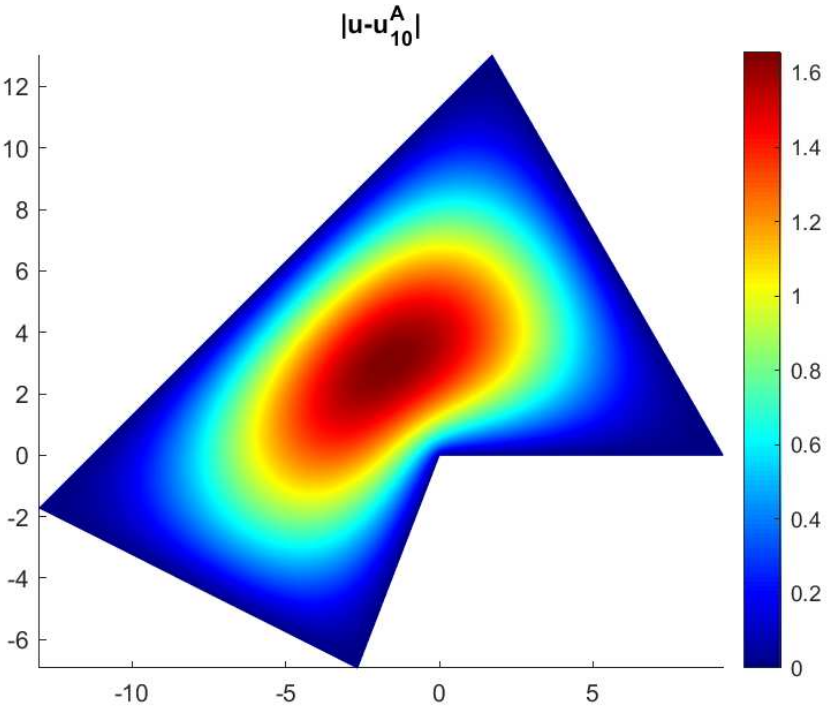}}
\caption{Example \ref{P1S} Test case 2:  (a) the domain and the initial mesh; (b) the ``spurious solution" $u$; (c) the direct mixed finite element solution $u^U_{10}$; (d) the difference $|u-u^U_{10}|$; (e) the solution $u^A_{10}$ from Algorithm \ref{femalg}; (f) the difference $|u-u^A_{10}|$.}\label{Ex1Q3}
\end{figure}

\begin{table}[!htbp]\tabcolsep0.03in
\caption{The $H^1$ error of the numerical solutions on quasi-uniform meshes.}
\begin{tabular}[c]{|c|c|c|c|c|}
\hline
\multirow{2}{*}{} & $j=7$ & $j=8$ & {$j=9$} & {$j=10$}  \\
\cline{2-5}
\hline
$\|u-u_j^U\|_{H^1(\Omega)}$  & 1.43517e-01 & 7.44186e-02 & 3.94988e-02 & 2.13310e-02 \\
\hline
$\|u-u_j^A\|_{H^1(\Omega)}$  &  4.08611 & 4.08457 & 4.08383 & 4.08329 \\
\cline{1-5}
\end{tabular}\label{Ex1Q3Err}
\end{table}

\noindent\textbf{Test case 3.} Consider the polygonal domain $\Omega$ with vertices $Q(0,0)$, $Q_1(\frac{16\sqrt{3}}{3},0)$, $Q_2(\frac{16-8\sqrt{2}}{1+\sqrt{3}}, \frac{16-8\sqrt{2}}{1+\sqrt{3}}+8\sqrt{2})$, $Q_3(-8\frac{\sqrt{2}+\frac{2\sqrt{3}}{2}}{1+\frac{1}{\sqrt{3}}}, 8\sqrt{2}-8\frac{\sqrt{2}+\frac{2\sqrt{3}}{2}}{1+\frac{1}{\sqrt{3}}})$ and $Q_4(4, -8\sqrt{3})$. Then we have $\omega=\angle Q_1QQ_4 \approx 1.589 \pi\in (\frac{3\pi}{2}, 2\pi)$. The domain $\Omega$ with the initial mesh is shown in Figure \ref{Ex1Q4}(a), and the ``spurious solution" $u$ is shown in Figure \ref{Ex1Q4}(b).

The direct mixed finite element solution $u_{10}^U$ and the difference $|u-u_{10}^U|$ are shown in Figure \ref{Ex1Q4}(c) and Figure \ref{Ex1Q4}(d), respectively. The error $\|u-u_j^U\|_{H^1(\Omega)}$ is shown in Table \ref{Ex1Q4Err}. These results continue to indicate that the direct mixed finite element solution converges to the ``spurious solution" $u \not \in H^3(\Omega)$.
On the other hand, since $\omega \in (\frac{3\pi}{2},2\pi)$, it follows $N=3$ in Algorithm \ref{femalg}. The solution $u_{10}^A$ of Algorithm \ref{femalg} and the difference $|u-u_{10}^A|$ are shown in Figure \ref{Ex1Q4}(e) and Figure \ref{Ex1Q4}(f), respectively. The error $\|u-u_j^A\|_{H^1(\Omega)}$ is shown in Table \ref{Ex1Q4Err}. These results confirm that the solution of Algorithm \ref{femalg} does not converge to the ``spurious solution".

\begin{figure}[h]
\centering
\subfigure[]{\includegraphics[width=0.29\textwidth]{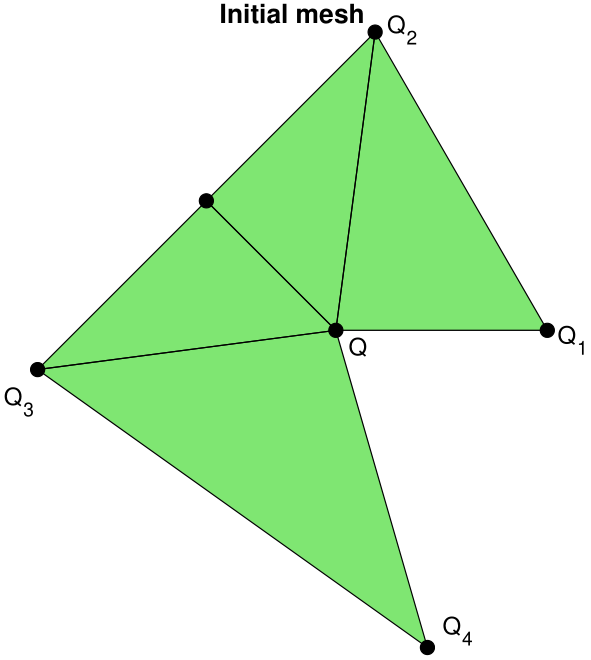}}
\subfigure[]{\includegraphics[width=0.30\textwidth]{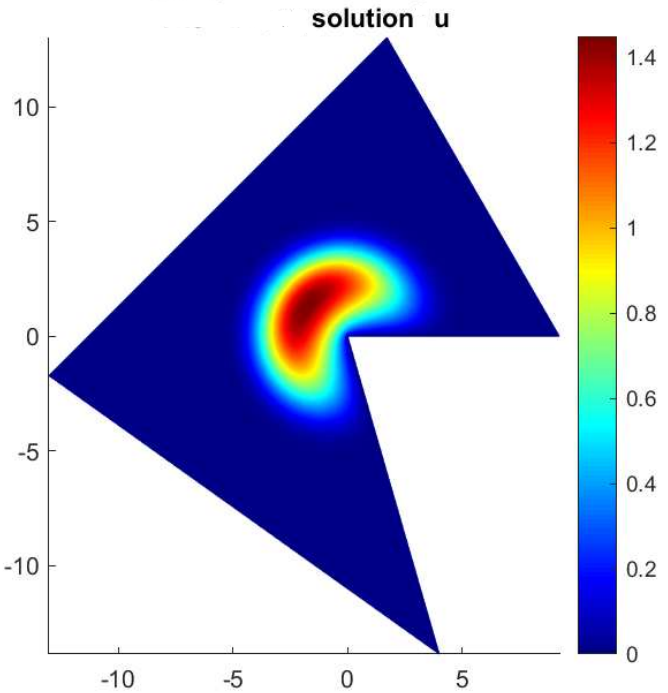}}
\subfigure[]{\includegraphics[width=0.30\textwidth]{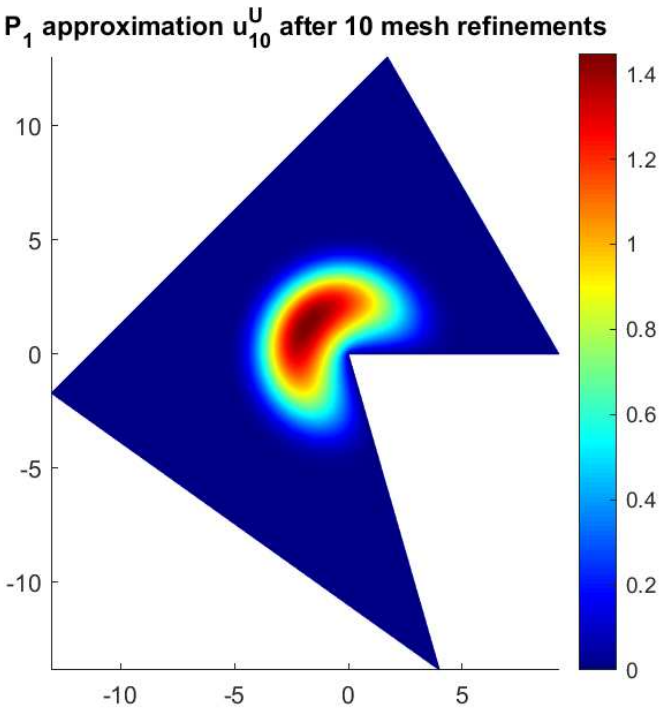}}
\subfigure[]{\includegraphics[width=0.30\textwidth]{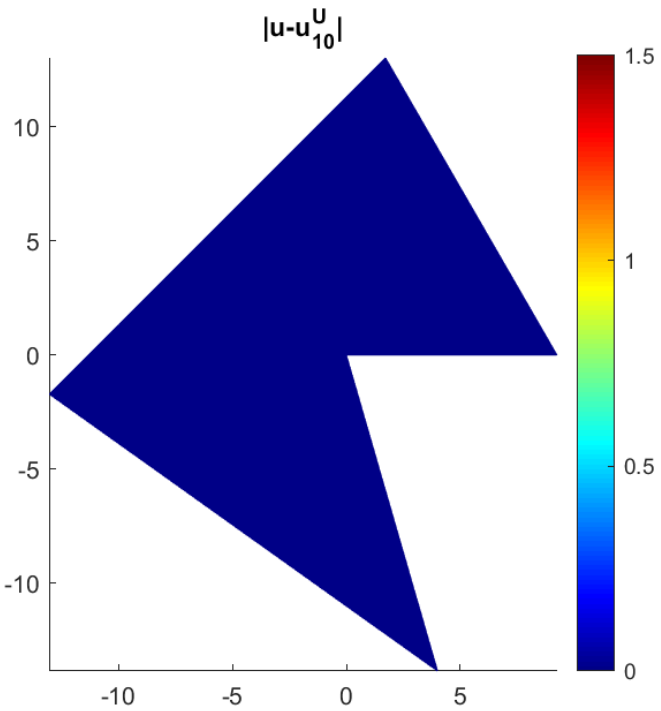}}
\subfigure[]{\includegraphics[width=0.30\textwidth]{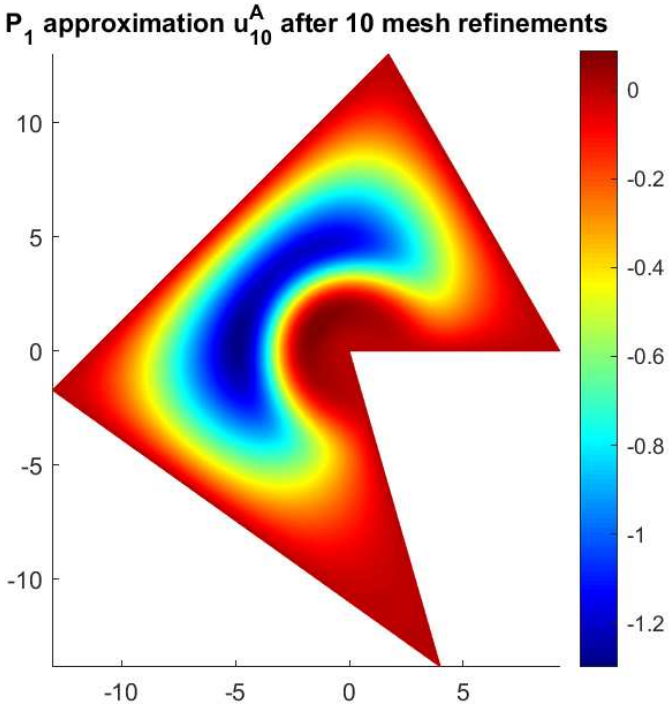}}
\subfigure[]{\includegraphics[width=0.30\textwidth]{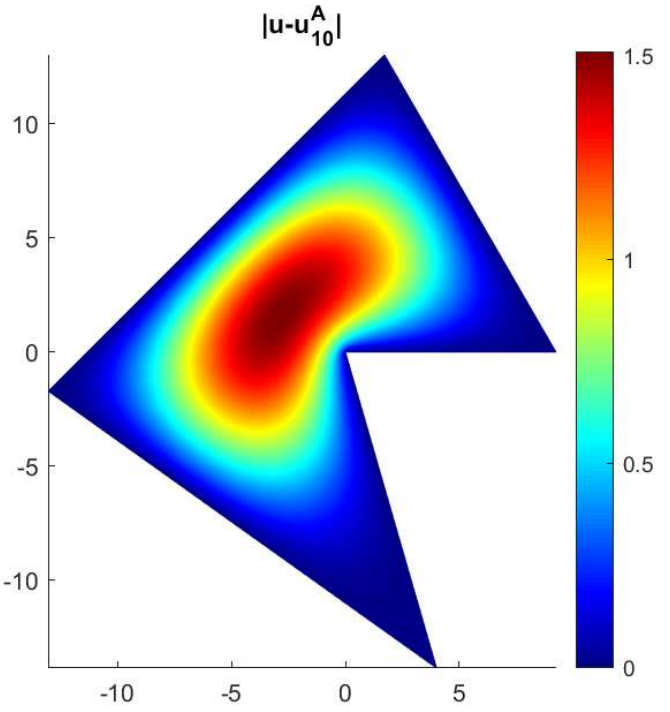}}
\caption{Example \ref{P1S} Test Case 3:  (a) the domain and the initial mesh; (b) the ``spurious solution" $u$; (c) the direct mixed finite element solution $u^U_{10}$; (d) the difference $|u-u^U_{10}|$; (e) the solution $u^A_{10}$ from Algorithm \ref{femalg}; (f) the difference $|u-u^A_{10}|$.}\label{Ex1Q4}
\end{figure}

\begin{table}[!htbp]\tabcolsep0.03in
\caption{The $H^1$ error of the numerical solutions on quasi-uniform meshes.}
\begin{tabular}[c]{|c|c|c|c|c|}
\hline
\multirow{2}{*}{} & $j=7$ & $j=8$ & {$j=9$} & {$j=10$}  \\
\cline{2-5}
\hline
$\|u-u_j^U\|_{H^1(\Omega)}$  & 1.474223e-01 & 8.67096e-02 & 5.25520e-02 & 3.25455e-02 \\
\hline
$\|u-u_j^A\|_{H^1(\Omega)}$  &  3.863711 & 3.85981 & 3.85832 & 3.85767 \\
\cline{1-5}
\end{tabular}\label{Ex1Q4Err}
\end{table}
\end{example}

\begin{example}\label{P1S2}
We solve the triharmonic problem in Example \ref{P1S} again using the direct mixed finite element method and Algorithm \ref{femalg} on quasi-uniform meshes.
Here, we take the solution $u_{ex}$ of the following Poisson problem as the exact solution,
\be\label{exactsolu+}
-\Delta u_{ex}= f_0 - \sum_{i=1}^N c_i \sigma_i \text{ in } \Omega, \qquad u_{ex} =0 \text{ on } \partial \Omega,
\ee
where
$$
f_0= -\Delta \left( \tilde{\eta}(r; \tau, R) r^{\frac{N\pi}{\omega}}\sin\left(\frac{N\pi}{\omega} \theta\right) \right) \in H_0^1(\Omega),
$$
with $ \tilde{\eta}(r; \tau, R)$ given in (\ref{cutoffexact}), $\sigma_i$ given in (\ref{H1part}),
and $c_i$ is the solution of the linear system (\ref{scoef0}). Note that the function $f_0=-\Delta u$ for $u$ in (\ref{wrongsolu}).
By Lemma \ref{Poissonreg}, we have $u_{ex}\in H^3(\Omega)$ and it satisfies
$$
-\Delta^3 u_{ex} = -\Delta^2 ( \Delta u_{ex}) =  \Delta^2 f_0 - \sum_{i=1}^N c_i \Delta^2 \sigma_i = \Delta^2 f_0 = -\Delta^2 (\Delta u)= -\Delta(\Delta (\Delta u)) = f,
$$
where we have used the result in Lemma \ref{slemH1}. Here, the source term $f$ is the same as that in Example \ref{P1S}.
The purpose of this example is to test the convergence of the direct mixed finite element method and Algorithm \ref{femalg} to the exact solution $u_{ex}$ in (\ref{exactsolu+}). From Test case 2 to Test case 4, we will use the finite element method solution $u_{exn}$ 
(instead of using the complicated notation $u_{ex,n+1}$) of (\ref{exactsolu+}) 
on mesh $\mathcal{T}_{n+1}$ as an approximation of $u_{ex}$.

\noindent\textbf{Test case 1.} Take $\Omega$ as the triangle $\vartriangle QQ_1Q_2$ with $Q(0,0)$, $Q_1(8,0)$ and $Q_2(4, 4\sqrt{3})$. In this case, the exact solution
$u_{ex} = u$ for a given $u$ in (\ref{wrongsolu}), and its contour is given in Figure \ref{Ex2Q1}(a). Here, $\omega=\angle Q_1QQ_2=\frac{\pi}{3} \in (0,\frac{\pi}{2})$. Thus, Algorithm \ref{femalg} coincides with the direct mixed finite element method.
The solution $u_{10}^A(=u_{10}^U)$ from Algorithm \ref{femalg} and the difference $|u-u_{10}^A|$ are shown in Figure \ref{Ex2Q1}(b) and Figure \ref{Ex2Q1}(c), respectively. The error $\|u-u_j^A\|_{H^1(\Omega)}$ and convergence rate $\mathcal{R}$ are shown in Table \ref{Ex2Q1Err}. These results show that the solution of Algorithm \ref{femalg} converges to the exact solution in the optimal convergence rate $\mathcal{R}=1$, which coincides with the result in Theorem \ref{thmuerr} or Table \ref{uh1errtab}.

\begin{figure}[h]
\centering
\subfigure[]{\includegraphics[width=0.30\textwidth]{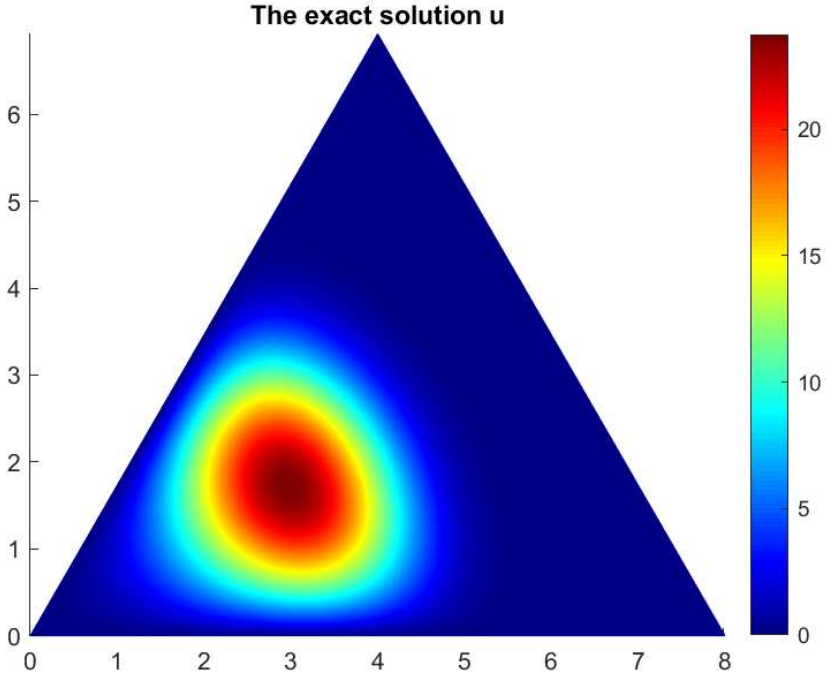}}
\subfigure[]{\includegraphics[width=0.30\textwidth]{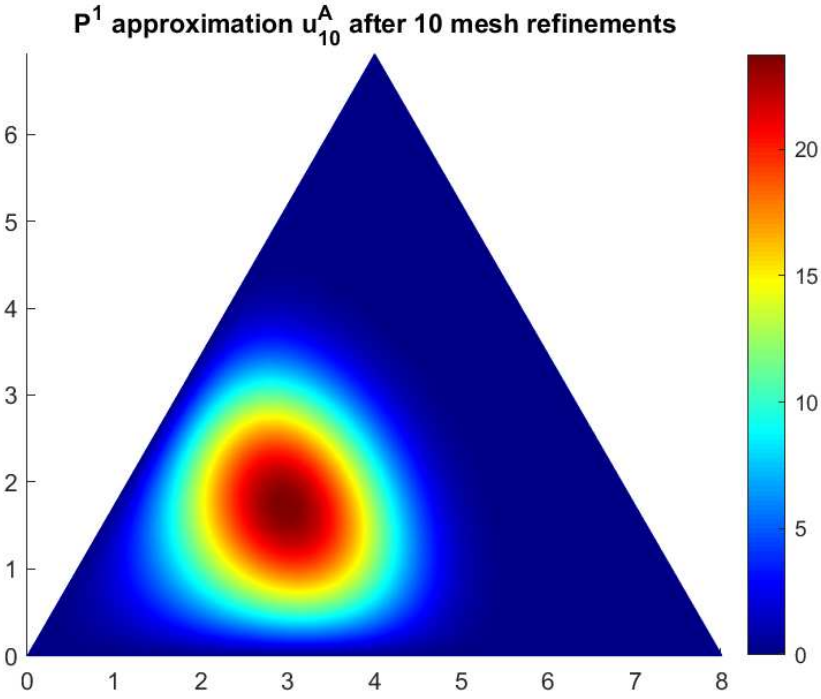}}
\subfigure[]{\includegraphics[width=0.30\textwidth]{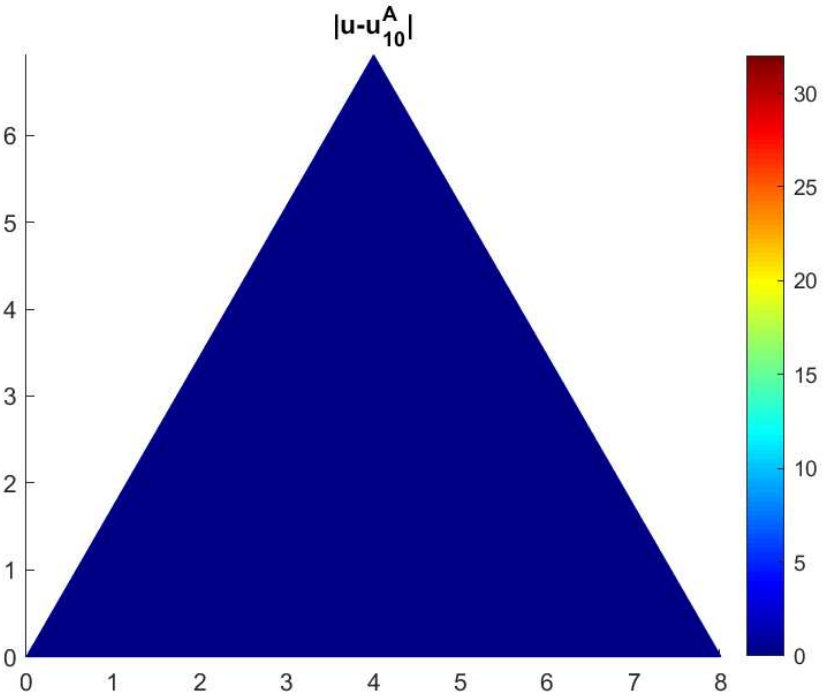}}
\caption{Example \ref{P1S2} Test case 1: (a) the exact solution $u$; (b) the solution $u^A_{10}$ from Algorithm \ref{femalg}; (c) the difference $|u-u^A_{10}|$.}\label{Ex2Q1}
\end{figure}

\begin{table}[!htbp]\tabcolsep0.03in
\caption{The $H^1$ error and convergence rate $\mathcal R$ for Example \ref{P1S2} Test case 1.}
\begin{tabular}[c]{|c|c|c|c|c|}
\hline
\multirow{2}{*}{} & {$j=5$} & $j=6$ & $j=7$ & {$j=8$} \\
\cline{2-5}
\hline
$\|u-u_j^A\|_{H^1(\Omega)}$  &  1.09202 & 5.45465e-01 & 2.72663e-01 & 1.36323e-01 \\
\hline
$\mathcal R$  &  $--$ & 1.00 & 1.00 & 1.00 \\
\cline{1-5}
\end{tabular}\label{Ex2Q1Err}
\end{table}

\noindent\textbf{Test case 2.} We consider the same domain and initial mesh (see Figure \ref{Ex1Q2}(a)) as Test case 1 in Example \ref{P1S}. Note that $\omega=\angle Q_1QQ_2=\frac{2\pi}{3} \in (\frac{\pi}{2},\pi)$.
The finite element solution $u_{exn}$ of the exact solution $u_{ex}$ is shown in Figure \ref{Ex2Q2}(a).
The direct mixed finite element solution $u_{10}^U$ and the difference $|u_{exn}-u_{10}^U|$ are shown in Figure \ref{Ex1Q2}(c) and Figure \ref{Ex2Q2}(b), respectively. The error $\|u_{exn}-u_j^U\|_{H^1(\Omega)}$ is shown in Table \ref{Ex2Q2Err}. These results indicate that the direct mixed finite element solution does not converge to the exact solution.
Note that $N=1$ in Algorithm \ref{femalg}, the solution $u_{10}^A$ from Algorithm \ref{femalg} and the difference $|u_{exn}-u_{10}^A|$ are shown in Figure \ref{Ex1Q2}(e) and Figure \ref{Ex2Q2}(c), respectively. The error $\|u_{exn}-u_j^A\|_{H^1(\Omega)}$ is shown in Table \ref{Ex2Q2Err}. These results imply that the solution of Algorithm \ref{femalg} converges to the exact solution.

\begin{figure}[h]
\centering
\subfigure[]{\includegraphics[width=0.30\textwidth]{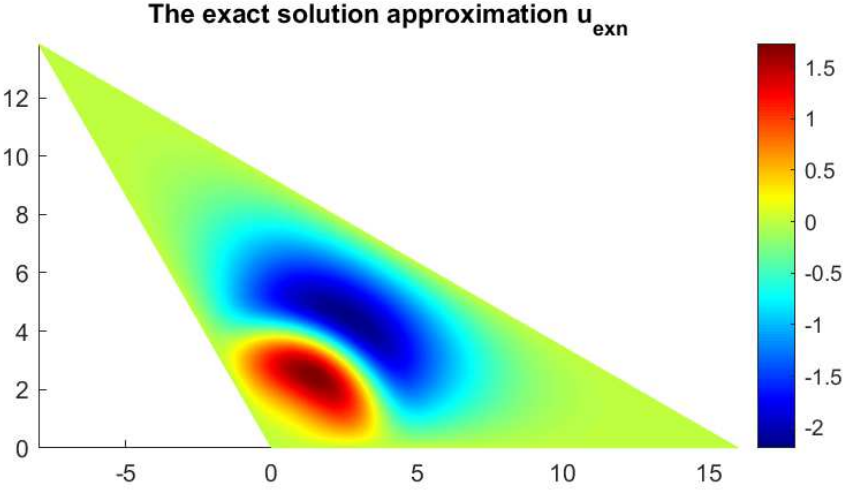}}
\subfigure[]{\includegraphics[width=0.30\textwidth]{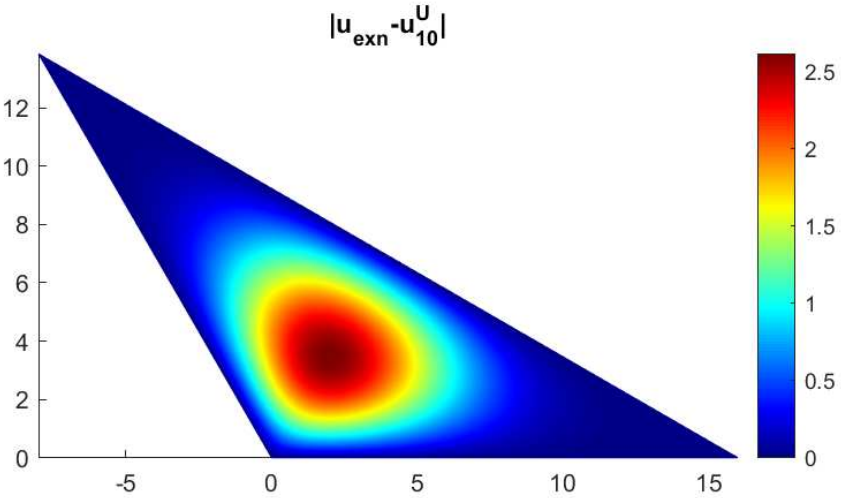}}
\subfigure[]{\includegraphics[width=0.30\textwidth]{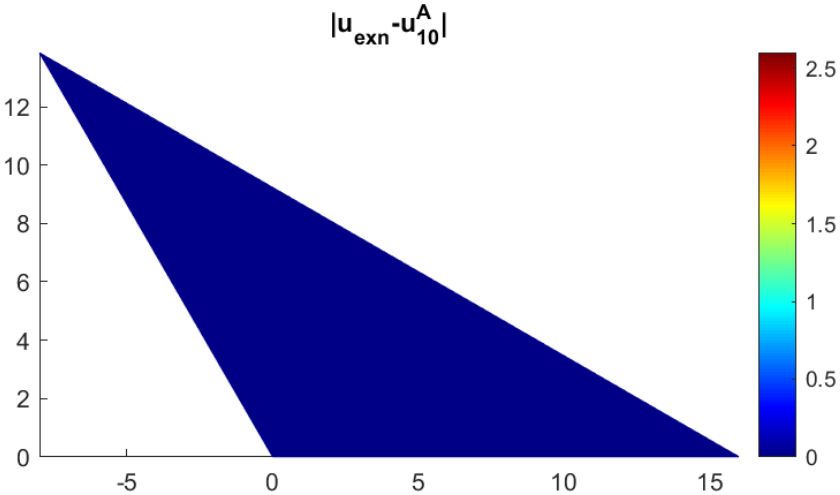}}
\caption{Example \ref{P1S2} Test case 2: (a) the exact solution approximation $u$; (b) the difference $|u_{exn}-u^U_{10}|$; (c) the difference $|u_{exn}-u^A_{10}|$.}\label{Ex2Q2}
\end{figure}

\begin{table}[!htbp]\tabcolsep0.03in
\caption{The $H^1$ error of the numerical solutions on quasi-uniform meshes.}
\begin{tabular}[c]{|c|c|c|c|c|}
\hline
\multirow{2}{*}{} & {$j=6$} & $j=7$ & $j=8$ & {$j=9$} \\
\cline{2-5}
\hline
$\|u_{exn}-u_j^U\|_{H^1(\Omega)}$  & 5.98206 & 6.01120 & 6.00363 & 5.99948 \\
\hline
$\|u_{exn}-u_j^A\|_{H^1(\Omega)}$  &  5.67208e-02 & 1.47272e-02 & 6.62074e-03 & 3.43917e-03 \\
\cline{1-5}
\end{tabular}\label{Ex2Q2Err}
\end{table}

\noindent\textbf{Test case 3.} We consider the same domain and initial mesh (see Figure \ref{Ex1Q3}(a)) as Test case 2 in Example \ref{P1S}. Recall that $\omega=\angle Q_1QQ_4 \approx 1.383\pi \in (\pi, \frac{3\pi}{2})$. The exact solution $u_{exn}$ is shown in Figure \ref{Ex2Q3}(a).
The direct mixed finite element solution $u_{10}^U$ and the difference $|u_{exn}-u_{10}^U|$ are shown in Figure \ref{Ex1Q3}(c) and Figure \ref{Ex2Q3}(b), respectively. The error $\|u_{exn}-u_j^U\|_{H^1(\Omega)}$ is shown in Table \ref{Ex2Q3Err}. These results indicate that the direct mixed finite element solution does not converge to the exact solution.
Note that $N=2$ in Algorithm \ref{femalg} in this case. The solution $u_{10}^A$ of Algorithm \ref{femalg} and the difference $|u_{exn}-u_{10}^A|$ are shown in Figure \ref{Ex1Q3}(e) and Figure \ref{Ex2Q3}(c), respectively. The error $\|u_{exn}-u_j^A\|_{H^1(\Omega)}$ is shown in Table \ref{Ex2Q3Err}. These results also imply that the solution of Algorithm \ref{femalg} converges to the exact solution.

\begin{figure}[h]
\centering
\subfigure[]{\includegraphics[width=0.30\textwidth]{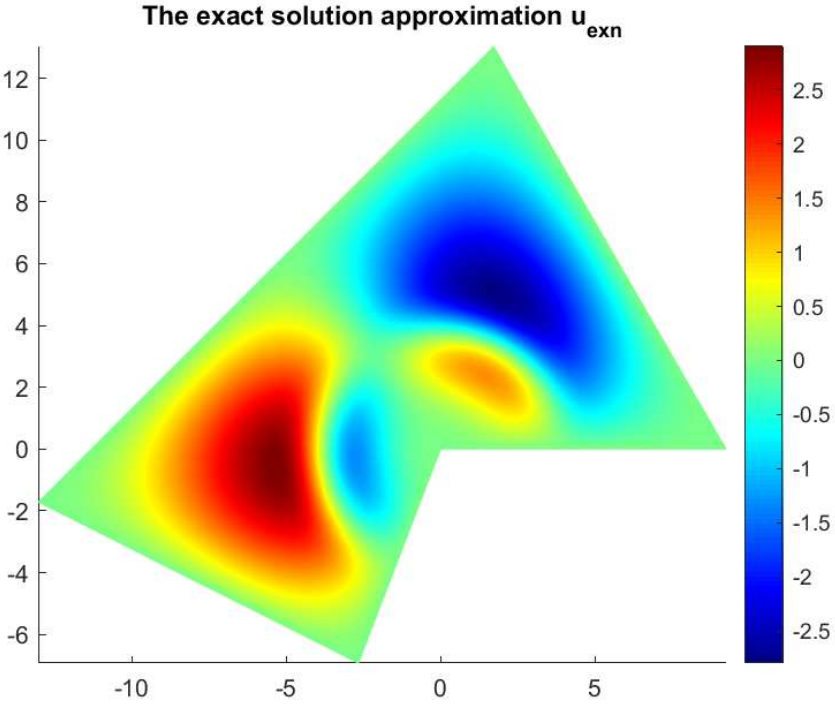}}
\subfigure[]{\includegraphics[width=0.30\textwidth]{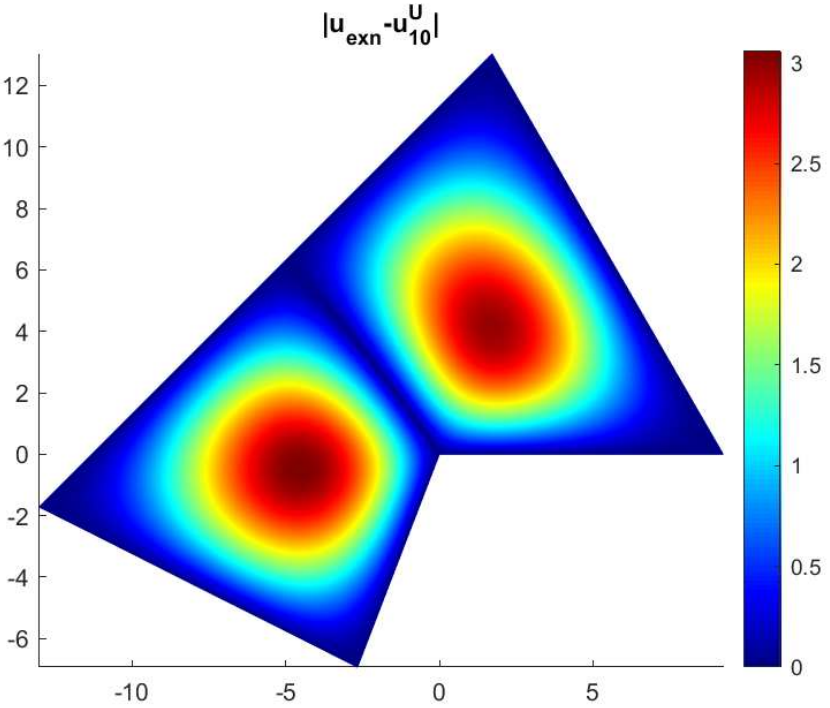}}
\subfigure[]{\includegraphics[width=0.30\textwidth]{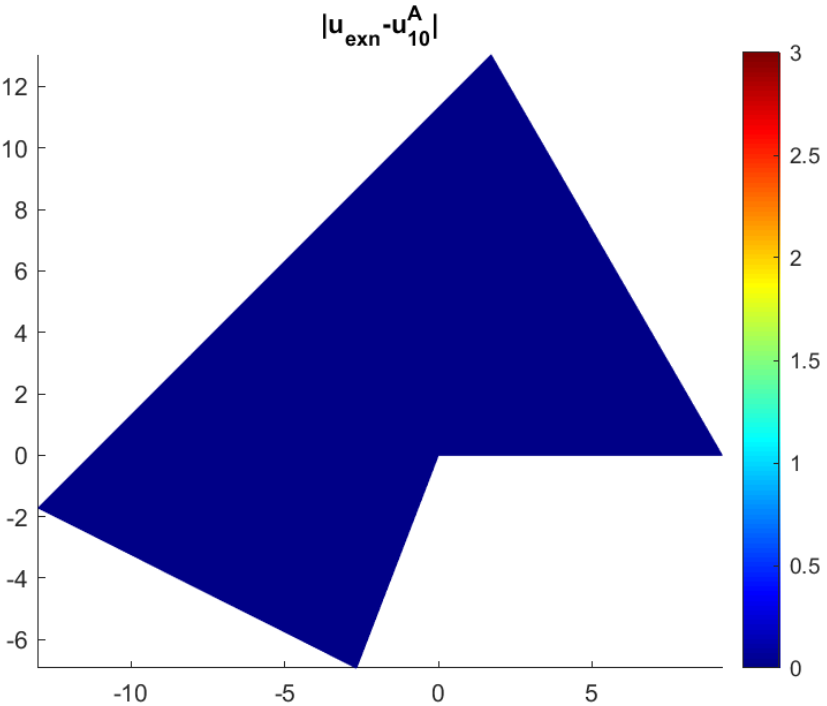}}
\caption{Example \ref{P1S2} Test case 3: (a) the exact solution u; (b) the difference $|u_{exn}-u^U_{10}|$; (c) the difference $|u_{exn}-u^A_{10}|$.}\label{Ex2Q3}
\end{figure}

\begin{table}[!htbp]\tabcolsep0.03in
\caption{The $H^1$ error of the numerical solutions on quasi-uniform meshes.}
\begin{tabular}[c]{|c|c|c|c|c|}
\hline
\multirow{2}{*}{} & {$j=6$} & $j=7$ & $j=8$ & {$j=9$} \\
\cline{2-5}
\hline
$\|u_{exn}-u_j^U\|_{H^1(\Omega)}$  & 9.67666 & 9.64665 & 9.63404 & 9.63164 \\
\hline
$\|u_{exn}-u_j^A\|_{H^1(\Omega)}$  & 5.27303e-02 & 2.09405e-02 & 1.01081e-02 & 4.20655e-03 \\
\cline{1-5}
\end{tabular}\label{Ex2Q3Err}
\end{table}

\noindent\textbf{Test case 4.} We consider the same domain and initial mesh (see Figure \ref{Ex1Q4}(a)) as Test case 3 in Example \ref{P1S}. Recall that $\omega=\angle Q_1QQ_4 \approx 1.589 \pi\in (\frac{3\pi}{2}, 2\pi)$. The approximation $u_{exn}$ of the exact solution is shown in Figure \ref{Ex2Q4}(a).
The direct mixed finite element solution $u_{10}^U$ and the difference $|u_{exn}-u_{10}^U|$ are shown in Figure \ref{Ex1Q4}(c) and Figure \ref{Ex2Q4}(b), respectively. The error $\|u_{exn}-u_j^U\|_{H^1(\Omega)}$ is shown in Table \ref{Ex1Q4Err}. These results continue to indicate that the direct mixed finite element solution does not converge to the exact solution.
Note that $N=3$ in Algorithm \ref{femalg}, the solution $u_{10}^A$ of Algorithm \ref{femalg} and the difference $|u_{exn}-u_{10}^A|$ are shown in Figure \ref{Ex1Q4}(e) and Figure \ref{Ex2Q4}(c), respectively. The error $\|u_{exn}-u_j^A\|_{H^1(\Omega)}$ is shown in Table \ref{Ex2Q4Err}. These results confirm that the solution of Algorithm \ref{femalg} converges to the exact solution.

\begin{figure}[h]
\centering
\subfigure[]{\includegraphics[width=0.30\textwidth]{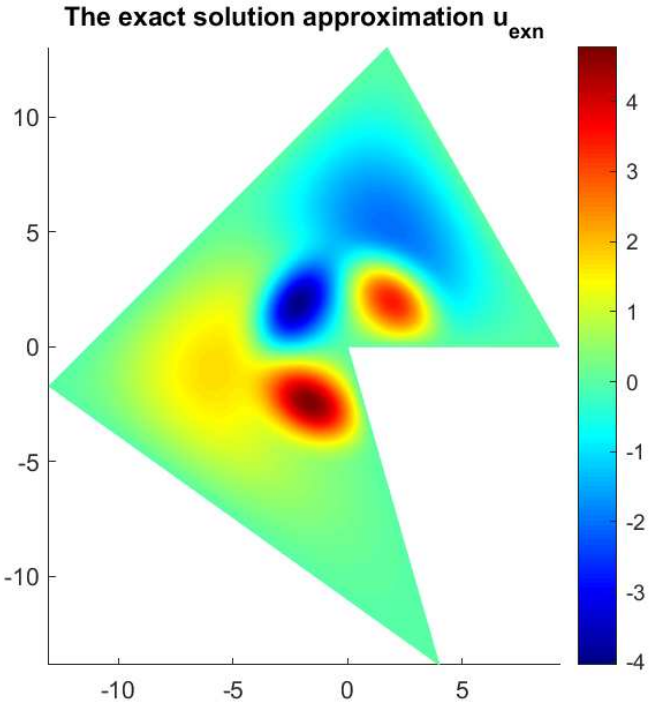}}
\subfigure[]{\includegraphics[width=0.30\textwidth]{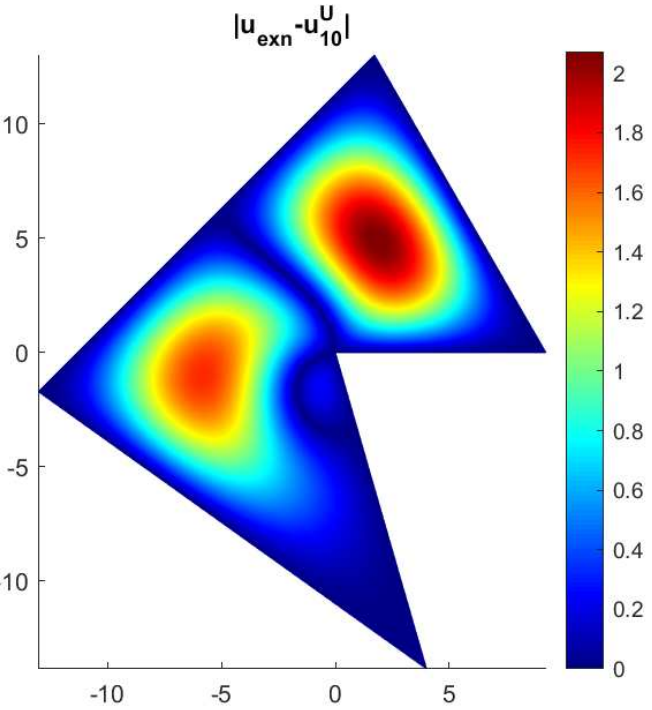}}
\subfigure[]{\includegraphics[width=0.30\textwidth]{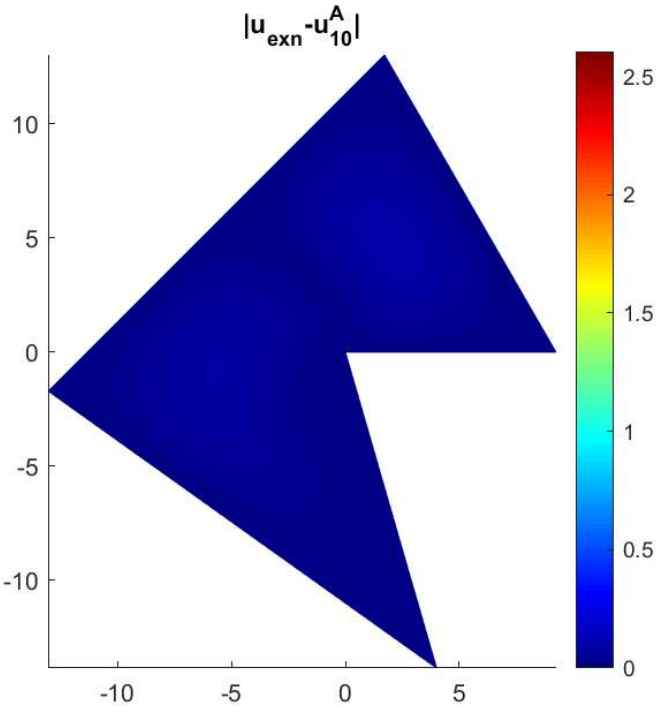}}
\caption{Example \ref{P1S2} Test case 4: (a) the exact solution $u$; (b) the difference $|u_{exn}-u^U_{10}|$; (c) the difference $|u_{exn}-u^A_{10}|$.}\label{Ex2Q4}
\end{figure}

\begin{table}[!htbp]\tabcolsep0.03in
\caption{The $H^1$ error of the numerical solutions on quasi-uniform meshes.}
\begin{tabular}[c]{|c|c|c|c|c|}
\hline
\multirow{2}{*}{} & {$j=6$} & $j=7$ & $j=8$ & {$j=9$} \\
\cline{2-5}
\hline
$\|u_{exn}-u_j^U\|_{H^1(\Omega)}$  & 7.47470 & 6.98223 & 6.60342 & 6.31031 \\
\hline
$\|u_{exn}-u_j^A\|_{H^1(\Omega)}$  & 6.79611e-01 & 4.98616e-01 & 3.78626e-01 & 2.93364e-01 \\
\cline{1-5}
\end{tabular}\label{Ex2Q4Err}
\end{table}

\end{example}

\begin{example}\label{P1h}  In this example, we investigate the convergence of Algorithm \ref{femalg} by considering equation  (\ref{eqnbi}) with $f=\sin\left(\frac{N\pi}{\omega}\theta\right)$ on different domains with angle $\omega$ categorized in Theorem \ref{thmuerr} or Table \ref{uh1errtab}, where $N$ is shown in Table \ref{TabomegaN}. For $\omega < \frac{\pi}{2}$, the numerical test on convergence rate can be found in Example \ref{P1S2} Test case 1. In the rest of this example, we focus on $\omega>\frac{\pi}{2}$.

\noindent\textbf{Test case 1.}  Take $\Omega$ as the triangle $\vartriangle QQ_1Q_2$ with $Q(0,0)$, $Q_1(16,0)$ and $Q_2(16x_0, 16\sqrt{1-x_0^2})$ for some $|x_0|<1$. The convergence rates for different $\omega= \angle Q_1QQ_2 \in (\frac{\pi}{2}, \pi)$ determined by choosing different $x_0$ are shown in Table \ref{Ex3Q2Err}. Here, $R=\frac{24}{5}$, $\tau=\frac{1}{8}$ are used when $x_0=-0.8$, and default values are used for other cases. The results show that the convergence rate is not optimal when $\omega<\frac{2\pi}{3}$, and it is optimal when $\omega \in [\frac{2\pi}{3}, \pi)$. These results are consistent with the expected convergence rate $\mathcal{R}$ in Theorem \ref{thmuerr} or Table \ref{uh1errtab} for $\omega \in (\frac{\pi}{2}, \pi)$.

\begin{table}[!htbp]\tabcolsep0.03in
\caption{The $H^1$ error for $\omega\in (\frac{\pi}{2}, \pi)$ on quasi-uniform meshes.}
\begin{tabular}[c]{|c|c|c||c|c|c|c|}
\hline
parameter $x_0$ & $\omega$ & expected rate & $j=7$ & $j=8$ & {$j=9$} & {$j=10$}  \\
\hline
$-0.2$ & $\approx 0.56409\pi$ & $0.46$  & 0.75 & 0.67 & 0.59 & 0.54 \\
\hline
$-0.4$ & $\approx 0.63099\pi$ & $0.83$ &  0.96 & 0.95 & 0.94 & 0.93 \\
\hline
$-0.5$ & $\frac{2\pi}{3}$ & $1.00$  & 1.03 &  1.01 & 1.00 & 1.00 \\
\hline
$-0.6$ & $\approx 0.70483\pi$ & 1.00 &  1.01 & 1.01 & 1.01 & 1.00 \\
\hline
$-0.8$ & $\approx 0.79517\pi$ & 1.00 &  1.02 & 1.01 & 1.01 & 1.00 \\
\hline
\end{tabular}\label{Ex3Q2Err}
\end{table}

\noindent\textbf{Test case 2.}
We consider the polygon $\Omega$ with vertices $Q(0,0)$, $Q_1(16,0)$, $Q_2(-8, 8\sqrt{3})$, and $Q_3(-8, -8y_0\sqrt{3})$ for some $y_0\in (0,1]$, which gives $\omega=\angle Q_1QQ_3 \in (\pi, \frac{4\pi}{3}]$. We then consider the domain $\Omega$ (see Figure \ref{Ex1Q4}(a)) presented in Example \ref{P1S} Test case 2, and the corresponding angle $\omega\in (\frac{4\pi}{3}, \frac{3\pi}{2})$. The convergence rates for different $\omega \in (\pi,\frac{3\pi}{2})$ are shown in Table \ref{Ex3Q3Err}. The results show that the convergence rate is not optimal when $\omega<\frac{4\pi}{3}$, and it is optimal when $\omega \in [\frac{4\pi}{3}, \frac{3\pi}{2})$. These results are consistent with the expected convergence rate in Theorem \ref{thmuerr} or Table \ref{uh1errtab} for $\omega \in (\pi,\frac{3\pi}{2})$.

\begin{table}[!htbp]\tabcolsep0.03in
\caption{The $H^1$ error for $\omega\in (\pi,\frac{3\pi}{2})$ on quasi-uniform meshes.}
\begin{tabular}[c]{|c|c|c||c|c|c|c|}
\hline
parameter $y_0$ or domain & $\omega$ & expected rate & $j=6$ & $j=7$ & {$j=8$} & {$j=9$}  \\
\hline
$0.2$ & $\approx 1.10615\pi$ & $0.38$  & 0.82 & 0.72 & 0.62 & 0.53 \\
\hline
$0.6$ & $\approx 1.25612\pi$ & $0.82$ &  0.96 & 0.95 & 0.94 & 0.93 \\
\hline
$0.8$ & $\approx 1.30101\pi$ & $0.93$ &  0.98 & 0.98 & 0.98 & 0.98 \\
\hline
$1.0$ & $\frac{4\pi}{3}$ & 1.00  & 1.00 &  1.00 & 1.00 & 1.00 \\
\hline
$\Omega$ in Figure \ref{Ex1Q4}(a) & $\approx 1.38305\pi$ & 1.00 &  1.02 & 1.02 & 1.01 & 1.01 \\
\hline
\end{tabular}\label{Ex3Q3Err}
\end{table}

\noindent\textbf{Test case 3.}
We consider the polygon $\Omega$ with vertices $Q(0,0)$, $Q_1(\frac{16\sqrt{3}}{3},0)$, $Q_2(\frac{16-8\sqrt{2}}{1+\sqrt{3}}, \frac{16-8\sqrt{2}}{1+\sqrt{3}}+8\sqrt{2})$, $Q_3(-8\frac{\sqrt{2}+\frac{2\sqrt{3}}{2}}{1+\frac{1}{\sqrt{3}}}, 8\sqrt{2}-8\frac{\sqrt{2}+\frac{2\sqrt{3}}{2}}{1+\frac{1}{\sqrt{3}}})$, and $Q_4(x_1, -8\sqrt{3})$ for some $x_1 \in (0, 8\sqrt{3}]$, which generates $\omega =\angle Q_1QQ_4 \in (\frac{3\pi}{2}, \frac{7\pi}{4}]$.
The convergence rates for different $\omega \in (\frac{3\pi}{2},2\pi)$ are shown in Table \ref{Ex3Q4Err}. These results are consistent with the expected convergence rate in Theorem \ref{thmuerr} or Table \ref{uh1errtab} for $\omega \in (\frac{3\pi}{2},2\pi)$. 

\begin{table}[!htbp]\tabcolsep0.03in
\caption{The $H^1$ error for $\omega\in (\frac{3\pi}{2}, 2\pi)$ on quasi-uniform meshes.}
\begin{tabular}[c]{|c|c|c||c|c|c|c|}
\hline
$x_1$ or domain & $\omega$ & expected rate & $j=6$ & $j=7$ & {$j=8$} & {$j=9$}  \\
\hline
$x_1=4$ & $\approx 1.58946\pi$ & $0.23$  & 0.87 & 0.76 & 0.63 & 0.50 \\
\hline
$x_1=8$ & $\frac{5\pi}{3}$ & $0.40$ &  0.83 & 0.75 & 0.65 & 0.60 \\
\hline
$x_1=8\sqrt{3}$ & $\frac{7\pi}{4}$ & $0.57$  & 0.87 & 0.82 & 0.77 & 0.71 \\
\hline
\end{tabular}\label{Ex3Q4Err}
\end{table}

\end{example}

\section*{Availability of supporting data} 
Enquiries about data availability should be directed to the authors.

\section*{Competing interests} 
The authors declare that they have no conflict of interest.

\section*{Funding}
H. Li was supported in part by the National Science Foundation Grant DMS-2208321 and by the Wayne State University Faculty Competition for Postdoctoral Fellows Award. P. Yin was supported by the University of Texas at El Paso Startup Award.

\section*{Authors' contributions}
The authors contribute equally.

\section*{Acknowledgments}
The authors thank the anonymous referees who provided valuable comments resulting in improvements in this paper.

\bibliography{LYZ20}

\bibliographystyle{plain}

\end{document}